\documentclass[10pt,a4paper]{article}
\usepackage[utf8]{inputenc}
\usepackage[T1]{fontenc}
\usepackage{lmodern}
\usepackage{dsfont}
\usepackage[colorlinks=True]{hyperref}        
\usepackage[a4paper]{geometry}
\usepackage{authblk}
\usepackage{amsmath, amsfonts, amssymb}
\usepackage{amsbsy}
\numberwithin{equation}{section}
\usepackage{enumerate}
\usepackage{graphicx}
\usepackage{pdfrender,xcolor}
\usepackage{ntheorem}
\usepackage{sectsty}
\usepackage{pdfpages}
\usepackage[french,english]{babel}
\usepackage{mathrsfs}
\usepackage{microtype}
\usepackage{bm}
\usepackage{ulem}

\DeclareMathAlphabet{\mathpzc}{OT1}{pzc}{m}{it}
\DeclareMathAlphabet{\mathbfcal}{OT1}{cmsy}{b}{n}

\theoremstyle{break}
\newtheorem{theorem}{Theorem}
\newtheorem{corol}[theorem]{Corollary}
\newtheorem{lemme}[theorem]{Lemma}
\newtheorem{prop}[theorem]{Proposition}
\newtheorem{definit}[theorem]{Definition}
\newtheorem{rem}[theorem]{Remark}
\newtheorem{exemple}[theorem]{Example}

\newtheorem{hyp}[theorem]{Hypothesis}
\numberwithin{theorem}{section}
\def\Proof{\paragraph{Proof:}}

\allsectionsfont{\centering\mdseries\normalsize\scshape}

\newcommand{\HHRZ}[1]{}

\newcommand{\bs}[1]{\boldsymbol{#1}}

\renewcommand{\L}{\mathcal{L}}
\def\hL{\widehat{\L}}

\newcommand{\Lou}{\mathcal{L}_{\omega,u}}
\newcommand{\Lnou}{\mathcal{L}^{(n)}_{\omega,u}}
\newcommand{\B}{\mathcal{B}}
\def\U{U}
\newcommand{\M}{\mathcal{M}}
\newcommand{\calN}{\mathcal{N}}
\newcommand{\Z}{\mathbb{Z}}
\newcommand{\R}{\mathbb{R}}
\newcommand{\T}{\mathbb{T}}
\newcommand{\N}{\mathbb{N}}
\newcommand{\CX}{\mathbb{C}}
\newcommand{\F}{\mathcal{F}}
\def\O{\mathcal{O}}
\def\P{\mathcal{P}}
\newcommand{\J}{\mathcal{J}}
\newcommand{\C}{\mathcal{C}}

\newcommand{\calH}{\mathcal{H}}
\newcommand{\Id}{\mathds{1}}
\newcommand{\V}{\mathcal{V}}

\newcommand{\A}{\mathcal{A}}
\newcommand{\calA}{\mathcal{A}}

\def\Cveca{\C_{\vec{a}}}
\def\undr{{{\underline r}}}

\def\unds{{\,{\underline s}}}
\def\undt{{\,{\underline t}}}

\def\myD{n}
\def\locinv{\vartheta}
\def\tlocinv{\widehat{\vartheta}}
\def\teta{\widehat{\eta}}
\def\tpsi{\widehat{\psi}}

\def\whdel{\widehat{\delta}}

\def\tg{\widehat{g}}
\def\tT{\widehat{T}}

\def\hatg{\widehat{g}}
\DeclareMathOperator*{\supess}{ess\,sup}

\def\felem{f}
\def\field{{\mathbf f}}

\def\hrr{\hookrightarrow}
\def\lhrr{\hookrightarrow}
\def\nc{\newcommand}

\nc{\mylabel}[1]{\label{#1}}

\nc{\myremark}[1]{{\mbox{\textcolor{red}{\bf [#1]}}} }

\def\ds{\displaystyle}
\def\eps{\epsilon}
\def\DS{{\cal S}}
\def\Bj{\bs{j}}
\def\Bp{\bs{p}}
\def\BQ{\bs{Q}}
\def\Bz{\bs{z}}

\def\BC{\mathscr{C}}

\def\DE{\mathds{E}}
\def\BL{\mathscr{L}}
\def\hBL{\widehat{\mathscr{L}}}
\def\Bpi{\bs{\pi}}

\def\BLambda{\bs{\Lambda}}

\def\Bf{\bs{f}}
\def\Bhf{\bs{\widehat{f}}}
\def\Bphi{ {\bs{\varphi} }}
\def\Bhphi{\bs{\widehat{\phi}}}

\def\Bpsi{\bs{\psi}}

\def\rr{\rightarrow}

\def\Halmos{\hfill$\square$}
\def\diam{{\rm diam}}
\def\Aut{{\rm Aut}}
\def\bfone{{\bf 1}}
\def\calF{{\mathcal F}}

\def\Jd{\J}

\def\Card{{\rm Card}}
\def\myh{{h}}
\def\myp{{p}}

\def\hphi{\widehat{\phi}}

\title{Regularity of characteristic exponents and linear response
for transfer operator cocycles}
\author{Julien Sedro\footnote{Laboratoire de Probabilit\'es, Statistique et Mod\'elisation (LPSM), Sorbonne Universit\'e, Universit\'e de Paris, 4 Place Jussieu, 75005 Paris, France. \textit{E-mail}: sedro@lpsm.paris},~Hans Henrik Rugh\footnote{Institut de Mathématiques d'Orsay (IMO), Université Paris-Saclay, Site d'Orsay 91405 Orsay Cedex France. \textit{E-mail}: Hans-Henrik.Rugh@math.u-psud.fr }}

\date{\today}
\begin{document}
	\maketitle
\begin{abstract}
We consider cocycles obtained by composing sequences of
transfer operators with positive weights, associated with
uniformly expanding maps (possibly having countably many branches) and depending upon
parameters. Assuming $C^k$ regularity with respect to
coordinates and parameters, we show that when the sequence is picked within
a certain uniform family the top characteristic exponent and generator of top Oseledets space of the cocycle are $C^{k-1}$
in parameters. 
As applications, we obtain a 
linear response formula for the equivariant measure associated with random products of uniformly expanding maps,
and we study the regularity of the Hausdorff dimension of
a repeller associated with
random compositions of one-dimensional cookie-cutters.
\end{abstract}

\section{Introduction}

The Oseledets-Lyapunov spectrum for transfer operators cocycles associated with random composition of maps plays a key role in the study of its ergodic properties, similar in many aspects to the one played by the spectrum of the transfer operator in the study of statistical properties of autonomous dynamical systems:
this is the general philosophy arising from a series of papers generalizing Oseledets M.E.T to a cocycle of quasi-compact operators on a separable Banach space \cite{FLQ10,FLQ13,GTQ14}.
In a deterministic setting, a deep connection between stability of the statistical properties, such as statistical stability, linear and higher-order response, and stability of the transfer operator spectrum has been established: it then appears natural to wonder what remains of such a connection in the case of random composition of maps.  
\medskip

Indeed, the regularity with respect to parameters
of the top characteristic exponent of operator
cocycles has been thoroughly studied, mostly in analytical regularity or for random product of matrices:
Starting with Ruelle's seminal paper \cite{Ruelle82}, in which
he shows real-analyticity of the characteristic exponent for
a random product of positive matrices (within a certain
compact family contracting $\R_+^n$)
to Le Page \cite{LePage89}, establishing under less strict hypotheses,
Hölder and smooth regularity in the case of an i.i.d product of matrices, or Hennion's paper \cite{Hennion91}, which gives sufficient conditions
under which the characteristic exponent of an i.i.d product of matrices is differentiable.
In an infinite-dimensional setting, Dubois \cite{Dubois08} extended Ruelle's result to more general cone-contractions
and showed real-analyticity with
respect to parameters of
the top characteristic exponent for a sequence of 
real-analytical contracting operators. 
A further generalization to complex cone-contractions
was given in \cite{Rugh10}. One may also mention \cite{D17}, where genericity of analyticity of Lyapunov exponents for random bounded linear operators is shown.

For transfer operator cocycles, the literature on characteristic exponents stability is more recent: in the context of finite regularity, we may include \cite{FGTQ14,FGTQ15,DF18,GTQ18}, as well as the recent \cite{Crim19},
which present a generalization of the Keller-Liverani approach to discrete spectrum stability \cite{KL98} to the cocycle case. However, the results in \cite{Crim19} only gives continuity w.r.t parameters of the discrete spectrum, and do not give any explicit estimate on the modulus of continuity.  
\\In the analytical setting, an investigation on characteristic exponent regularity w.r.t parameters
was also recently used to establish a random analogue of
the Nagaev-Guivarc'h method (see \cite{Drag2018,DFGV19}).

\medskip

The problem we consider here concerns a family of cone-contracting transfer operators, specifically the problem of quenched (linear) response for a perturbed cocycle of expanding maps, i.e the regularity w.r.t parameters of the equivariant measure of this random dynamical system. This is a natural extension of the response problem for deterministic dynamics, which was investigated both for expanding maps (\cite{BCV16,Ba14,Ba16,BahsounSaussol16,BaladiSmania08,BaladiSmania12,BaladiSmaniaalt10,BaladiTodd16,GalatoloPollicott}) and more general hyperbolic systems (\cite{Ru97,Ru97erratum,Ru98,GL06,GL08,Dol04,Ba16}).
\\In the random setting, the annealed response (i.e. response for the invariant measure of the associated Markov chain) has also received extensive attention: besides the already mentioned \cite{GL06}, see the seminal \cite{HairerMajda10}  or the recent \cite{BRS17,GalatoloGiulietti17,GalatoloSedro20}. However, we emphasize that the above-mentioned
works only cover the annealed case,
and that 
to the authors knowledge, the present paper is
the first time that quenched response is addressed. 
It is certainly of major interest in many applications, notably for climate science: see \cite{CSG08,CSG11,LucariniL12}.
\\We note that even for a single deterministic system, in the context of finite regularity, the study of the
regularity of the physical measure upon a parameter is a considerable task, as the map associating the system to its transfer operator is 
not even continuous in the usual operator norm. The methods to overcome this problem may be roughly divided in three categories: first, Ruelle's original approach \cite{Ru97,Ru97erratum,Ru98}, 
relying on structural stability. This
allows for changes of variables that, although not smooth themselves, depend smoothly on parameters. Second, various forms of quantitative stability results for fixed points of linear Markov operators, as in \cite{BRS17,GalatoloPollicott, GalatoloGiulietti17,GalatoloSedro20}. Third, the so-called {\it weak spectral perturbation theory},
originally developed by
Keller and Liverani \cite{KL98}, who managed under conditions of
quasi-compactness 
to deduce H\"older continuity of the discrete part of
the spectrum. Gou\"ezel and Liverani \cite{GL06} then extended this
to differentiability and higher-order regularity of the discrete spectrum, also establishing explicit formulae for the $n$-th derivative of this spectral data.
However, this last approach relies heavily upon
the study of the resolvent operator, which makes it
unsuitable when considering non-autonomous dynamics, sequential or random compositions. 

\medskip

The approach we describe in this paper 
is suited both
to deterministic dynamics and to
random (non-autonomous) situations, including the quenched case. 
In the deterministic setting, our methods are easily adapted to the case of 
sequential dynamical systems (see Remark \ref{rem:sequential} below).
 This was considered, though not treated in full details,
 in Ruelle's seminal paper \cite{Ru97}. As for the random case, our results cover the case of a possibly non i.i.d random product chosen among an infinite, non-countable set, of uniformly expanding maps having an infinite number of branches, defined on any manifold that admits such maps. In the i.i.d case, we notice that our linear response result for the stationary measure of the skew-product associated to the random product (Theorem \ref{annealedlinrespthm}) 
 implies a linear response for the annealed invariant measure of a random composition of uniformly expanding maps, whence parts of the results\footnote{It should be noted that we only recover those results when one perturbs the maps themselves, which is not exactly the setup of \cite{BRS17}, where also perturbations of the probability measure on the base are considered (see also \cite[p.5, footnote 3]{BRS17}).} of \cite{BRS17}.
\\Our approach will rely on two main ingredients: Birkhoff's cone contraction theory to construct equivariant measures for random product of expanding maps\footnote{This step is quite classical, see e.g. \cite{BKS96}. For an occurrence of this idea in a context pertaining to linear response for non-autonomous systems, but somewhat different, see \cite{SelleyTanzi20}.} and, to study their dependence on parameters, a new idea, first appearing in
\cite{Sedro} in which a quantitative stability result for fixed point of a graded scheme 
of (possibly non-linear) operators was proven (we refer to the discussion at the beginning of §\ref{sec:loss-of-reg} for more on graded schemes of operators).
The main idea is that a fixed
point of a contracting $C^k$ map, when viewed in $C^{k-1}$ may gain
differentiability in parameters. We will apply this principle
to Banach spaces which will consist of 
$C^k$-sections of a fiber bundle  and show that the section 
becomes differentiable when viewed as a $C^{k-1}$ section: this implies
the wanted regularity of characteristic exponents.
Once one has become familiar with the general idea,
the actual implementation is quite simple, 
although somewhat hampered by a complicated
notation.
\\Our approach looks close in spirit to the 'quantitative stability for fixed points' family of results, but has two major differences: first, we may study fixed points of non-linear operators (notably fractional linear maps like in \eqref{eq:fractlinmap}). Second, we may study regularity of our fixed point beyond differentiability.
\medskip

The paper is organized as follows: in Section \ref{sec:cone-con}, we consider a cone-contracting cocycle of linear operators and construct the associated fractional linear map \eqref{eq:fractlinmap}. We construct its fixed point (Lemma \ref{thm:fixedptsection}), study some of its regularity properties (Lemma \ref{lem:Qop}), and relate it to the top characteristic exponent of the linear operator cocycle (Theorem \ref{thm:convenientformula}).
\\In Section \ref{sec:loss-of-reg}, we consider families of such cone-contracting cocycle of linear operators, depending upon a parameter and acting on a scale of Banach spaces. In Theorem \ref{thm:field} we show that, under a suitable regularity assumption for the cocycle, both the fixed point map (seen in the right Banach space) and top characteristic exponent are smooth. As a by-product, we get an implicit-function like formula for the derivative of the fixed point map.
\\In Section \ref{sec:dynapp}, we present an application
to the problem of linear
response for the equivariant measure associated to a random product of uniformly expanding maps. Here, issues of uniformity with respect to various parameters are paramount and we give a detailed construction of the $c^{(k,\alpha)}$-structure on a general Riemann manifold. In §\ref{subsec:cocyclemeasurability}, we construct a distance on a space of couples $(T,g)$ of uniformly expanding maps and weights, satisfying some uniformity and regularity conditions. This allows us to establish strong Bochner measurability of the transfer operator cocycle. Based on a "parameter-extraction" approach, we also establish the regularity of the transfer operator cocycle with respect to parameters. In §\ref{subsec:cones}, we exhibit a family of Birkhoff cones that are uniformly contracted by the (weighted) transfer operator of a $c^{r}$ uniformly expanding map. In §\ref{subsec:linresp} and §\ref{subsec:Hausdorffdimandstuff}, we establish a linear response formula for the equivariant (Theorem \ref{quenchedlinrespthm}) and stationary (Theorem \ref{annealedlinrespthm}) measures of a random product of smooth expanding maps.
As an example,
we show that the Hausdorff dimension of the repeller associated to one-dimensional cookie-cutters depends smoothly on the system.
\\Finally, in a series of appendices, we recall or prove various results useful to our purposes: in Appendix \ref{appendix:conecontraction}, we recall various results in real cone contraction theory. In Appendix \ref{appendix:measurability}, we have regrouped some material on (strong) Bochner measurability. In Appendix \ref{app:diffcalc}, we prove some results in classical differential calculus, instrumental in our approach; in particular, we state and prove several "regularity extraction" theorems, central to the main sections of this paper. Finally, in Appendix \ref{app:graded-diff-calc}, we establish a generalization of the classical Leibniz principle to families of operators having some "loss of regularity" property.\medskip

We end this introduction by giving two fairly simple but
non-trivial examples to which the results of the present paper apply.
Below, let $(\Omega,\mathcal{F},\mathds P)$ be a probability space and
let $\tau: \Omega\rr \Omega$ be an invertible
$\mathds P$-ergodic map (see however Remark \ref{rem:sequential}).

\begin{exemple}
	\label{ex:intro} 
	We equip $\R^d$ with the standard Euclidean norm
	and we let $\T^d=\R^d/\Z^d$ be the $d\geq 1$ torus with the
	induced Riemannian metric. For $\theta_0<1$ we consider
	the subset $G(\theta_0)\subset M_d(\Z)$
	 of invertible integer matrices $M$ such that
	the  $\|M^{-1}\|_{L(\R^d;\R^d)}\leq \theta_0$.
	Finally let $0<\theta_1 < \frac{1}{\theta_0}-1$.

	Consider measurable maps
	$\omega\in \Omega \mapsto   M_\omega \in G(\theta_0)$
	and\footnote{By $c^r$ we mean the closure of $C^\infty$
	functions in a space of $C^r$ functions,
	cf.\ Definition  \ref{def:cka}}
	$\omega \in\Omega \mapsto \kappa_\omega
	\in c^{r} (\R \times \T^d; \R)$, $r>2$
	 for which
	\ $\sup_{(\omega,u,x)}\|
	   \partial_x \kappa_\omega(u,x)\|_{L(\R^d;\R)} \leq \theta_1$ \
	and \ $\sup_\omega \|\kappa_\omega\|_{C^r} \leq K<+\infty$.
	We define the parametrized random expanding map:
	$T_{\omega,u}:\T^d\to\T^d$, $\omega\in\Omega$, $u\in \R$ by
	\begin{equation*}
	  %\label{elementaryrandex}
	T_{\omega,u}(x):=M_\omega x+\kappa_\omega(u,x)\mod  \Z^d
	\end{equation*}
	%We consider the composition of random maps
	%\eqref{elementaryrandex} along orbits of $\tau$:
	%\begin{equation*} 
	%T^{(n)}_{\omega,u}:=
	%T_{\tau^{n-1}\omega,u} \circ \cdots \circ T_{\omega,u}.
	%\end{equation*}
\end{exemple}

%\begin{exemple}
	%In one dimension we obtain the following somewhat simpler special case:
	%Let $\omega\in \Omega \mapsto d_\omega \in \Z\setminus\{-1,0,1\}$
	%and $\omega \in\Omega \mapsto \kappa_\omega \in C^{r} (\R \times \S^1; \R)$, $r>2$
	%%$\omega\in \Omega \mapsto a_\omega\in S^1=\R/\Z$
	%be Bochner measurable maps for which
	%\ $\sup_{(\omega,u,x)}|\partial_x \kappa_\omega(u,x)| \leq \theta<1$ \
	%and \ $\sup_\omega \|\kappa_\omega\|_{C^r} \leq K<+\infty$.
	%We define the parametrized random expanding map:
	%$T_{\omega,u}:\S^1\to\S^1$, $\omega\in\Omega$, $u\in \R$ by
	%\begin{equation}\label{elementaryrandex}
	%T_{\omega,u}(x):=d_\omega x+\kappa_\omega(u,x)\mod  \Z
	%\end{equation}
	%We consider the composition of random maps
	%\eqref{elementaryrandex} along orbits of $\tau$:
	%\begin{equation*} 
	%T^{(n)}_{\omega,u}:=
	%T_{\tau^{n-1}\omega,u} \circ \cdots \circ T_{\omega,u}.
	%%\end{equation*}
%\end{exemple}

\begin{exemple}
	\label{ex:introtwo} 
	Let $\omega \mapsto \kappa_\omega \in c^{r} (\R \times [0,1]; [0,1])$, $r>2$
	%$\omega\in \Omega \mapsto a_\omega\in S^1=\R/\Z$
	be a measurable map for which 
	$\sup_\omega\|\kappa_\omega\|_{C^r}  \leq K<+\infty$,
	$(x\mapsto \kappa_\omega(u,x)) \in {\rm Diff}^r_+([0,1])$ and
	\ $\inf_x \partial_x \kappa_\omega(u,x) \geq \theta > \frac12$
	for all
	$u\in \R$ and $\omega\in\Omega$.
	We define in this case:
	$T_{\omega,u}:D_{\omega,u} \to(0,1)$, $\omega\in\Omega$, $u\in \R$ by
	\begin{equation*}
	   %\label{elementaryrandex2}
	T_{\omega,u}(x):=
	\frac{1}{\kappa_\omega(u,x)}
	- \left\lfloor\frac{1}{\kappa_\omega(u,x)} \right\rfloor,
	\end{equation*}
	where the domain of definition is
	$D_{\omega,u}=\{ x\in (0,1): 
	1/\kappa_\omega(u,x)\notin \N \}$ (which is co-countable for fixed $\omega,u$).
	Note that when $\kappa_\omega(u,x)\equiv x$, $T_{\omega,u}$ is the standard Gauss map.
\end{exemple}
For either of the two
examples above we consider the composition of random maps 
%\eqref{elementaryrandex}
along orbits of $\tau$:
\begin{equation*} 
T^{(n)}_{\omega,u}:=
T_{\tau^{n-1}\omega,u} \circ \cdots \circ T_{\omega,u}.
\end{equation*}
In the first  example this is defined for all $x\in \T^d$ 
 and in the
second example on a set of full Lebesgue measure in $(0,1)$.

\begin{theorem}
For the above two examples we have with
$M=\T^d$ and  $M=(0,1)$, respectively:
	\label{exquenchedlinresponse}
	\begin{itemize}
		\item For every $u\in \R$, the \emph{skew-product} 
		$F_u(\omega,x):=(\tau\omega,T_{\omega,u}(x))$ admits a unique
		stationary probability measure $\nu_u$ on $\Omega\times M$,
		whose decomposition along
		the marginal $\mathds{P}$, $\nu_{\omega,u}$, 
		has Lebesgue-density
		$h_{\omega,u}\in C^{r-1}(M)$. 
		Furthermore,
		$\supess_{\omega\in\Omega}\|h_{\omega,u}\|_{C^{r-1}(M)}<+\infty$.
		\item For any $1\le s<r-2$, the map $u\in \R\mapsto h_u\in L^\infty(\Omega,C^{s}(M))$ 
		is differentiable, and one has the 
		following \emph{linear response formulae} for any
		$\psi\in L^1(M)$:
		%\begin{align}
		%\partial_u\left[\int_{\mathbb{S}^1}\phi h_{\omega,u}dm\right]_{u=0}&=
		%\sum_{n=0}^\infty\int_{\mathbb{S}^1}
		%\phi\circ T^{(n)}_{\tau^{-n}\omega,0}
		%\partial_u\L_{\tau^{-(n+1)}\omega,0}h_{\tau^{-(n+1)}\omega,0}dm\\
		%\partial_u\left[\int_{\Omega}\int_{\mathbb{S}^1}
		%\phi h_{\omega,u}dmd\mathds{P} \right]_{u=0}&=
		%\sum_{n=0}^\infty\int_\Omega\int_{\mathbb{S}^1}
		%\phi\circ T^{(n)}_{\omega,0}
		%\partial_u\L_{\tau^{-1}\omega,0}h_{\tau^{-1}\omega,0}dmd\mathds{P}
		%\end{align}
		\begin{align}
		\partial_u\left[\int_{M}\psi h_{\omega,u}dm\right]_{u=u_0}&=
		\sum_{n=0}^\infty\int_{M}
		\psi\circ T^{(n)}_{\tau^{-n}\omega,u_0}
		\;\partial_u\L_{\tau^{-(n+1)}\omega,u_0}h_{\tau^{-(n+1)}\omega,u_0}dm
	\label{first-quenched}	
		\\
		\partial_u\left[\int_{\Omega}\int_{M}
		\psi h_{\omega,u}dmd\mathds{P} \right]_{u=u_0}&=
		\sum_{n=0}^\infty\int_\Omega\int_{M}
		\psi\circ T^{(n)}_{\omega,u_0}
		\;\partial_u\L_{\tau^{-1}\omega,u_0}h_{\tau^{-1}\omega,u_0}dm\,d\mathds{P}
	\label{first-annealed}	
		\end{align}
	\end{itemize}
\end{theorem}

Proofs for the above two examples are given in Section \ref{sub:examples}.

\begin{rem} \label{rem:sequential} We have chosen to focus upon the
  linear response for an invertible ergodic random dynamical system. 
  However, the method presented here
  allows for various other possibilities. We mention
  the following:
  \begin{itemize}
  \item
  In the article, by "quenched linear response"
  we mean a linear response when the random environment is fixed.
  This requires in the  context
  of Formula \eqref{first-quenched},
  that $\tau$ is invertible since our construction of $h$ is 
  computed along a (unique) backward orbit under $\tau$.
   It does not, however, rely upon
  the existence of a $\tau$-invariant probability measure on $\Omega$.
  Thus, modulo some obvious modifications, the formula
   also holds when  taking
  $\Omega=\Z$ and
  letting $\tau$ be the shift. This corresponds to what is known as
  a sequential dynamical system.
  \item  
  For our annealed formula, \eqref{first-annealed}, measurability 
  of $\tau$ and
  invariance of the measure $\mathds{P}$ is necessary.
    In order to have unique characteristic exponents,
    ergodicity of $\tau$
    is usually important but in the present situation it is not.
   This boils down to the characteristic
  exponent being identically zero in the above examples.
  % the ergodic map $\tau$ need not be invertible (as in Oseledets original
  % work).
  \end{itemize}
\end{rem}

\section{A cone-contracting operator cocycle}
\label{sec:cone-con}

We will here give an abstract formulation
of the above type of problems. We consider a Banach space $(E,\|.\|)$ and 
a collection $M\subset L(E)$
of bounded linear operators which contracts uniformly
a regular convex cone $\C\subset E$. We
make the following assumptions  on $\C$ and $M$
(see Appendix \ref{appendix:conecontraction} for further details):
\begin{enumerate}
	\item[(H1)] (outer regularity) There is $\ell\in E'$ of norm one 
	and $K\in[1,+\infty)$ so that for every $\phi\in \C$:
	\begin{equation*}
	\langle \ell,\phi \rangle \geq \frac{1}{K} \|\phi\|
	\label{outer reg}
	\end{equation*}
	\item[(H2)] (inner regularity) There is $\rho\in(0,1]$ so that the sub-cone
	\begin{equation*}
	\C(\rho)= \{ \phi\in \C : B(\phi,\rho\|\phi\|) \subset \C \}
	\label{inner reg}
	\end{equation*}
	is non-empty (in particular, $\C$ has non-empty interior). We fix such a $\rho$ in the following.
	We define:
	\begin{equation*}
	\C_{\ell=1} = \{ \phi \in \C : \langle \ell,\phi \rangle = 1\} 
	\ \ \mbox{and} \ \
	\C_{\ell=1}(\rho) = \{ \phi \in \C(\rho) : 
	\langle \ell,\phi \rangle = 1\} .
	\end{equation*}
	\item[(H3)] (uniform contraction) We require that for every $L\in M$:
	\begin{equation*} 
	L(\C^*)\subset \C^*(\rho),
	\end{equation*}
	Here, we have written
        $\C^*:=\C-\{0\}$ and $\C^*(\rho):=\C(\rho)-\{0\}$ for the punctured cones. 
	\item[(H4)] (uniform bounds) 
	There is $1\leq \vartheta<+\infty$ so that for every $L\in M$
	\begin{equation*}
	\frac{1}{\vartheta} \leq  \|L\| \leq \vartheta.
	%\label{eq:unifbounds}
	\end{equation*}
\end{enumerate}

Such operators enjoy strong contraction properties for the
projective Hilbert metric. In the following we will give a short
summary of the consequences which will be needed below.
As this is fairly standard we will be quite brief and refer
to e.g. \cite{Ba00,Rugh10} for standard properties
of cones and the above mentioned notions
of regularity and to \cite{Dubois08,Ruelle82,Rugh08} 
for more details on the standard construction of the
cone-contracting cocycles given below.

First, $\C(\rho)$ has finite
diameter in $\C$. From Corollary \ref{corol dC bound} one gets:
\[\ds \Delta = {\rm diam}_{\C} (\C(\rho)^*)
\leq 2\log \left( 1+\frac{2 K}{\rho} \right).\]
This yields a strict Birkhoff contraction factor:
$\eta =\tanh\frac{\Delta}{4}\leq \frac{K}{K+\rho}<1$.
Also, the above four conditions imply
(see Lemma \ref{lemma:upperlower}) that for any $L\in M$,
$\phi\in \C(\rho)$:
\begin{equation}
\frac{\rho}{\vartheta K} 
\; \|\phi\| \ \leq  \
\langle \ell, L \phi \rangle \  \leq  \ \vartheta \; \|\phi\| 
\ \ \mbox{and} \ \
\frac{\rho}{\vartheta K} 
\; \langle \ell,\phi \rangle \ \leq  \
\langle \ell, L \phi \rangle \  \leq  
\ \vartheta \; K  \langle \ell, \phi \rangle .
\label{L equiv}
\end{equation}
Let $(\Omega,\mathcal{F},\mathds{P})$ be a probability space and
$\tau:\Omega\to\Omega$ a measure-preserving, invertible, and ergodic map.
We will define an operator valued cocycle over $\Omega$ in the following
way: let $X =L^\infty(\Omega,E)$ denote the
set of uniformly bounded $\calF$-Bochner measurable sections
$\Bphi=(\phi_\omega)_{\omega\in \Omega} : \Omega\rr E$ of
the product bundle
$\Omega\times E$.
We write 
$\|\Bphi\|_X=\sup_{\omega\in\Omega} \|\phi_\omega\|<\infty$ for the norm:
$(X,\|\cdot\|_{X})$ is then a Banach space.
We will also let $\BC=\C(\Omega)\subset X$ denote the 
bounded $\calF$-Bochner measurable sections taking values in the cone $\C$. 
We define in a similar way \[\BC(\rho):=\{\Bphi\in X,~\forall\omega\in\Omega,~\phi_\omega\in\C(\rho)\},\] 
\[\BC_{\ell=1}:=\{\Bphi\in X,~\forall\omega\in\Omega,~\phi_\omega\in\C_{\ell=1}(\rho)\},\] 
\[\BC_{\ell=1}(\rho):=\{\Bphi\in X,~\forall\omega\in\Omega,~\phi_\omega\in\C_{\ell=1}(\rho)\}.\]
We write $\phi\in\BC^*$ 
for a section that is nowhere vanishing, 
	i.e.\ $\phi_\omega \in \C^*(\rho)$ for every $\omega\in \Omega$.
\\

We consider a family of operators of the form
$\BL = (\L_\omega)_{\omega\in\Omega} : \Omega \rr M$, i.e.
each $\L_\omega\in M$ is a uniform cone-contraction
satisfying (H3) and (H4).
We let this family act as a bundle map 
in the following way ($\Bphi\in X$ and $\omega\in \Omega$):
\begin{equation*}
\BL : X \rr X , \ \ \
(\BL \Bphi)_\omega = \L_{\tau^{-1} \omega} \phi_{\tau^{-1} \omega} . 
\end{equation*}
We will further assume that $\BL$ is strongly Bochner measurable,
i.e. it preserves Bochner
measurability of sections
(cf. Appendix \ref{appendix:measurability}).
The $n$-th iterate is given by:
$$(\BL^n \Bphi)_\omega = \L_{\tau^{-1}\omega} 
\circ \cdots \circ \L_{\tau^{-n} \omega} \phi_{\tau^{-n}\omega}.$$
For $\ell \in E'$ satisfying (H1), we define the section of normalization factors:
\begin{equation*}
\Lambda : X \rr L^\infty(\Omega;\R_+), \ \ \
(\Lambda(\Bphi))_\omega  =
\langle \ell, \L_{\tau^{-1}\omega} \phi_{\tau^{-1}\omega} \rangle.
\end{equation*}
and the fractional linear map $\ds \Bpi : \BC^* \rr \BC_{\ell=1}(\rho)$, defined for $\omega\in\Omega$ by
\begin{equation}
(\Bpi (\Bphi))_\omega = 
\frac{\L_{\tau^{-1} \omega} \phi_{\tau^{-1} \omega}}
{\langle \ell, \L_{\tau^{-1}\omega} \phi_{\tau^{-1} \omega}\rangle}.
\label{eq:fractlinmap}
\end{equation}
Here, as well as below, we will adopt the convention that
a complex or real valued section over $\Omega$
acts upon $X$ by fiberwise multiplication. Similarly,
a nowhere vanishing section has an inverse given by taking fiberwise
inverse. The map \eqref{eq:fractlinmap} may then be written 
in the following compact way:
\begin{equation*}
\Bpi (\Bphi) = \Lambda(\Bphi)^{-1} \BL(\Bphi).
\end{equation*}
We may then establish the following:
\begin{lemme}\label{thm:fixedptsection}
	The map \eqref{eq:fractlinmap} admits a unique fixed point section
	$\field=\Bpi\circ \field \in \BC_{\ell=1}(\rho)$.
\end{lemme}
\Proof We see from Corollary \ref{CORA8}  that
for $m\geq n\geq 1$, $\Bphi,\Bpsi\in \BC^*$:
\begin{equation*}
\left\| \Bpi^{n}(\Bphi) - \Bpi^{m}(\Bpsi)\right\|_X
\leq \frac{K}{2} \eta^{n-1} \Delta.
\label{pi contraction}
\end{equation*}
In particular, $(\Bpi^n (\Bphi))_n$ is Cauchy in 
$\BC_{\ell=1}(\rho) \subset X$, whence converges to some  
$\field=(\felem_\omega)_{\omega\in \Omega} \in \BC_{\ell=1}(\rho)$,
which is a fixed point
of $\Bpi$. \Halmos
\medskip

Let us set $\Bp=\Lambda(\field) = \langle \ell, \BL \field\rangle : \Omega\rr \R^*_+$, 
which by \eqref{L equiv} verifies
\begin{equation*}
	p_\omega = \langle \ell, \L_{\tau^{-1}\omega}
	\felem_{\tau^{-1}\omega} \rangle \in 
	\left[ \frac{\rho}{\vartheta K}, \vartheta K\right].
	\label{norm factor}
\end{equation*}
For fixed $\omega\in \Omega$ we define the characteristic exponent
of the cocycle operator as follows:
\begin{equation}
	\chi_\omega =
	\limsup_n \frac{1}{n} \log \| (\BL^n)_\omega\|
	\label{char exp}
\end{equation}
We then have:
\begin{theorem}\label{thm:convenientformula}
	Let $\BL=(\L_\omega)_{\omega\in \Omega}$ be a 
	strongly Bochner measurable
	cocycle consisting of 
	uniform cone-contractions
	as defined above with
	$\ds \int_\Omega \left| \log \left\|\L_\omega\right\| \right| d\mathds P(\omega)
	<+\infty$.  Then
	\begin{enumerate}
		\item 
		The real valued section
		$\Bp =  (p_\omega)_{\omega\in \Omega} =
		\langle \ell, \BL \field\rangle : \Omega\rr \R^*_+$, 
		is measurable
		and $\log \Bp$ is $\mathds P$-integrable.
		\item The characteristic exponent \eqref{char exp}
		of the cocycle equals a.s.
		\begin{equation}\label{eq:convenientformula}
			\chi(\BL) = \DE\left( \log \Bp \right) =  \int_\Omega \log p_\omega d \mathds P(\omega) ,
		\end{equation}
	\end{enumerate}
\end{theorem}
\Proof
The first item is an easy consequence of \eqref{L equiv}. 
\\For the second item, note that we have the identity $\L_{\tau^{-1}\omega} \felem_{\tau^{-1}\omega} =
p_\omega \felem_{ \omega}$.
By the uniform bounds in \eqref{L equiv}, we see that \eqref{char exp} is equivalent to:
\begin{equation*}
	\chi_\omega=
	\limsup_{n\to\infty} \frac{1}{n} \log   \langle \ell, (\BL^n \field)_\omega \rangle =
	\limsup_{n\to\infty} \frac{1}{n} \sum_{k=0}^{n-1} \log p_{\tau^{-k} \omega}.
\end{equation*}
We are thus reduced to consider Birkhoff averages of $\log p_\omega$.
Since $\log \bs{p}\in L^1(\mathds P)$, we conclude by Birkhoff's theorem that a.e. the limsup is 
in fact a limit and (by ergodicity)
a constant equal to the integral of $\log \bs{p}$.\Halmos
\medskip

By \ref{CORA5} and \ref{CORA8} we see that
for $\phi\in \BC_{\ell=1}(\rho)$ and $z\in X$ small enough:
\begin{equation}
\left\| \Bpi^{n}(\Bphi) - \Bpi^{n}(\Bphi+\Bz)\right\|_X
\leq \eta^{n-1} \frac{K}{\rho } \left( \|\Bz\| + o(\|\Bz\|)\right) .
\label{pi n devel}
\end{equation}
This entails the following:

\begin{lemme}
	\label{lem:Qop}
	There is $C=C(\rho,K,\vartheta)>0$
	so that for $\Bphi\in \BC_{\ell=1}(\rho)$ and $\Bz\in X$
	with $\|\Bz\|< \frac{\rho}{2\vartheta^2 K}$:
	\begin{equation}
	\| \Bpi(\Bphi+z) - \Bpi(\Bphi) - \BQ(\Bphi) . \Bz \|
	\leq C  \| \Bz \|^2,
	\label{expand:Bpi}
	\end{equation}
	Here, $\BQ(\Bphi)$ 
	is the derivative of
	$\Bpi(\Bphi)$ 
	which is given by the expression:
	\begin{equation}
	\BQ(\Bphi).\Bz :=  \partial_{\Bphi} \Bpi(\Bphi). \Bz =
	\Lambda(\Bphi)^{-1} \left( \BL \Bz - \Lambda(\Bz) \Bpi(\Bphi)\right).
	\label{piderivative}
	\end{equation}
    At the fixed point $\field\in\BC_{\ell=1}(\rho)$,
the operator $\bfone-\BQ(\field)$ is invertible
	%\HHR{Isn't it true only at the fixed point $\field$ ?}
	and satisfies the bound:
	\begin{equation*}
	\left \|\left( \bfone- \BQ(\field) \right)^{-1} \right\|
	\leq 1 + \frac{K/\rho}{1-\eta}.
	%\label{bound:oneeta}
	\end{equation*}
\end{lemme}
\Proof
In the following consider $\Bphi$ and $\Bz$ as in the Lemma.
From \eqref{L equiv} we deduce that
$\Lambda(\Bphi): \Omega \rr 
\left[ \frac{\rho}{\vartheta K},\vartheta K \right]\subset \R_+^*$,
i.e.\ are uniformly bounded from above and below.
\\When $\|\Bz\| < \frac{\rho}{2 \vartheta^2 K}$, we have
$|\Lambda(\Bz)| < \frac{\rho}{2 \vartheta K}\leq \frac12 \Lambda(\Bphi)$
implying that
\begin{equation}\label{eq:boundnormfactor}
\Lambda(\Bphi +\Bz)=\Lambda(\Bphi )+\Lambda(\Bz) \geq \frac{\rho}{2\vartheta K}.
\end{equation}
Thus, 
$\Lambda(\Bphi +\Bz)$ is invertible
in the Banach algebra $L^\infty(\Omega;\R)$ through a Neumann series.
and $\Bpi(\Bphi+\Bz)$ is analytic in $\Bz$ and uniformly
bounded by $2\vartheta^2K/\rho$.
One has the identity:
\begin{align*}
\Bpi(\Bphi+\Bz)- \Bpi(\Bphi) &=
\Lambda(\Bphi+\Bz)^{-1} \left( \BL \Bz - \Lambda(\Bz) \Bpi(\Bphi) \right) \\
&=  (\bfone- \Lambda(\Bphi+\Bz)^{-1} \Lambda(\Bz)) Q(\Bphi).\Bz,
\end{align*}
with $Q(\Bphi)$ as in \eqref{piderivative}.
This implies, by \eqref{L equiv}, \eqref{pi n devel} and \eqref{eq:boundnormfactor} 
\begin{align*}
\|\Bpi(\Bphi+\Bz)- \Bpi(\Bphi)- Q(\Bphi).z\|  &\leq 
\| \Lambda(\Bphi+\Bz)^{-1} \; \Lambda(\Bz) \; Q(\Bphi) \; \Bz\|\\
& \leq  \frac{2\vartheta K}{\rho}\vartheta\|\Bz\|\frac{2K}{\rho} \; \| \Bz\| =
\frac{4 \vartheta^2 K^4}{\rho^2} \|\Bz\|^2.
\end{align*}
By \eqref{pi n devel} we get, at the fixed point $\field$, the uniform bound
$\|(\BQ(\field))^n\|\leq \frac{K}{\rho}\eta^{n-1}$, $n\geq 1$
so by a Neumann
series we obtain invertibility of $1-\BQ(\field)$ 
as well as the stated bound. \Halmos \\

Our goal in the following will be
to consider a family of measurable cocycles
of uniform cone-contractions depending 
on a parameter 
$u$ (in a Banach space), and see how the regularity
of the characteristic exponent and also the
fixed point section w.r.t. $u$ depends upon the regularity of
the cocycle. We emphasize that since the work of Ruelle \cite{Ruelle82}
on analytic matrix co-cycles,
this is well-understood 
also when the cocycle is a family of transfer
operators depending in an analytic way upon $u$, see e.g.
\cite{Dubois08,Ruelle82,Rugh08,Rugh10}.
However, for transfer operators associated with compositions
of maps with finite regularity, results are sparse and 
very incomplete due to an inherent 'loss of regularity'.
It seems appropriate to first treat a slightly more abstract setup and
then apply it to a cocycle of transfer operators. The readers who wish
to have a concrete example in mind may want to consult Section \ref{sec:dynapp} first.

\section{Cone-contracting cocycles with  loss of regularity}\label{sec:loss-of-reg}

Let $\B$ be a Banach space and $\U \subset\B$ an open and convex
subset. An element $u\in\U$ will be considered as a parameter for our problem.
Let $r_0>0$ and  consider a family, indexed by $r\in (0,r_0]$
of Banach spaces with associated
regular cones as described in the previous section, i.e. 
\[ \DS_r=(E_r,\C_r,\rho_r,K_r,\vartheta_r,\ell_r), \ \  r\in (0,r_0]  .\]
We assume that the family is 'graded' in the sense that there are
continuous linear  injections 
$j_{s,r}: E_{r} \hrr E_{s}$, 
for every $0<s<r\leq r_0$, satisfying uniform bounds with
respect to $s$ and $r$, i.e there is some constant $C>0$ independent from $r$ and $r$ such that $\|j_{s,r}\|_{E_r\to E_s}\leq C$, and verifying a
transitivity condition:  
$\forall 0<t<s<r\leq r_0$: $j_{t,s}j_{s,r} = j_{t,r}$.
\\As a natural example, the reader may think of $E_r$ as being
$C^r$-functions on a manifold, the 'downgrading' as being
the natural injection from $C^r$ into $C^s$ for $0<s<r\leq r_0$, 
and $\C_r$ to be a cone of nonnegative $C^r$ functions of the form described 
in subsection \ref{subsec:cones}.
%This will precisely be the case for our application in section \ref{sec:dynapp}.
% When the source space is clear from the context (or irrelevant)
%we sometimes write 
%$j_{s}: E_{r} \hrr E_{s}$ for this injection.
\\For simplicity, we will assume that the injections preserve cones and the linear forms
above, i.e.  $j_{s,r}(C_{r})\subset C_s$ and 
$\ell_r =   \ell_{s}\circ j_{s,r}$.
In this construction 
we also assume that all constants and norms are uniformly bounded on compact
subsets of $(0,r_0]$ (in our applications some constants
diverge as $r\rr 0^+$).\\

Given a probability space
$(\Omega,\mathcal{F},\mathds{P})$  as above we 
construct Bochner-measurable sections for each system $\DS_r$.
In this way, we obtain spaces $X_r=L^\infty(\Omega; E_r)$ 
as well as bounded cone fields 
$\BC_r \equiv  \C_r(\Omega)=\{\Bphi : \Omega \rr C_r
\mbox{\ Bochner measurable} \}$
and the associated slices 
$\BC_{r,\ell_r=1} \equiv \C_{r,\ell_r=1}(\Omega)$ etc. as above.
The injection $j_{s,r}$ induces  by pointwise action maps on
Bochner-measurable sections:
$\Bj_{s,r}: \BC_{r}=\C_{r}(\Omega) \hrr \BC_{s}=\C_{s}(\Omega)$
and $\Bj_{s,r} : X_{r} \hrr X_{s}$.
This family will again be transitive.

For each system $\DS_r$, $r\in (0,r_0]$ we assume given 
a family of cone-contracting cocycles,
$\BL_{r,u}$ acting upon $X_r$ and depending upon the
parameter $u\in U$.
We assume that the family
commutes with our injections in the natural way:
($0< s<r\leq r_0$): \ \ 
$\Bj_{s,r} \circ \BL_{r,u} = 
\BL_{s,u} \circ \Bj_{s,r}$. When written out over fibers:
\begin{equation*}
j_{s,r} \L_{r,u,\tau^{-1}\omega} \phi_{r,\tau^{-1} \omega} =
\L_{s,u,\tau^{-1}\omega} j_{s,r} \phi_{r,\tau^{-1}\omega}, 
\ \ \ u\in U,
\omega\in \Omega,
\end{equation*}
when acting upon $\Bphi_r\in X_r$.
We will often be using the normalizing fields obtained by acting with
$\ell_r$:
\begin{equation*}
\Lambda_{r,u,\omega}(\phi_{r,\omega}) := 
{\langle \ell_{r}, \L_{r,u,\omega} 
	(\phi_{r,\omega}) \rangle }  \in \R,  \ \ \omega\in \Omega.
\end{equation*}
Omitting mentionning $\omega$ explicitly, we write:
\begin{equation*}
\BLambda_{r,u}(\Bphi_r) := 
{\langle \ell_{r}, \BL_{r,u} 
	(\Bphi_{r}) \rangle }  \in L^{\infty}(\Omega).
\end{equation*}

When $\Bphi_{r}\in \BC_{r,\ell_r=1}$ the normalizing field
$\BLambda_{r,u}(\Bphi_r)$ takes values in 
$\left[ \frac{\rho_s}{\vartheta K}, \vartheta K\right]$, so
is invertible in the Banach algebra $L^\infty(\Omega)$.
Therefore, one has
\[ \Bpi_{s,u}(\Bj_{s,r} \Bphi_{r})_{\tau \omega} = 
\frac
{\L_{s,u,\omega} (j_{s,r} \phi_{r,\omega}) }
{\langle \ell_s, \L_{s,u,\omega}
	(j_{s,r} \phi_{r,\omega} \rangle }
=\frac
{j_{s,r} (\L_{r,u,\omega} \phi_{r,\omega}) }
{\langle \ell_{r}, \L_{r,u,\omega} 
	(\phi_{r,\omega}) \rangle }
= j_{s,r} \Bpi_{r,u}(\Bphi_{r})_{\tau \omega}.
\]
Thus, we also have $\Bpi_{s,u} \circ \Bj_{s,r} = \Bj_{s,r} \circ \Bpi_{r,u}$.

A family of fields, $\Bphi = (\Bphi_r)_{r\in (0,r_0]}$ with
each $\Bphi_r\in X_r$, is said to be {\em consistent} if 
for every $0<s<r\leq r_0$, $\Bphi_s = \Bj_{s,r} \Bphi_r$.
Note, that given a consistent family one has for every $0<s<r\leq r_0$: 
\begin{equation*}
\BLambda_{r,u}(\Bphi_r(u))=
\langle \ell_s, \Bj_{s,r} \BL_{r,u} \Bphi_r(u) \rangle =
\langle \ell_s, \BL_{s,u} \Bj_{s,r} \Bphi_r(u) \rangle =
\BLambda_{s,u}(\Bphi_s(u)).
\end{equation*}
We sometimes write 
$\BLambda_u(\Bphi)$ for the common value in this case.
%In the following, given a consistent family, we sometimes write
%$\BL_{s,u} \Bphi$ instead of $\BL_{s,u} \Bphi_s$ and
%$\Bpi_{s,u}(\Bphi)$ (when defined) instead of 
%$\Bpi_{s,u}(\Bphi_s)$. The index $s$ is inferred from the context.

\begin{lemme}
	\label{lem:fixedfield}
	There exists a unique parametrized consistent family 
	$\field(u) = (\field_r(u))_{r\in (0,r_0]}$, $u\in U$ 
	for which  
	\begin{equation*}
	\field_r (u) = \Bpi_{r,u}(\field_r(u)) \in \BC_{r,\ell_r=1}, \ \ \
	r\in (0,r_0],  \ \ u\in U.
	\end{equation*}
\end{lemme}
\Proof Pick a
section $\Bphi_{r_0} \in \BC_{r_0,\ell_{r_0}=1}$ and set
$\Bphi_s= j_{s,r_0} \Bphi_{r_0}\in \BC_{s,\ell_s=1}$, $0<s \leq r_0$
and iterate 
as in the previous section in each member of the family. 
We obtain fixed fields,
$\field_r(u) = \Bpi_{r,u}(\field_r(u))$
and since $\Bj_{s,r} (\Bpi_{r,u})^n (\Bphi_{r,u})= (\Bpi_{s,u})^n
(\Bj_{s,r}\Bphi_{r,u})$ each iterate constitutes a consistent 
family and this relationship 
holds while taking the limits. Note, however, that the convergence
may not be uniform over the scale $(0,r_0]$.
\Halmos
\\

Up to now we haven't assumed anything about
the regularity of
$u\mapsto \BL_{r,u} \in L(X_r)$, whence of the fixed fields
at the individual levels, $r\in (0,r_0]$.
In fact in the applications below, 
there is a priori no regularity of $u \mapsto \BL_{r,u} \in L(X_r)$
w.r.t. $u$, apart from being uniformly bounded.
The idea in the following is that by 'downgrading' the target space
within the family we recover regularity of the map 
$u\mapsto \Bj_{s,r} \BL_{r,u} = \BL_{s,u} \Bj_{{s,r}} \in L(X_{r},X_s)$,
$0<s<r\leq r_0$
which in turn we transform into
regularity of the fixed field section 
$u \mapsto \field_s(u)\in X_s$.
Here is the main theorem of this section:

\begin{theorem} 
	\label{thm:field}
	Let $(\field_r(u))_{r\in (0,r_0]}$, $u\in U$ be the fixed
	point family from Lemma
	\ref{lem:fixedfield}.
	Suppose that for every $0<s<r\leq r_0$
	and $t\in[0, r-s)$,
	the map
	\begin{equation*}
	u \in U \mapsto \BL_{s,u} \Bj_{s,r} \in L(X_r,X_s)
	\end{equation*}
	is $C^{t}$.
	Then, the mappings taking $u\in U$ to:
	\begin{equation*}
	\field_s(u)\in  X_s     , \ \ \ \
	\BLambda_u(\field) \in L^\infty(\Omega,\mathbb R_+) \ \ \
	\mbox{and} \ \ \
	\int_\Omega \log
	\Lambda_u(\field)_\omega \, d\mathds P(\omega)
	\end{equation*}
	are all $C^t$
	for $0< t<r_0-s$.
	In particular, it follows that the last expression
	(the characteristic exponent)
	is $C^{r_0-\epsilon}$ for
	any $\epsilon>0$.
\end{theorem}
\  \\

The proof will take the rest of Section \ref{sec:loss-of-reg} and
is divided into the following
steps taking us through our
scale of Banach spaces $X_r$, $r\in (0,r_0]$:
\begin{enumerate}
	\item
	First, by contraction and transversality
	we show H\"older 
	continuity 
	of $\field_s$.
	\item 
	This continuity will then allow us to 'bootstrap' and show 
	differentiability  
	of $u\mapsto \field_{s}$ when $0<s<r_0-1$.
	\item   Finally,
	a recursive argument will yield higher order 
	differentiability as well as the optimal
	H\"older exponent of the
	last derivative. 
	% Because of notational complications
	%we provide details in the case of the second derivative and
	%only a sketch of the proof for higher order.
\end{enumerate}

%\begin{rem}
%For simplicity of exposition we have assumed throughout 
%(H4) that our family of cocycles is uniformly bounded
%from below and above.
%This condition 
%may be relaxed quite a lot but at the expense of
%making the necessary
%regularity conditions 
%on the co-cycle somewhat illegible.
%\end{rem}

Step one, yielding H\"older-continuity, 
is done through 
\begin{lemme}
	Let $0<s<r_0$.  Then
	$u\in\U\mapsto \field_{s}(u) \in X_s$  is $\alpha$-H\"older continuous
	%for $\alpha \in (0,1] \cap (0, r_0-s)$.
	for $\alpha = 1 \wedge (r_0-s)$.
	\label{lem:alphacont}
\end{lemme}

\paragraph{Proof:}  
For $n\geq 1$ we make use of the following telescopic identity,
putting the downgrading at the right place:
\begin{equation*}
\Bj_{s,r_0}(\BL_{r_0,u+h}^n - \BL_{r_0,u}^n) = \sum_{k=0}^{n-1}
\BL_{s,u+h}^k
\left( \Bj_{s,r_0} (\BL_{r_0,u+h}-\BL_{r_0,u}) \right)
\BL_{r_0,u}^{n-k-1} \in L(X_{r_0};X_s).
\end{equation*}
The middle term is $\alpha$-H\"older from the assumptions 
in Theorem \ref{thm:field} and the other factors are bounded,
so $u\rr \Bj_{s,r_0} \BL^n_{r_0,u}$ is $\alpha$-H\"older.
When $\Bphi_{r_0}\in \BC_{r_0,\ell_{r_0}=1}$,
the same is true for $\Lambda_{u}^{(n)}(\Bphi_{r_0}) =
\langle \ell_s ,  \Bj_{s,r_0} \BL^n_{r_0,u} \Bphi_{r_0})$
and since the section is uniformly bounded
from below, also $u\in\U \mapsto \Bj_{s,r_0} \Bpi_{r_0,u}^{n}(\Bphi_{r_0})$
is $\alpha$-H\"older.
For $h\in \B$ small (so that $u+h\in\U$)
we get
by \eqref{CORA5} and \eqref{CORA8} in the appendix:
\begin{align*}
\|\field_s(u+h) &-\field_s(u)\|_{X_s}   = 
\| \Bpi_{s,u+h}^{n}(\field_s(u+h))- \Bpi_{s,u}^{n}(\field_s(u))\|_{X_s}
\\
&\leq \| \Bpi_{s,u+h}^{n}(\field_s(u+h))- \Bpi_{s,u}^{n}(\field_s(u+h))\|_{X_s}+
\| \Bpi_{s,u}^{n}(\field_s(u+h))- \Bpi_{s,u}^{n}(\field_s(u))\|_{X_s} \\
& \leq 
\| \Bj_{s,r_0}\Bpi_{r_0,u+h}^{n}(\field_{r_0}(u))-
\Bj_{s,r_0} \Bpi_{r_0,u}^{n}(\field_{r_0}(u))\|_{X_s}+
\frac{K}{\rho}\eta^{n} 
\|\field_s(u+h) -\field_s(u)\|_{X_s} .
\end{align*}
Now, fix 
a value of $n$ such that $\frac{K}{\rho}\eta^{n} \leq \frac12$ to  obtain 
\begin{equation*}
\|\field_s(u+h) -\field_s (u)\|_{X_s} \leq 2 \|
\Bj_{s,r_0}\Bpi_{r_0,u+h}^{n}(\field_{r_0}(u))-
\Bj_{s,r_0}\Bpi_{r_0,u}^{n}(\field_{r_0}(u))\|_{X_s}.
\end{equation*}
%Since $\Bh^{0}=\Bj^{0}(\Bh^{(-1)})$ 
And here the RHS is $\alpha$-H\"older continuous as we showed above.
\Halmos\\

\subsection{Differentiability and beyond}
For the second step we will base our proof on arguments 
already given in \cite{Sedro}.
For $s\in(0,r_0]$, $\Bphi_s\in\BC_{s,\ell_s=1}$
and $\Bz_s\in X_s$ we set:
\begin{equation}
Q_{s,u}(\Bphi_s).\Bz_s := \partial_{\Bphi} \Bpi_{s,u}(\Bphi_s).z_s =
\Lambda_u(\Bphi_s)^{-1} \left( 
\BL_{s,u} \Bz_s-  \Lambda_u(\Bz_s)  \Bpi_{s,u} (\Bphi_s) \right).
\label{eqforQ}
\end{equation}
The fact that this is well-defined follows from Lemma \ref{lem:Qop}
where it is also shown that
$\bfone - Q_{s,u}(\Bphi_s) \in L(X_s)$ is invertible
with the bound
$\|(\bfone- Q_{s,u})^{-1} \|_{L(X_s)}
\leq 1 + \frac{K_s/\rho_s}{1-\eta_s}$ and that
for $\|\Bz_s\|_s < \rho/2\vartheta K$:
\begin{equation*}
\| \Bpi_{s,u}(\Bphi_s+\Bz_s) -\Bpi_{s,u}(\Bphi_s)
- Q_{s,u}.\Bz_s  \|_s \leq  C(\|\Bz_s\|_s^2) .
\label{exp:BphiBpi}
\end{equation*}
with $C=C(\rho_s,\vartheta_s,K_s)<+\infty$.
When downgrading the regularity we recover regularity
as stated in 
the following
\begin{prop}
	\label{prop:piCt}
	Given $0<s<r\leq r_0$, $0\leq t <r-s$ we have 
	for $\Bphi_r\in \BC_{r,\ell_r=1}$ that %and $\Bz_r\in X_r$,
	%$\|\Bz_r\|<\rho/2\vartheta K$:
	\begin{equation*}
	u \in U \mapsto \Bj_{s,r} \Bpi_{r,u} (\Bphi_r) \in X_s
	\ \ \mbox{and} \ \ 
	u \in U \mapsto 
	\Bj_{s,r} Q_{r,u}(\Bphi_r) \in L(X_r;X_s)
	\end{equation*}
	are both  $C^t$.
\end{prop}
\Proof 
By hypothesis, the map $ u\in U \mapsto 
\Bj_{s,r} \BL_{r,u} \in L(X_r,X_s)$
is $C^t$.
By the uniform bounds 
$u\in U \mapsto 
\BLambda_{u}(\Bphi_r) := 
\langle \ell_s, \Bj_{s,r} \BL_{r,u} \Bphi_r \rangle$
is  a $C^t$ section over $\Omega$ taking values in
$\left[ \frac{\rho_s}{\vartheta K}, \vartheta K\right]
\subset (0,+\infty)$. So the same is true for its
(fiber-wise) inverse
$ \left( \BLambda_{u}(\Bphi_r)  \right)^{-1}$
and then also
for the fractional linear bundlemap 
$\Bpi_{s,u}(\Bj_{s,r}\Bphi_r) = 
\Bj_{s,r} \Bpi_{r,u}(\Bphi_r) = 
\BLambda_u(\Bphi_r)^{-1}
\Bj_{s,r} \BL_{r,u} \Bphi_r$.
From the expression (\ref{eqforQ}) for 
$\Bj_{s,r}Q_{r,u} = Q_{s,u}\Bj_{s,r}$ we see that
each term is $C^t$ and the result follows from the Leibniz principle.
\Halmos\\

\begin{lemme}
	\label{lem:Plem}
	Let  $0<s<r-1\leq r_0-1$ and $\Bphi_r\in \BC_{r,\ell_r=1}$.
	Then 
	\begin{equation}
	\Bj_{s,r} \left(\Bpi_{r,u+h}(\Bphi_r) -\Bpi_{r,u}(\Bphi_r)\right)
	= P_{s,r,u}(\Bphi_r) .h + \O_{s,r} \left(|h|^{1+\alpha}\right)
	\label{exp:uBpi}
	\end{equation}
	with $\alpha = (r-1-s)\wedge 1$ and
	\begin{equation}\label{eqforP}
	P_{s,r,u}(\Bphi_r) = \partial_u \left(\Bj_{s,r} \Bpi_{r,u}(\Bphi)\right) =
	\Lambda_u(\Bphi)^{-1} 
	\left( 
	\partial_u(\Bj_{s,r}\BL_{r,u} \Bphi_r )-
	(\partial_u \Lambda_u(\Bphi)) \Bpi_{s,u}(\Bphi)
	\right)
	\end{equation}
	Furthermore,
	the map $u\in \U \mapsto  P_{s,r,u}( \Bphi_r) \in L(\B,X_s)$ is
	$C^{t}$ for $t\in[0,r-s-1)$.
\end{lemme}
\paragraph{Proof:}  
From the hypothesis in Theorem \ref{thm:field} we
have that $\Bj_{s,r}\BL_{r,u} \in L(X_r;X_s)$
is $C^{1+\alpha}$. So applying the 
MVT (recall that $\U$ is convex) we obtain:
\begin{align*}
\|\BL_{s,u+h} \Bphi_s - & \BL_{s,u} \Bphi_s - 
h\cdot\partial_u \left( \Bj_{s,r}  \BL_{r,u} \right)
\Bphi_r\|_{X_s}    \\
& \leq  
\sup_{0\leq t\leq 1} |h|_\B
\; \| 
\| \partial_u  \Bj_{s,r} (\BL_{r,u+th} -
\BL_{r,u}) \|\; \| \Bphi_r\|_{X_r}  
=  \O_{s,r} \left(|h|^{1+\alpha}\right) \|\Bphi_r\|_{X_r}.
\end{align*}
Since
$\partial_u\Lambda_u(\Bphi)= \partial_u (\Lambda_{s,u}(\Bj_{s,r}\Bphi_r))$
we also have the expansion:
\begin{equation*}
\Lambda_{u+h}(\Bphi) - \Lambda_{u}(\Bphi)
- h\cdot (\partial_u \Lambda_u ) (\Bphi) 
= \O_{s,r}\left(|h|^{1+\alpha}\right) .
\end{equation*}
As $\Bphi$ is in the normalized cone,
$\Lambda_u(\Bphi)^{-1}$ is  bounded  and 
(\ref{exp:uBpi}) now follows from expanding the LHS using the
two previous estimates, giving the wanted form for $P_{s,r,u}(\Bphi_r)$.
\\Finally, from the hypothesis in \ref{thm:field}, 
$\partial_u(\Bj_{s,r}\BL_{r,u}) \Bphi_r$
is $C^{t}$ for $t\in[0,r-s-1)$,  whence also $\partial_u\Lambda_u(\Bphi)$ and
the claim follows.  \Halmos\\

\begin{theorem}[cf.\ Theorem 1 in \cite{Sedro}]
	\label{sedro thm}
	Let $0<s<r_0-1$. Then $u\in U\mapsto \field_s(u) \in X_s$
	is differentiable.
	The derivative is given by:
	\begin{equation}
	D_u \field_s(u) = \left( 1-Q_{s,u}(\field_s(u)) \right)^{-1} 
	P_{s,r_0,u}(\field_{\,r_0}(u)) \in L(\B,X_s),
	\label{formula:deriv}
	\end{equation}
\end{theorem}
\paragraph{Proof:} Fix $u\in U$
and set $\Bphi_t=\field_t(u)$ (all $0<t\leq r_0$).
For $h\in \B$ small (in particular,
$u+h\in U$) we set $\Bz_t=\Bz_t(h)=\field_t(u+h)-\field_t(u)$.
From the fixed point property, consistency and the above estimates we deduce:
\begin{align*}
\Bz_s &= \Bpi_{s,u+h} (\Bphi_s+\Bz_s)- \Bpi_{s,u}(\Bphi_s) \\
&=   \Bpi_{s,u+h}(\Bphi_{s}+\Bz_{s})- 
\Bpi_{s,u+h} (\Bphi_{s}) + 
\Bj_{s,{r_0}} \left(\Bpi_{r_0,u+h} (\Bphi_{r_0})-
\Bpi_{r_0,u}(\Bphi_{r_0}) \right)\\
&= 
Q_{s,u+h}(\Bphi_s).\Bz_s  + \O_s(\|\Bz_s\|^2) + P_{s,r_0,u}(\Bphi_{r_0}).h +
\O_{s,r_0}(|h|^{1+\alpha}) 
\end{align*}
By
Lemma \ref{lem:alphacont} and since $s<r_0-1$
we already know that $\Bz_s=\O_s(|h|)$
so the term  $\O_s(\|\Bz_s\|^2)$ can be neglegted.
We still need to convert $Q_{s,u+h}$ into $Q_{s,u}$. 
For $s+\alpha<t<r_0-1$ we have
$\Bz_t(h)=\O_t(|h|)$ and then
\begin{align*}
(Q_{s,u+h}(\Bphi_s)-Q_{s,u}(\Bphi_s)).\Bz_s
&= 
\Bj_{s,t}(Q_{t,u+h}(\Bphi_t)-Q_{t,u}(\Bphi_t)).\Bz_t(h)\\
&=
\O_{s,t}(|h|^{\alpha})  \times \O_{t}(|h|)= 
\O_{s,t}(|h|^{1+\alpha}),
\end{align*}
where we have inserted a downgrading operator to exploit
$\alpha$-H\"older continuity from the second term 
(Proposition \ref{prop:piCt})
and the Lipschitz continuity shown in \eqref{lem:alphacont} at the fixed point $\Bphi_s$ for the first term.
Thus,
$\Bz_s = Q_{s,u}(\Bphi_s).\Bz_s  +  P_{s,r_0,u}(\Bphi_{r_0}).h +
\O_{s,r_0}(|h|^{1+\alpha})$ and
using invertibility of $\bfone-Q_{s,u}$ we get:
\begin{equation*}
\Bz_s =
(\bfone-Q_{s,u}(\Bphi_s))^{-1} \left( P_{s,r_0,u}(\Bphi_{r_0}).h \right) \; 
+\;\O_{s,r_0} (|h|^{1+\alpha}) .
\end{equation*}
This shows differentiability of $\field_s(u)$ with the stated
formula for the derivative.
%To see that $h\mapsto D_u \field_s(u)$ is itself
%$\alpha$-H\"older continuous one has to joggle with downgrading
%as in the proof of Lemma \ref{lem:alphacont}. The proof, however,
%follows the same lines as the general proof of
%our main theorem, which we give below.
\Halmos\\

For higher order regularity, we will bootstrap through the
formula (\ref{formula:deriv}). Showing that the RHS has some regularity
means that $f_s(u)$ has the same regularity plus one.
Obtaining  a $C^{1,\alpha}$-regularity this way is not very hard, 
using equicontinuity and telescoping as above.
But already for $C^2$ the algebra starts getting quite involved.
To treat the general case, we have developped what
to us appears to be a new
algebraic differentiation tool:
A graded Leibnitz principle.
This allows us to use a simple inductive argument. 
A detailed description of the notation involved and
the relevant technical 
lemmas have been relegated to  Appendix \ref{app:graded-diff-calc}.
Our inductive hypothesis is the following:

\begin{hyp} 
	\label{hyp:gamma}
	For $\gamma \in (0,r_0]$ we let $\calH(\gamma)$ denote the
	following property: $\forall \; 0 < s<r_0$ and
	$0\leq t < \gamma\wedge (r_0-s)$, the map:\ 
	$ u \in \U \mapsto \field_s(u)= \Bj_{s,r_0} \field_{r_0}(u)
	\ \in X_s$ \ 
	is $C^{t}$.
\end{hyp}

What we have shown above is that $\calH(1)$ holds and that
in addition when $0<s<r_0-1$ then $u\mapsto \field_s(u)$ is differentiable
with a derivative given by formula (\ref{formula:deriv}).
For $0<s<r_0-1$,
by equivariance, the expression
$P_{s,u}(\field(u))= P_{s,r,u}(\field_r(u))$
is independent of $r\in (s+1,r_0]$.

\begin{lemme}
	\label{lem:induc}
	Under hypothesis $\calH(\gamma)$
	with $\gamma\in (1,r_0]$:
	\begin{enumerate}
		\item 
		The map:
		$u \mapsto  P_{s,u}(\field(u)) \in L(\B;X_s)$
		is $C^{t}$ for $0\leq t < \gamma\wedge(r_0-s-1)$.
		\item  The map:
		$u \mapsto 
		%R_{s,u}(\field_s(u)) :=
		\left( \bfone - Q_{s,u}(\field_s(u)) \right)^{-1}
		\Bj_{s,r} \in  L(X_r;X_s)$
		is $C^{t}$ for $0\leq t < \gamma\wedge (r-s)$.
	\end{enumerate}
\end{lemme}
The proof of this
lemma is given at the end of Appendix \ref{app:graded-diff-calc}.

\paragraph{Proof of Theorem \ref{thm:field}:} 
Write $r_0=k_0+\alpha_0$ with $k_0\in \N$ and $\alpha_0\in (0,1]$.
By the above $\calH(\alpha_0)$
is satisfied.
We will show that when $k\in \{0,1,\ldots,k_0-1\}$
then hypothesis $\calH(k+\alpha_0)$ implies 
$\calH(k+1+\alpha_0)$. This will prove that $\calH(r_0)$ holds
which is precisely the conclusion of the theorem.

Thus, assume that we have shown 
that $\calH(k+\alpha_0)$ holds
for some $k\in \{0,1,\ldots,k_0-1\}$.
Let $\eps>0$ be small, $r=r_0-1-\eps$ and $0<s<r$.
Setting 
$R_{s,u}:=\left( \bfone - Q_{s,u}(\field_s(u)) \right)^{-1} \in L(X_s)$
we may write 
\begin{equation*}
D_u \field_s(u) = R_{s,u}(\field_s(u)) \Bj_{s,r}
P_{r,u}(\field(u)) \in L(\B,X_s),
\end{equation*}

Here,
in terms of Appendix \ref{app:graded-diff-calc},
$M^1_{s,r}(u)=R_{s,u}(\field_s(u))\Bj_{s,r}$, $0<s<r\leq r_0$
defines an equivariant $(\gamma,0)$-regular
family of operators
and
$N_s(u) = \Bj_{s,r} P_{r,u}(\field(u))$  defines a left-equivariant
$(\gamma-1,1)$-regular family. 
By the graded Leibniz principle,
Proposition \ref{prop:prod},
the above product
is $(\gamma-1,1)$-regular so in particular
$u\mapsto D_u \field_s(u)$ is $C^t$ for
$0\leq t < \gamma\wedge (r-s-1)$. Therefore,
$u\mapsto \field_s(u)$ is $C^t$ for 
$0\leq t < (\gamma \wedge (r-s))+1 = (\gamma+1)\wedge(r_0-s-\eps)$.
We may here let $\eps$ go to zero
and conclude that $\calH(\gamma+1)$ holds. 
As already mentioned this
implies $\calH(r_0)$, thus concluding the proof. \Halmos\\

\section{Dynamical applications}
\label{sec:dynapp}

We illustrate the above results through 
applications to dynamical systems.
Let $(M,g)$ be an
$n$-dimensional\footnote{ The manifold could also be
	modeled over any Hilbert space. Finite dimension is however used
	when dealing with physical measures.}
$C^\infty$ connected  Riemannian
manifold without boundaries and of bounded diameter.
Also let $\lambda>1$ be a fixed expansion constant.
Further geometric conditions on an atlas for $M$ will be given below.
We will work with parametrized (Bochner)-measurable families of 
$C^r$ uniformly $\lambda$-expanding maps on $M$
and a continuous scale 
of Hölder spaces $(C^s(M))_{0<s\leq r_0}$ for some $r_0>1$.
We consider parameters taking
values in an open non-empty convex subset $\U$ of a Banach space $\B$.
Our strategy will be to
exhibit a regular Birkhoff cone in every $C^s$ space
which is contracted in a strict and uniform way by the transfer
operators associated with the expanding maps. 
This will allow us to construct the fixed point
$\field_u\in L^\infty(\Omega,C^s(M))$, $u\in \U$ of
a transfer operator cocycle, and then to study
regularity with respect to $u\in \U$
of the stationary measure and the top characteristic exponent
with the tools of the previous sections.
\medskip

When $M$ is a torus or an interval (as in our examples in the introduction)
the tangent bundle is trivializable, which simplifies the discussion 
when dealing with expanding maps.  Our theorems covers more
general situations and e.g.\ 
in the case of expanding maps on (infra-)nil-manifolds\footnote{
	As shown by Gromov in \cite{GR81} 
	any expanding self-map of a compact manifold is topologically conjugate
	to an infra-nil-endomorphism.}, 
the tangent bundle may not be trivializable and more care is needed.
The construction of  $C^{(k,\alpha)}$ structures on general manifolds
for $k\geq 1$  and $\alpha\in (0,1]$
($\alpha$-H\"older continuous $k$'th derivative) 
is more delicate
due to the non-locality of the H\"older-continuity.
%It is convenient  to set up good local coordinates.
Eldering, in \cite{El13},
uses the exponential map to provide almost 'flat'
local coordinates around base-points. 
The approach we take here was mostly inspired
by Ruelle \cite{Ru90}.

%%%%%%%%%%%%%%%%%%%%%%%%%
%     No longer needs path-distance ?? 
% For an open connected subset $W\subset M$
% we write $d_W$ for the path distance of points, i.e.\
% the infimum Riemannian length of a $C^1$-curve contained
% in $W$ and joining two 
% given points in $W$. 
% In particular, with $W=M$ it is just
% the manifold distance.

\begin{hyp}
	\label{hypAtlas}
	There exists a $C^\infty$-atlas 
	$\calA_M = \{(\V_j,\eta_j) : j\in\J \}$
	for $(M,g)$ with the following properties:
	
	\begin{enumerate}
		
		\item The atlas consists of simply connected charts and
		is uniformly $C^\infty$.
		By 'uniformly' we 
		mean that there is a sequence $(c_k)_{k\geq 1}$ of positive reals
		so that
		every transition map between overlapping charts
		(thus a map between open sets in Banach spaces):
		\begin{equation*}
		\eta_j\circ \eta_m^{-1}  : \eta_m(\V_j\cap \V_m) \rr \eta_j(\V_j\cap \V_m)
		\ \ (\subset \R^n)  \ , \ \ j,m\in \J
		\end{equation*}
		has a $C^k$ norm bounded by $c_k$. 
		
		\item
		\label{hypJ}
		There is $\delta_0>0$ so
		that  $\forall y\in M$ the set
		$\ds \J(y) \stackrel{\rm def}{:=}  \{ j\in \J : B(y,\delta_0)\subset \V_j \}$
		is non-empty.
		
		\item 
		\label{hypConvex}
		Each chart $V_j=\eta_j(\V_j)$, $j\in\J$ is convex in $\R^n$.
		
		\item
		\label{equivEuclid}
		There is $L<+\infty$ so that $\forall j\in\J$,
		$x,x'\in \V_j$, $x_j=\eta_j(x), x_j'=\eta_j(x')$:
		\begin{equation*}
		\frac{1}{L} \|x_j-x'_j\|_{\R^n}
		\leq d_M(x,x') 
		%\leq d_{\V_j}(x,x')
		\leq L  \| x_j - x'_j \|_{\R^n}.
		\end{equation*}

		In the following we consider a fixed atlas verifying the above conditions and
		we fix a value for the Lebesgue number $\delta_0>0$.
		We set $\delta_1=\delta_0/3$.
		
	\end{enumerate}
\end{hyp}

\begin{rem}
	\begin{enumerate}
		\item
		Our functional analysis
		will be built with respect to this fixed atlas, but
		our main results are
		independent of the choices made
		(though constants in estimates may change).
		We also note that the existence of
		an atlas verifying the above is automatic when $M$
		is a compact $C^\infty$ manifold without boundary.
		\item 
		The convexity condition \ref{hypAtlas}.\ref{hypConvex} on the $V_j$'s may be relaxed, but to find a good geometric condition is complicated, notably for the result \ref{cor:compodiff} in the Appendix to hold.
		\\We refer to e.g.\ De La Llave and Obaya in \cite[section 6]{RDLL99} for sufficient conditions in this direction. We have here opted for the simplest hypothesis.
	\end{enumerate}
\end{rem}

As mentioned above
$\U\subset \B$ is a non-empty open convex
subset (our parameter space)
in a Banach space. Then $U\times M$ inherits a manifold structure
with an atlas consisting of charts $(U\times \V_j, \teta_j)_{j\in \J}$ 
in which
\begin{equation*}
\teta_j : (u,y) \in U\times \V_j \mapsto (u,\eta_j(y)) \in U \times V_j.
%\label{def:teta}
\end{equation*}

Throughout this section we equip $\R^n$ with the Euclidean metric
and $\B\times \R^n$ with the norm 
$\|(u,x)\|= |u|_{\B} \vee |x|_{\R^n}$ with "$\vee$" meaning max.

Let  $k\geq 0$ be
an integer and $\alpha\in(0,1]$. 
Let $Z$ be a Banach space and for $j\in\J$ let
$\phi_j : \U \times V_j \rr Z$ be a  $C^k$ map 
(thus from an open convex subset
of a Banach space into a Banach space).
We then define  (Fr\'echet) derivatives of $\phi_j$
in the usual way and 
for  $0\leq q\leq k$ the
associated semi-norm of $\phi_j$:
\begin{equation*}
\myp^{(q)}_{}(\phi_j)  = 
\sup\{ \|D^q \phi_j (u,x) \| : u\in U, x \in V_j \}
\end{equation*}

Recall that $\delta_1=\delta_0/3$ (Hypothesis \ref{hypAtlas}.\ref{hypJ}).
We will define a $\delta_1$-local H\"older-continuity of $D^k\phi_j$ 
in the chart $\U\times V_j$. 
For $\xi=(u,x_j),\xi'=(u',x_j') \in U\times V_j$ we set:
\begin{equation*}
d_j (\xi,\xi') = |u-u'|_{\B} \vee |x_j - x'_j|_{\R^n} ,
%\label{induceddist}
\end{equation*}
We say that  $D^k \phi_j$ is $(\alpha;\delta_1)$-H\"older continuous
if there is $C<+\infty$ so that:
$\forall\; \xi,\xi'\in U\times V_j$ with
$d_j(\xi,\xi')\leq \delta_1$:
\begin{equation}
\|D^k \phi_j (\xi) -D^k \phi_j (\xi')\|  \leq
C \ d_j(\xi,\xi')^\alpha.
\label{eq:holderineq}
\end{equation}
We let
$\myh^{(k,\alpha)} (\phi_j)
:=\myh^{(k,\alpha)}_{\delta_1,d_j} (\phi_j)$ be
the smallest constant $C$ for which (\ref{eq:holderineq}) holds
and define the norm:
\begin{equation}
\|\phi_j\|_{(k,\alpha)} := 
\|\phi_j\|_{(k,\alpha)}^{\delta_1,d_j} = 
\max_{0\leq q\leq k}
\myp^{(q)}(\phi_j) 
\  \vee \ \myh^{(k,\alpha)}_{\delta_1,d_j}(\phi_j).
\label{def:chartCka}
\end{equation}

As shown in Lemma \ref{holdereq} making another choice of $\delta_1$
and the metric $d_j$, if equivalent to the above,
will lead to an equivalent $\|\cdot\|_{(k,\alpha)}^{\delta_1,d_j}$-norm.
Unless explicitly stated otherwise $\delta_1$ and the metric $d_j$ are
as above and will be omitted in the notation for the Banach spaces and
the norms.\\

\begin{definit}
	We let $C^{(k,\alpha)}(\U\times V_j;Z)$ denote the
	space of functions $\phi_j : U\times V_j \rr Z$ for which the norm 
	\eqref{def:chartCka} is finite. $(C^{(k,\alpha)}(\U\times V_j;Z),\|\cdot \|_{(k,\alpha)})$ is a Banach space.
	\label{def:Ckad}
\end{definit}

A function on the parametrized manifold, $\Phi:\U\times M \rr Z$, induces a collection $(\phi_j)_{j\in\J}$ of partial maps on charts: 
for each  $j\in\J$, $\phi_j:=\Phi\circ (\teta_j)^{-1} : \U \times V_j\rr Z$.
This leads to the following:
\begin{definit}
	\label{def:cka}
	The $(k,\alpha)$-norm of  $\Phi:\U \times M\rr Z$ is defined as:
	\[\|\Phi\|_{(k,\alpha)}= \sup_{j\in \J} 
	\|\phi_j\|_{(k,\alpha)}.\]
	In the following, $C^{(k,\alpha)} (\U\times M;Z)$ denotes the Banach space of functions $\Phi$ for which the above norm is finite, and $c^{(k,\alpha)}(\U\times M;Z)$ denotes the closure of $C^\infty$ functions in $C^{(k,\alpha)} (\U\times M;Z)$.
	\\$C^{(k,\alpha)} (M;Z)$ and $c^{(k,\alpha)}(M;Z)$ are defined in the same way (omitting $\U$ in the construction).
\end{definit}

\begin{rem}
	\label{Cs-structure}
	The above construction yields the
	following description of 
	$C^{(k,\alpha)}(\U\times M)$ (and similarly for the other
	spaces in the above definition):
	A family $(\phi_j)_{j\in \J}$
	with each $\phi_j\in C^{(k,\alpha)}(\U\times V_j)$ corresponds to a (unique)
	element $\Phi\in C^{(k,\alpha)}(\U\times M)$ if and only if we have the following
	identification:
	\begin{equation*}
	\forall j,m\in\J, \ u\in \U, \  y \in
	\V_j\cap \V_m
	\ : \ \phi_j (u, \eta_j(y)) 
	= \phi_m(u,  \eta_m(y)) .
	\end{equation*}
	Denoting 
	this identification by $\sim$, we may write:
	%\begin{equation*}
	$\ds          C^{(k,\alpha)}(\U\times M) \simeq 
	\left( \prod_{j\in\J} C^{(k,\alpha)}(\U\times V_j) \right) / \sim $.
	%\end{equation*}
\end{rem}
We also want to define a local $C^{(k,\alpha)}$-structure (and a $c^{(k,\alpha)}$-structure) 
for parametrized locally defined maps of $M$ into itself which satisfies a uniform
Lipschitz condition. Again we shall do this in terms of our fixed atlas.
\\More precisely, for each $j\in \J$ consider  maps of the type $\locinv: \U \times \V_j \rr M$ that verifies the following uniform Lipschitz-condition:
\begin{equation}
d_M(\locinv(u,y),\locinv(u',y')) \leq  
|u-u'|_{\B}+\frac{1}{\lambda} d_M(y,y')
\label{ineq:psi}
\end{equation}
for all $(u,y),(u',y') \in \U \times \V_j$.
Thus, the map is $1/\lambda$-contracting in the 
manifold variable and $1$-Lipschitz in the $u$-variable.
More generally, it suffices to assume that the map
is locally  Lipschitz in the $u$-variable.
But by restricting to a smaller subset of $U$
(our results are only local with respect to $U$) and
scaling the norm in $\B$ we may always reduce to the above case.
Given two maps 
$\locinv,\tlocinv : \U\times\V_j\rr M$ we define their uniform
distance to be:
\begin{equation*}
d_0(\locinv,\tlocinv) = \sup_{\zeta\in \U\times \V_j} 
d_M(\locinv(\zeta),\tlocinv(\zeta)).
\end{equation*}
The following lemma will allow us to compare in a finer way
uniformly close maps:

\begin{lemme}
	\label{lem:inclusion}
	For $j\in \J$ 
	let $\locinv,\tlocinv : \U\times\V_j\rr M$ 
	verify (\ref{ineq:psi}) and $d_0(\locinv,\tlocinv)\leq \delta_1$.
	Then for any $\zeta_0=(u_0,y_0)\in \U\times \V_j$
	and $m\in \Jd(\locinv (\zeta_0))\cup \Jd(\tlocinv (\zeta_0))$ we have:
	\begin{equation*}
	\locinv, \tlocinv : 
	B_{\B\times {\V_j}} (\zeta_0,\delta_1) \rr 
	B_{M} (\locinv (\zeta_0),\delta_0) \cap
	B_{M} (\tlocinv (\zeta_0),\delta_0)
	\subset \V_m.
	\end{equation*}
\end{lemme}
\Proof $\locinv,\tlocinv$ are $2$-Lipschitz and $d_0(\locinv,\tlocinv)\leq \delta_1$. So
the result follows from $0<\delta_1<\delta_0/3$ and Hypothesis
\ref{hypAtlas}
on the atlas.\Halmos\\

Thus, when
$\locinv,\tlocinv$ are uniformly close,
we may go to the same local coordinates  in the source
and in the image. 
Let us write 
$\psi_j := \locinv \circ (\teta_j)^{-1}
: U\times V_j \rr M$ and
$\tpsi_j := \tlocinv \circ (\teta_j)^{-1}
: U\times V_j \rr M$
when expressing the above two maps in local coordinates in the source
(but not the image).
\\Given $\xi_0=(u_0,y_0)\in \U\times V_j$ (thus in local coordinates)
let us write
\begin{equation}
B_j(\xi_0) = \{(u,y)\in \U \times V_j :
d_j((u,y),(u_0,y_0))  <\delta_1 \}
\label{def:Bxiball}
\end{equation}
for the product ball in charts.
%Let us write
%$\psi_j=\psi(\cdot,\eta_j^{-1}(\cdot)): \U\times V_j \rr  M$
%for $\psi$ expressed in local coordinates in the domain.
For any $m\in \J(\psi_j(\xi_0))\cup \J(\tpsi_j(\xi_0))$ we may 
go to local coordinates also in the image if we restrict the domain
to $B_j(\xi_0)$ and define
$\psi_{m,j}=(\eta_m\circ \psi_{j})_{ | B_j(\xi_0)}$
and 
$\tpsi_{m,j}=(\eta_m\circ \tpsi_j)_{| B_j(\xi_0)}$.
Then  $\psi_{m,j}, \tpsi_{m,j} :  B_j(u_0,y_0) \rr V_m\subset \R^n$
are local maps between 
(the same) convex subsets of normed
vector spaces, whence may be compared. \\

Given $\locinv,\tlocinv:\U\times \V_j\rr M$, or equivalently,
$\psi_j,\tpsi_j:\U\times V_j\rr M$ for which
$d_0(\psi,\tpsi)
\leq  \delta_1$ we define
their $(k,\alpha)$-distance:
\begin{equation*}
d^{(k,\alpha)}_{j}
(\psi_j,\tpsi_j) =
\sup_{\xi\in U\times V_j} \ \ 
\sup_{m\in \Jd(\psi_j(\xi)) \cup \Jd(\tpsi_j(\xi))}
\|\psi_{m,j}-\tpsi_{m,j}\|_{C^{(k,\alpha)}(B_j(\xi);V_m)},
\label{def:kaj-dist}
\end{equation*}
with local maps defined as in the previous paragraph.
When $(k,\alpha)=(0,1)$ that's it.
For $k\geq 1$, however, we also introduce the following 'gauge'
on $\locinv$ (or equivalently, $\psi_j)$:

\begin{equation}
p^{(k,\alpha)}_j(\psi_j) :=
%p_{{(k,\alpha)}}(\psi_j) =
%\|\psi\|_{C^{(k,\alpha)}(\U\times \V_j)}
\sup_{\xi\in U\times V_j} \ \ \sup_{m\in \Jd(\psi_j(\xi))  }
\|D\psi_{m,j}\|_{C^{(k-1,\alpha)}(B_j(\xi);V_m)}.
\label{eq:psikaj}
\end{equation}
We denote by
$\left( C^{(k,\alpha)}(\U\times \V_j;M),
d^{(k,\alpha)}_j \right) $
the  metric space of maps verifying (\ref{ineq:psi})  
and of finite gauge in the case $k\geq 1$.
Again $c^{(k,\alpha)}(\U\times \V_j;M)$ is 
the subspace
of maps which may be uniformly approximated by
$C^\infty$ functions. \\

\subsection{Uniformly expanding maps and associated weights}
\label{sec:expmaps}
In order to specify regularity of an expanding map, e.g.\
for the examples in the introduction,
it is convenient
to do so for its inverse branches rather than the map itself.
The following definition may look like an overkill but it allows for 
countably many branches and the presence of (certain types of)
singularities.

\begin{definit}
	Let  $\calA$ be a finite or countable index set.
	For each $j\in\J$ let $\Psi_j =\{\locinv_j^{\,i} : i\in \calA\}$
	be a family of maps $\locinv_j^{\,i} : \U \times \V_j \rr M$, $i\in\calA$
	satisfying (\ref{ineq:psi}). 
	
	\begin{enumerate}
		\item We say that two such
		families $\Psi_j$ and $\Psi_m$ are 
		compatible if $\;\forall u\in \U, y \in \V_j\cap \V_m$,
		there exist $\sigma\in \Aut(\calA)$ and an open neighborhood
		$W\subset \U \times (\V_j\cap \V_m)$ of $(u,y)$ 
		such that
		\begin{equation*}
		\forall i\in\calA : \ \ \
		\locinv_j^{\,i} \equiv \locinv_m^{\,\sigma(i)}    \ \ {\rm on} \ W.
		\end{equation*}
		\item
		Let the collection $\Psi=(\Psi_j)_{j \in\J}$ consist of
		compatible families. For any given $(u,y)\in \U\times M$
		there is $j\in \J$ so that
		$y\in \V_j$.
		By compatibility the set
		$D(u,y):= \{\locinv_j^{\,i}(u,y):i\in \calA\}$
		is then independent of the choice of $j$.
		The  collection is said to be separating
		if for every $u\in \U$: 
		\begin{enumerate}
			\item $D(u,y)\cap D(u,y')=\emptyset$ whenever $y\neq y'\in M$ and
			\item The map $i\in\calA \mapsto \locinv_j^{\,i}(u,y)\in M$ is 
			injective for every $j\in\J$ and $y\in \V_j$.
		\end{enumerate}
	\end{enumerate}
\end{definit}
For fixed $u\in\U$ the definition is designed
so that
$\ds G_u= \bigcup_{y\in M} D(u,y) \times \{y\}\subset M\times M$
is the graph of a map $T_u$ over some domain.
\\Consider the two natural projections 
$\pi_1: G_u\rr D_u:=\bigcup_{y\in M} D(u,y)$
and
$\pi_2:G_u \rr M$. By separability, the first is injective
so indeed $G_u$ is the graph of the map
$T_u = \pi_2 \circ \pi_1^{-1} : D_u \rr M$, defined on the domain $D_u$.
By construction, it is a covering map of $\M$ of
degree $\Card\, \calA$ which is
locally $\lambda$-expanding.
For each $j\in \J$ the family
$\locinv_j^{\,i}:\U\times \V_j\rr M$, $i\in \calA$ is precisely the 
(parametrized) $\frac{1}{\lambda}$-contracting inverse branches
of $T_u$ on the chart $\V_j$.
\\

Consider a compatible and separating collection of maps $\Psi$ 
and the associated dynamical system $T_u:D_u\rr M$  from
the above construction. As before we write
$\psi_j^i := \locinv^i_j \circ (\teta_j)^{-1} : U\times V_j\rr M$
with the source expressed in local coordinates.
For $k\geq 1$ we now make
a $(k,\alpha)$-regularity assumption on the branches:
We require that there is a constant $k_T<+\infty$ so that
our gauges, cf. (\ref{eq:psikaj}), are all uniformly bounded:
\begin{equation}
\forall i\in \A, \; j\in \J \ : \ p^{(k,\alpha)}_j(\psi^i_j) \leq k_T. \label{eq:kT}
\end{equation}
We will also associate 
a weight (not to be confused with the metric tensor)
$g_u:D_{u} \rr \R$ defined on the domain of $T_u$.
The most convenient way to specify regularity of the weight $g$ is through
its composition with
inverse branches. Again we do so using the charts in our fixed atlas.
So with $j\in\J$ and $\psi^{\,i}_j: \U\times V_j \rr M$ 
being a local
inverse branch we define
$g_{j}^i(u,y):=
g(u, \psi^i_j(u,y))$ for 
$(u,y)\in \U\times V_j$.
\\We will assume that each $g_j^i\in c^{(k,\alpha)}(U\times V_j;\R)$
and that when $k\geq 1$
there is a constant $k_g<+\infty$ so that
the derivative verifies the uniform bound:
\begin{equation}
\| D(g_j^i) \|_{(k-1,\alpha)} \leq k_g, 
\ \ i\in\calA, j\in\J.
\label{eq:condg}
\end{equation}
Finally, we require that there is a uniform
constant $1\leq k_\Sigma<+\infty$ so that for $(u,y)\in \U\times M$:
\begin{equation}
\frac{1}{k_\Sigma}\leq  \sum_{x : T_u(x)=y } e^{ g(u,x)} \leq k_\Sigma.
\label{eq:condS}
\end{equation}

\begin{rem}
	\begin{enumerate}
		\item
		Below we shall make use of the following observation: 
		Let $\xi_0=(u_0,y_0)\in \U\times V_j$.
		Each $g_j^i$ is by hypothesis $k_g$-Lipschitz
		and if $\xi=(u,y)\in B_j(\xi)$
		then $d_j(\xi,\xi_0)\leq \delta_1$.
		Thus, 
		\begin{equation}
		\sum_{i\in\A }\ \sup_{B_j(\xi_0)} \ e^{ g_j^i} 
		\ \leq \ k_\Sigma \ e^{\delta_1 k_g}.
		\label{eqn:unifexp}
		\end{equation}
		\item
		Finally, note that when $\calA$ is infinite, \eqref{eq:condS} implies that $g$ itself can not be uniformly bounded. 
		This is why we (have to) use the derivative of $g_j^i$ in the regularity condition \eqref{eq:condg}.
	\end{enumerate}
\end{rem}

\subsection{Pairs of maps and weights}\label{subsec:cocyclemeasurability}
We fix in the following a real value
of $r_0>1$. We set $\undr_0=(k_0,\alpha_0)$ with $(k_0,\alpha_0)\in \N\times (0,1]$ and $\delta_1>0$ as in the previous section.
We also choose arbitrary but fixed values of the associated constants  $k_T,k_g,k_\Sigma$ as described in \eqref{eq:kT}, \eqref{eq:condg} and \eqref{eq:condS}.
We consider the collection $\P=\P_{k_T,k_g,k_\Sigma}$ of $c^{(k_0,\alpha_0)}$ pairs $(T,g)$ with $T$ being a parametrized family of $\lambda$-expanding covering maps and $g$ is an associated weight as constructed above, for which the above regularity conditions are verified with the afore-mentioned constants.
\\We define a $(k,\alpha)$-distance $d_\P$ between pairs $(T,g), (\tT,\tg)\in \P$ in the following way:
First, if the two covering maps $T,\tT$ have different degrees then the distance is infinity. 
\\If they have the same degree we let $\calA$ be a common index set and consider for each $j\in\J$ the two sets of inverse branches (expressed in local coordinates) $\psi^i_j: U \times V_j\rr M$ and $\tpsi^i_j : U \times V_j\rr M$, $i\in\calA$.
We will compare these up to a (simultaneous) permutation of branches and weights. 
Let $P$ be the subset of 
permutations $\sigma\in {\rm Aut}(\calA)$ for which
$d_0(\psi^i_j, \tpsi^{\sigma(i)}_j)\leq \delta_1$  for all $i\in\A$.
\\For $\sigma\in P$ we set:
\begin{equation*}
d_{j,\sigma}=\sup_{i\in \A} \left\{d^{(k,\alpha)}_j (\psi^i_j, \tpsi^{\sigma(i)}_j)
\vee \|g_j^i-\tg_j^{\sigma(i)} \|_{(k,\alpha)}\right\}
\label{eq:djsigma}
\end{equation*}
and
\begin{equation}
d_{\P}((T,g),(\tT,\tg)) = \sup_{j\in\J} \ \inf_{\sigma\in P} d_{j,\sigma}.
\label{eq:dka-dist}
\end{equation}
If $P$ is empty the distance is infinity.
This provides a suitable notion of distance
on our collection
of pairs in $\P$\footnote{$d_{\P}$ does
	not necessarily verify the triangle inequality. One could turn it into a real
	metric by a Kobayashi-like construction, but this 
	does not seem to be of any particular use in our context.}.

\begin{lemme}
	\label{lemma:Lcont}
	Let $0<s\leq r_0$ and write $\unds=(k,\alpha)$ with 
	$s=k+\alpha$, $k\in\N$ and $\alpha\in(0,1]$.
	To every pair $(T,g)\in \P$ there is a well-defined
	bounded linear parametrized
	transfer operator $\hL_{T,g} \in L( c^\unds(M); c^\unds(\U\times M))$
	which for
	$\phi\in c^\unds(M)$ is given by the expression:
	\begin{equation*}
	(\hL_{T,g} \phi)(u,y) = \sum_{x\in D_u: T_u(x)=y} e^{g(u,x)} \phi(x),
	\ \ \ (u,y)\in \U\times M.
	\end{equation*}
\end{lemme}
\Proof For $j\in \J$ and $y\in V_j$  we may rewrite the expression as a sum
over inverse branches using local coordinates on $V_j$
(recall that $g_j^i(u,y) = g(u, \psi_j^i(u,y))$):
\begin{equation*}
(\hL_{T,g,j} \phi)(u,y)
= \sum_{i\in\A}
e^{g_j^i(u,y)} \phi( \psi_j^i(u,y)), \ \ (u,y) \in U \times V_j.
\end{equation*}
Let $(u_0,y_j)\in \U\times V_j$ and define 
as in \eqref{def:Bxiball}
$B_j=B_j(u_0,y_j) \subset \U\times V_j$.
It suffices to show that the $C^{\unds}$ norm
of $\hL_{T,g}\phi$ is uniformly bounded when restricted to
$B_j$ for arbitrary $(u_0,y_j)\in\U\times V_j$.
This may be done going to local coordinates.
For each $i\in\calA$, let $x^i=\psi_j^i (u_0,y_j)$ and pick
$m_i\in \Jd(x_i)$ (cf.\ Hypothesis \ref{hypAtlas}). 
The above sum then becomes
\begin{equation}
(\hL_{T,g,j} \phi)(u,y)
= \sum_{i\in\A}
e^{g^i_j(u,y)} \phi_{m_i} \circ \psi_{m_i,j}^i(u,y),
\label{eq:Linchartj}
\end{equation}
with $\phi_{m_i}=\phi\circ\eta_{m_i}^{-1}$, $\psi_{m_i,j}^i : B_j\rr V_{m_i}$ being the induced maps, now between open sets in Banach spaces.
Using Proposition \ref{prop:banalg} we see that all terms in the above sum are $c^{\unds}(B_j)$. 
From the estimates in that proposition and \eqref{eqn:unifexp} we get the following uniform bound:
\begin{eqnarray}
\|\hL_{T,g,j} \phi\|_s  &\leq &
\sum_{i\in\A}
C_s 
\left| e^{g^i_j}\right|_{C^0(B_{j})} 
\|\phi_{m_i}\|_s
(1+k_g)^s (1+k_T)^s   \nonumber\\
&\leq&
C_s 
(1+k_g)^s (1+k_T)^s  k_\Sigma \; e^{\delta_1 k_g } \;
\|\phi\|_s.
\label{eq:unifbdop}
\end{eqnarray}
where $C_s$ is a constant depending only upon $s$
and the constants related to the manifold and the choice of norms.
(and not to the specific choice of $(T,g)$).
The above bound is first for
each $B_{j}$, but then extends
to all of $\U\times M$ by definition of 
our $c^s(\U\times M)$-structure.
\Halmos \\

We then have the following result:
\begin{lemme}\label{lemme:continSOT}
	With the notation as in Lemma \ref{lemma:Lcont} above,
	for any fixed $\phi\in c^{\unds}(M,\CX)$,
	the map
	$(T,g)\in\P
	\mapsto
	\hL_{T,g}\phi\in c^\unds(\U\times M)$ is uniformly continuous.
\end{lemme}
\Proof Let $\eps>0$ and $\hphi\in c^\unds(M;\CX)$. 
First, by the uniform bounds
on the operator (\ref{eq:unifbdop}) and since
$\hphi$ is the $C^s$-limit of smooth functions,
we may find $\phi\in C^{s+1}(M,\CX)$
so that $\|\hL_{T,g}\phi -\hL_{T,g}\hphi\|_s \leq \eps/2$ for every
$(T,g)\in\P$. We fix this $\phi$.

Let $\Delta \in (0,\delta_1) \subset (0,1]$ and  consider two pairs 
$(T,g), (\tT,\tg)\in \P= \P_{k_T,k_g,k_\Sigma}$ of distance
$d_\P((T,g),(\tT,\tg))< \Delta$.
Let $j\in\J$
and for $(u_0,y_j)\in\U\times V_j$
set $B_{j}=B_j(u_0,y_j)$ as before.
By definition of $d_{\P}$ there is 
a common index set $\calA$ and
a permutation $\sigma$ such that
for every $i\in\A$
there is $m_i\in \J$  and local inverses:
$\psi_{m_i,j}^i,\tpsi_{m_i,j}^{\sigma(i)}:B_j\rr V_{m_i}$,
such that 
$d_j^{(k,\alpha)} (\psi_{m_i,j}^i,\tpsi_{m_i,j}^{\sigma(i)}) < \Delta$ 
and 
$\|g_j^i-\hatg_j^{\sigma(i)}\|_{(k,\alpha),j}  < \Delta$.
We may then write: 
\begin{align}
\hL_{T,g,j}\phi-\hL_{\tT,\tg,j}\phi &=
\sum_{i\in\calA}(e^{g_j^i}-e^{\tg_j^{\sigma(i)}})
\phi_{m_i}\circ\psi^i_{m_i,j} \label{eq:diff1}\\
\notag&+\sum_{i\in\A}e^{\tg_j^{\sigma(i)}}
(\phi_{m_i}\circ\psi^i_{m_i,j}-\phi_{m_i}\circ\tpsi_{m_i,j}^{\sigma(i)}).
\end{align}
Taking the $C^s$ norm and using Proposition \ref{prop:banalg} 
the first sum in \eqref{eq:diff1} is bounded by
\begin{equation*}
\sum_{i\in\calA}
C_{1,s}C_{4,s}|e^{1+g^i_j}|_{0}
(1+k_g)^s 
\Delta 
\  \
C_{3,s}
(1+k_\psi)^{s} 
\|\phi\|_{s} \ \leq\  C_1 \; \|\phi\|_{s+1} \;\Delta,
\end{equation*}
with $C_1=C_1(s,k_T,k_g,k_\Sigma,\delta_1)$
depending only upon the constants.
For the second term we note that by Lemma \ref{lem:inclusion}
the images of $\psi^i_{m_i,j}$ and $\tpsi^{\sigma(i)}_{m_i,j}$
are contained in the same convex set $V_m$.
Corollary \ref{cor:compodiff} then applies and implies:
\begin{equation*}
\sum_{i\in\calA}C_{2,s}|e^{\tg_j^{\sigma(i)}}|_{0}
(1+k_g)^s  \ \ C_{s,1} \|\phi\|_{s+1} \Delta
\leq C_2\; \|\phi\|_{s+1}\; \Delta,
\end{equation*}
again with $C_2$
depending only upon constants.
Choosing $\Delta$ small enough
each term can be made smaller than $\eps/4$,
thus concluding the proof. 
\Halmos\\

\begin{rem}
	How small we have to take 
	$\Delta$ depends in an
	essential way upon $\hphi$ through the norm $\|\phi\|_{s+1}$
	of the approximant.
	The theorem
	is in general false if $c^\unds$ is replaced by $C^\unds$.
	To see what may go wrong consider the following explicit example:
	For $s=(k,\alpha)$ with $k\in \N$ and $0<\alpha<1$,
	the function $\phi(x)=|x|^{k+\alpha}$
	belongs to  $\phi\in C^{\unds}([-1,1])$.
	However, the map $ u\in (-1,1)\mapsto 
	    \phi (\frac{u+x}{2} ) \in C^{\unds}([-1,1])$
	    is no-where continuous. We leave the verification of this fact
	    to the interested reader.

\end{rem}

\medskip

Let $(\Omega,\mathcal F,\mathds P)$ be a probability space, 
endowed with a measure-preserving, invertible and
ergodic map $\tau:\Omega\to\Omega$. 
Let $\omega\in\Omega\mapsto (T_\omega,g_\omega)\in \P=(\P_{k_T,k_g,k_\Sigma},d_{\P})$
be a Bochner measurable family of pairs.
\\By this we mean that for every $\Delta>0$
there is a countable measurable partition
$(\Omega_i)_{i\in\N}$ of $\Omega$ such
that the diameter in $\P$ for the
$d_{\P}$ distance (\ref{eq:dka-dist}) of each set 
$\{(T_\omega,g_\omega): \omega\in \Omega_i\}$ 
is smaller than $\Delta$.
In particular, the degree of $T_\omega$ is constant on each $\Omega_i$.
\\For $\omega\in\Omega$ and $u\in U$ we write $T_{\omega,u}(x) := T_\omega(u,x)$, $x\in  D_{\omega,u}$.
A random product of expanding maps above $(\Omega,\tau)$ is then given by
\begin{equation}\label{eq:randomprodexpmaps}
T^{(n)}_{\omega,u}:= T_{\tau^{n-1}\omega,u}\circ\dots\circ T_{\omega,u}: D^{(n)}_{\omega,u} \rr M
\end{equation}
with domain of definition given by $D^{(n)}_{\omega,u}=T_{\omega,u}^{-1} \circ\cdots\circ T_{\tau^{n-1}\omega,u}^{-1} (M)$. 

We associate to $\omega\in\Omega$ the transfer operator 
$\hL_{\omega}:=
\hL_{T_{\omega,u},g_{\omega,u}} 
\in L(c^{\unds}(M); c^{\unds}(U\times M))$,  $0<s\leq r_0$ 
and parametrized by $u\in U$.
For fixed $\phi \in c^{\unds}(M)$ we see by Lemma \ref{lemme:continSOT}
that the map
$\omega \in \Omega \mapsto 
\hL_{\omega} \phi \in c^\unds(U\times M)$ 
as a composition of a continuous and a Bochner measurable map
is Bochner measurable.
Thus, the map 
$\omega\mapsto \hL_\omega  \in L(c^\unds(M) ; c^\unds(\U\times M))$
is strongly Bochner measurable 
(in the sense of Definition \ref{def:Bmeas}).
\\Let
$X_s(M)=L^\infty(\Omega;c^\unds(M))$ and 
$X_s(\U\times M)=L^\infty(\Omega;c^\unds(\U\times M))$ 
denote the
$\Omega$-Bochner measurable sections 
of $c^\unds$-functions on $M$,
respectively $\U\times M$. The operator 
$\hL_{T,g} : c^\undr(M)\rr c^\undr(U\times M)$ has
norm bounded by $C_r$, a constant that only depends upon the
constants of $\P$ and the manifold. Our uniform choice of constants
implies that $\hL_\omega$
induces a well-defined, bounded operator 
$\hBL_s : X_s(M) \rr X_s(\U\times M)$
by declaring that for $\Bphi=(\phi_\omega)_{\omega\in \Omega}\in X_s(M)$:
\begin{equation*}
(\hBL_s \Bphi)_\omega (u,y) = 
(\hL_{\tau^{-1} \omega} \phi_{\tau^{-1}\omega}) (y),
  \ \ \ \omega\in \Omega, u\in U, y\in M.
\end{equation*}
This time we also make explicit in the notation
the dependency upon $s\in (0,r_0]$.
We may now
extract the variable $u$
as a parameter 
to obtain a parametrized 
operator co-cycle over $\Omega$:
\begin{equation*}
u\in U \mapsto \BL_{s,u} \in L(X_s(M)) \ , \ \ \
(\BL_{s,u} \Bphi)_{\omega}(y) 
= ( \hBL_s \Bphi)_\omega (u,y), \ \  \omega\in \Omega, \ \ \Bphi \in X_s(M).
\end{equation*}
The key-point of this section is that if we pre- or post-compose
with a down-grading operator we recover regularity 
from the above parameter extraction:
\begin{theorem}\label{thm:regoftransferop}
	For any $0<s < r\leq r_0$ and $t\in [0,r-s)$,
	the map 
	$u \in U \mapsto 
	\Bj_{s,r}\BL_{r,u}
	=\BL_{s,u} \Bj_{s,r}
	: X_{r}(M) \rr X_{s}(M)$ is 
	$C^{t}$.
\end{theorem}
\Proof  The core of the proof is
Proposition 
\ref{prop:Bparamextr} in the appendix, except that it only
deals with convex subset in Banach spaces.
In view of the description in remark \ref{Cs-structure}
of our $c^\undr(\U\times M)$-structure 
and the parameter extraction Theorem \ref{lemma:injection}
we have, however, the following injections:
\begin{align*}
c^\undr(\U\times M) 
\hrr \prod_{j\in\J} c^\undr(\U\times V_j) 
\hrr \prod_{j\in\J} C^{r-s} (\U;  c^\unds( V_j)) 
\simeq C^{r-s} (\U ; \prod_{j\in\J} c^\unds( V_j)).
\end{align*}
The identification in Remark \ref{Cs-structure} carries over to the
last expression and shows that
$c^\undr(\U\times M) 
\hrr C^{r-s} (\U;  c^\unds(M))$.
With similar identifications
for $X_r(M) = L^\infty(\Omega; c^\undr(M))$ and using 
Proposition \ref{prop:Bparamextr} 
we have 
$X_r(\U\times M) \hrr C^{r-s}(\U ; X_s(M))$.
Thus, $\hBL : X_r(M) \rr
X_r(\U\times M) \hrr C^{r-s}(\U ; X_s(M))$
and we infer that
$u\in\U \rr \Bj_{s,r} \BL_{r,u} : X_r(M) \rr X_s(M)$ is $C^t$
for $0\leq t< r-s$ as wanted.
\Halmos

The above result is precisely the regularity-hypothesis needed 
in Theorem \ref{thm:field}. What is lacking in order to
apply that theorem is to show that 
our operator cocycle is contracting a cone-field in a uniform
way. This will be the subject of the next section.

\subsection{A family of cones adapted to expanding systems}
\label{subsec:cones}

Our setup is such that we are able to use cone-families that do
not depend upon the parameter
$u\in U$ or the random variable $\omega\in \Omega$.
So in the following, and to simplify notation,
we will omit any reference to $u$ and $\omega$.\\

Let $r>0$ and write as usual $\undr=(k,\alpha)$
with $k\in\N$ and $\alpha\in(0,1]$.
Recall that a real-valued function 
$\phi\in c^\undr(M;\R)$ induces $c^\undr$ maps on our charts 
$\phi_j = \phi \circ \eta_j^{-1} : V_j\rr \R$, $j\in\J$.
We will construct a family of cones adapted to this $c^\undr$-structure.
Recall that when $\xi,\xi'\in\V_j$ we let
$d_j^g(x,x')=d_M(\xi,\xi')$ denote the induced Riemannian distance
between $x=\eta_j(\xi),x'=\eta_j(\xi')\in V_j$.
Also let $g^j_x$ be the induced Riemannian metric tensor in the chart $V_j$.
The length of $v\in \R^n \simeq T_x V_j$ for this metric is:
\begin{equation*}
\|v\|_{g^j_x} = \sqrt{g^j_x(v,v)} .
\end{equation*}
By Hypothesis \ref{hypAtlas} this norm is $L$-equivalent to
the Euclidean norm on $\R^n$.

Given a $k$-linear form  $A$ on $\R^n$ and a base point $x\in V_j$
we define its $g$-norm at $x$ to be:
\begin{equation*}
\|A\|_{g_x^j} = \sup \left\{ |A.(v_1,\ldots,v_k)| \ : \ \ 
\forall 1\leq i\leq k :
v_i\in \R^n, \|v_i\|_{g^j_x} \leq 1\right\} .
\end{equation*}
With this definition we have in particular:
\begin{equation}
\frac{1}{L^q}\|D^q \phi_j(x)\| \leq \|D^q \phi_j(x)\|_{g_x^j}
\leq L^q \|D^q \phi_j(x)\|
\label{equiv:seminorms}
\end{equation}

Given a vector of strictly positive numbers
$\vec{a}:=(a_1,\dots,a_k,a_{k,\alpha})\in(\mathbb{R}_+^*)^{k+1}$,
we consider the convex cone $\C_{\vec{a}}$ 
of non-negative functions
$\phi\in c^\undr(M)$, $\phi\geq 0$ satisfying the following conditions:
\begin{itemize}
	\item [(C1)] For $j\in\J$, 
	$q\in\{1,\dots,k\}$ and $x\in V_j$: 
	$\|D^q\phi_j(x)\|_{g^j_x}\leq a_q\phi_j(x)$. 
	\item [(C2)] For $j\in\J$ and $x,x'\in V_j$
	with $0<d_j^g(x,x')\leq \delta_1$:
	\[
	\frac{\|D^k\phi_j(x)-D^k\phi_j(x')\|_{g^j_x}}{d_j^g(x,x')^\alpha} 
	\leq a_{k,\alpha}
	\phi_j(x).\]
\end{itemize}

\begin{rem}
	Similar cones
	were considered in \cite{BKS96}, but the
	treatment in the case of a H\"older condition on derivatives
	was somewhat incomplete. 
	%(see the remark right before \cite[Lemma 3.2]{BKS96}).
\end{rem}

The semi-norms used to define the cone $\Cveca$ are tailored to
prove contraction properties for the Hilbert-metric. 
We also have to prove regularity of the cone and for that we
need the following lemma to connect norms in the Banach spaces and
the semi-norms used in the above definition.
For use elsewhere in the article
it is slightly more general than needed here.

\begin{lemme}
	\label{holdereq}
	Let $K,L\in [1,+\infty)$ and let
	$W\subset X$ be an open convex set in a Banach space.
	Let $d_1$ be the norm-metric on $W$ and let
	$d_2$ be another metric
	which is $L$-equivalent with $d_1$, i.e.\
	$\frac{1}{L} d_1 \leq d_2 \leq L d_1$ 
	as functions on $W\times W$. 
	Also let $\delta_1,\delta_2>0$ be $K$-equivalent constants.
	Then $E=C^{(k,\alpha)}_{\delta_1,d_1}(W ; Z)
	= C^{(k,\alpha)}_{\delta_2,d_2}(W ; Z)$  (with $Z$ any Banach space) and 
	for any  $\phi\in E$ the H\"older-exponents
	$h^{(k,\alpha)}_{\delta_1,d_1} (\phi)$ and
	$\myh^{(k,\alpha)}_{\delta_2,d_2} (\phi)$
	are  $\lceil LK \rceil^{1-\alpha} L^\alpha$-equivalent.
\end{lemme}
\Proof It is enough to show the lemma for $k=0$. 
\\Let $\xi,\xi'\in W$ with $\|\xi-\xi'\|\leq \delta_1$.
Let $N\geq 1$ be an integer.
The points  $\xi_\ell=\xi + \frac{\ell}{N} (\xi'-\xi)$
$0\leq \ell\leq n$ belong to $W$ (convex). So
\begin{equation*}
d_2(\xi_\ell,\xi_{\ell+1}) \leq
\frac{L}{N}  \| \xi - \xi'\| \leq 
\frac{L\delta_1}{N}  \leq
\frac{LK}{N}  \delta_2.
\end{equation*}
When $N=\lceil LK \rceil$  then $d_2(\xi_\ell,\xi_{\ell+1})\leq \delta_2$ so
by the $\delta_2,d_2$-H\"older estimate:
\begin{equation*}
|\phi(\xi_\ell) - \phi(\xi_{\ell+1})| \leq C_2 d_2(\xi_\ell,\xi_{\ell+1})^\alpha
\leq C_2 \left( \frac{L}{N} \right)^\alpha \|\xi-\xi'\|^\alpha .
\end{equation*}
Finally, $
|\phi(\xi) - \phi(\xi')| 
\leq C_2 N^{1-\alpha} L^\alpha \|\xi-\xi'\|^\alpha
\leq C_2 \lceil LK\rceil^{1-\alpha} L^\alpha d_1(\xi,\xi)^\alpha$,
implying one inequality. Interchanging $d_1$ and $d_2$ turns
out to yield a bound
with the same constant (even though the problem is not symmetric).
\Halmos\\

The above lemma shows that H\"older constants evaluated for the
Euclidean metric and for the induced manifold metric are equivalent
within uniform constants. It follows that our usual chart-norm  
(\ref{def:chartCka}) is equivalent to:
$\|\phi\|_g:= \max_{0\leq q\leq k} \sup_x \|D^q \phi_j(x)\|_{g_x^j} \vee
h_{\delta_1,d_j^g}^{(k,\alpha)}(\phi_j)$.
For the latter norm 
it is clear from the construction
that $\C_{\vec{a}}$ is a closed, convex cone
but this is then also the case in our usual chart-norm.
\medskip

\begin{lemme}
	\label{lem:bdratio}
	For every $\vec{a}$ there is a constant
	$R>0$ such that for $\phi\in \C_{\vec{a}}$: 
	\begin{equation*}
	\phi(x)\leq R \phi(y), \ \ \forall  x,y\in\M.
	\end{equation*}
	%Let $D$ denote the diameter of $M$. For $k\geq 1$ we may take 
	%for the value of this constant: $R=\exp(\frac{4}{3} a_1 D)$,
	%while if $k=0$ we may take:
	%$R = (1 + a_{(0,\alpha)})^N$ with $N = \frac{4D}{3\delta_1}+1$.
\end{lemme}
\Proof Consider $k\geq 1$ and $j\in\J$.
By Hypothesis \ref{hypAtlas} on the atlas
we have $\|D\phi_j(x)\|_{g_x^j}\leq a_1  \phi_j(x)$ 
for all $x\in V_j$. 
So if a geodesic lies in $\V_j$ and is joining
$y,y'\in \V_j$ we have
\begin{equation*}
\phi(y')\leq e^{a_1  d_M(y,y')} \phi(y).
\end{equation*}
Since the $d_M$-distance is realized by geodesics, we deduce that for
every $y,y'\in M$: $\phi(y') \leq e^{a_1 D} \phi(y)$
with $D$ being the diameter of $M$.

For $k=0$, the condition (C2) yields: 
$\phi(x') \leq (1+a_{0,\alpha} d_j^g(x,x'))^\alpha) \phi(x)$ whenever
$d_j^g(x,x')\leq \delta_1$. Given arbitrary $y,y'\in M$ we consider
a geodesic between the two points and chop it into $N$ pieces with $N=\lfloor \frac{D}{\delta_1}\rfloor+1$. Each piece has length at most $\delta_1$. For endpoints of such a piece we may thus apply the previous estimate showing the wanted bound with
$R=(1+a_{0,\alpha} \delta_1)^N$. \Halmos\\

\begin{lemme}
	\label{lem:Ca-inner-outer}
	Every cone $\C_{\vec{a}}$ is inner and outer regular.
	Furthermore, if $\sigma\in(0,1)$ then
	\begin{equation*}
	\diam_{\C_{\vec{a}}} (\C_{\sigma \vec{a}})^* 
	\leq 2 \log\left(\frac{1+\sigma}{1-\sigma}\right) + 2 \log R.
	\end{equation*}
\end{lemme}

\Proof The constant function $\bfone\in\C_{\vec{a}}$. 
Let $\phi=\bfone+\delta\psi\in c^\undr(M)$
with $\|\delta\psi\|_r\leq c$.
Then by the equivalence of seminorms \eqref{equiv:seminorms}, we have
for $j\in\J$, $x\in V_j$, $1\leq q\leq k$:
\[\|D^q\phi(x)\|_{g_x^j} \leq L^q  c \leq a_q (1-c) \leq  a_q\phi(x)\]
provided $c \leq \frac{a_q}{L^q+a_q}$. 
\\Similarly for the H\"older-condition using Lemma \ref{holdereq}.
Thus, $B_{c^{\undr}}(\bfone,\rho)\subset \C_{\vec{a}}$ whenever
$\ds 0<\rho< \min_{1\leq q\leq k} \frac{a_q}{1+a_q} \wedge 
\frac{a_{k,\alpha}}{1+a_{k,\alpha}}$,
showing inner regularity.\\

For outer regularity, one may
take $\ell$ to be any probability measure on $M$,
e.g.\ evaluation at a point. Then for $\phi\in \C_{\vec{a}}$ we have
by the previous lemma
$\phi(x) \leq R \langle \ell,\phi\rangle$ for all $x\in M$, whence
\begin{equation*}
\|\phi\|_r \leq 
K \left(\max_{1\leq q\leq k} \frac{a_q}{1+a_q} 
\vee \frac{a_{k,\alpha}}{1+a_{k,\alpha}}\right)
\sup_x \phi(x) \leq
R \left( \max_{1\leq q\leq k} \frac{a_q}{1+a_q} 
\vee \frac{a_{k,\alpha}}{1+a_{k,\alpha}} \right) \langle \ell,\phi\rangle
\end{equation*}
again with a constant $K$ depending upon the equivalence-constants of
semi-norms. This shows outer regularity.\\

In order to show that the inclusion is uniformly bounded
we will describe a 
general (and fairly standard)  way to obtain such bounds. Note that
in both (C1) and (C2) the left hand-side of the two inequalities
is a semi-norm
on $c^\undr(M)$. So (C1) and (C2)
may be formulated as a (in fact, uncountable) collection of elementary
conditions of the form:
\begin{equation}  p(\phi)\  \leq  \ a \; \phi(x)
\label{eq:simplecond}\end{equation}
where $a$ is one of the cone-constants, $x\in M$ and
$p$ is a semi-norm  on $c^\undr(M)$.

Now consider $\phi_1,\phi_2\in\C_{\sigma\vec{a}}^*$. In order to
calculate the projective Hilbert-distance of the two functions
in $\C_{\vec{a}}$ we need to calculate the infimum of $t,t'>0$
so that $t\phi_1-\phi_2,t'\phi_2-\phi_1\in \C_{\vec{a}}$.
However, when $\phi_1,\phi_2$ verifies the elementary
condition (\ref{eq:simplecond})
with $\sigma a$ instead of $a$ we have the following sequence
of inequalities:
\begin{align}
p( t\phi_1-\phi_2 )
& \leq  t 
p( \phi_1)  +
p( \phi_2)  \nonumber  \\
& \leq  \sigma a (t \phi_1(x) + \phi_2(x)) \nonumber \\
& \leq  a (  t \phi_1(x) - \phi_2(x)) \label{eq:tvalue}
\end{align}
where equality in \eqref{eq:tvalue} is satisfied as soon as we set:
\begin{equation*}
t =
\frac{1+\sigma}{1-\sigma} \sup_{x\in M} 
\frac{\phi_2(x)}{\phi_1(x)}.
\end{equation*}
Here, the semi-norm has completely disappeared in the last condition.
With a similar calculation for $t'$ we get for the Hilbert-distance:
\begin{equation*}
\exp (d_{\C_{\vec{a}}} (\phi_1,\phi_2)) \leq
tt' =
\left(\frac{1+\sigma}{1-\sigma} \right)^2 
\sup_{x,x'\in M} 
\frac{\phi_2(x)}{\phi_1(x)}.
\frac{\phi_1(x')}{\phi_2(x')}
\leq
\left(\frac{1+\sigma}{1-\sigma} R \right)^2 .
\end{equation*}
Taking log we get the wanted estimate.\Halmos\\

We claim that the (weighted) transfer operator 
is a \textbf{strict and uniform contraction} 
of the cones $\C_{\vec{a}}$, in the following sense:
\medskip

\begin{theorem}[\cite{BKS96} Lemma 3.2]
	\label{unifconecontraction}
	Let $r>0$ and $\undr=(k,\alpha)$ as above. If 
	$k=0$ take $\sigma \in (\frac{1}{\lambda^\alpha}, 1)$ and if
	$k\geq 1$ take instead $\sigma\in (\frac{1}{\lambda},1)$ with
	$\lambda>1$
	being the expansion constant  for our expanding maps.
	There exists a vector $\vec{a}=(a_1,...,a_k,a_{k,\alpha})
	\in (\R^*_+)^{k+1}$ depending only upon the constants of 
	$\P=\P^{(r)}_{k_T,k_g,k_\Sigma}$ 
	and the manifold such that for every pair $(T,g)\in\P$ 
	the  transfer operator  $\L_{T,g} \in L(c^\undr(M,\mathbb R))$, defined for $\phi\in c^\undr(M,\mathbb R)$ by
	\begin{equation*}
	   %\label{transferopexpmap}
	\L_{T,g}\phi(x)=\sum_{y,Ty=x}e^{g(y)}\phi(y),\ \ \  x\in M, 
	\end{equation*}
	satisfies:
	$\L_{T,g}(\C_{\vec{a}})\subset \C_{\sigma\vec{a}}$. 
	\\In particular,
	by the previous lemma, $\L_{T,g}$ is a uniform contraction of $\C_{\vec{a}}$
	for the Hilbert metric.
	Furthermore, there exist $\rho>0$, depending only on
	the constants of the problem
	such that for every $\phi\in(\C_{\vec{a}})^*$
	\begin{equation*}
	B(\L_{T,g}\phi,\rho\|\L_{T,g}\phi\|)\subset \C_{\vec{a}}.
	\end{equation*}	
\end{theorem}

\paragraph{Proof:}
(Our proof follows the same lines
as in the proof of \cite{BKS96} Lemma 3.2. Here, we provide
details of the uniform bounds also in the presence of a H\"older condition
on the top-most derivative).
Write $\L:=\L_{T,g}$. 
For the first statement, 
take a chart $V_j$, $j\in\J$.
The transfer operator takes the form  as in
\eqref{eq:Linchartj} (where as above we omit the $u$-variable):
\begin{equation*}
(\L_{j} \phi)(y)
= \sum_{i\in\A}
e^{g^i_j(y)} \phi_{m_i} \circ \psi_{m_i,j}^i(y), \ \ y\in V_j
\end{equation*}
When $\phi\in \C_{\vec{a}}$ we will determine 
suitable values of $a_1,...,a_k$ by induction starting with $a_1$.
\\Assume that for $1\leq q<k$  we have determined $a_1,...,a_{q-1}$.
We then get by the Fa\`a di Bruno formula and algebra for the $q$'th derivative:
\begin{equation*}
D^q (\L_j \phi)(y) =
\sum_{i\in\A} e^{g^i_j} 
\left[ D^q\phi_{m_i}\circ \psi_{m_i,j}^i (y).
(D\psi^i_{m_i,j},\ldots,D\psi^i_{m_i,j}) + R_{q-1}(y) \right]
\label{eq:Lphideriv}
\end{equation*}
where $R_{q-1}$ is linear in $\phi,...,D^{q-1}\phi$ and multinomial
in derivatives up to order $q$ of $\psi^i_{m_i,j}$ and $g_j^i$.
Note that the only term with a $q$'th derivative of $\phi$
has the form $D^q\phi\circ\psi(y).(D\psi(y),\ldots,D\psi(y))$. 
Since the $\psi$'s are $1/\lambda$-contractions in the Riemannian metric,
the $\|\cdot\|_{g_y^j}$-norm of this term is
bounded by $\frac{a_{q}}{\lambda^q}  \phi \circ \psi(y)$.
Inserting in the above expansion, using
the induction hypothesis and the telescopic principle \eqref{eq:telescopic} we get:
\begin{equation*}
\|D^q \L_j \phi(y)\|_{g^j_y} \leq 
\left( \frac{a_q}{\lambda^q} + \sum_{\ell=1}^{q-1}
a_\ell  \ p_{q,\ell}(k_T,k_g) 
\ + \ p_{q,0}(k_T,k_g)
\right) \L_j\phi(y), \ \ y\in V_j, 
\end{equation*}
where $p_{q,\ell}$, $0\leq \ell\leq q-1$
are polynomials in the constants (except for $k_\Sigma$) of the problem.
Given $\frac{1}{\lambda}<\sigma<1$ 
and having chosen $a_1,...,a_{q-1}$ it then suffices to choose
$a_q$ large enough to ensure that
\begin{equation*}
\|D^q \L_j \phi(y)\|_{g^j_y} \leq 
\sigma {a_q}
\L_j\phi(y), \ \ y\in V_j.
\end{equation*}
For the H\"older-condition  on the $k$'th derivative consider
$y,y'\in V_j$, $0< d^g_j(y,y')\leq \delta_1$.
In order to achieve cone-contraction for
the H\"older condition we notice that in the difference
$D^k \L_j \phi(y')- D^k \L_j \phi(y)$ there is only one term 
which involves the H\"older bound for the $k$'th derivative
of $\phi$. It takes the following form (we may treat the other terms
using telescoping):
\begin{equation*}
e^{g_j^i(y)}\left( D^q\phi \circ\psi(y') - D^q\phi \circ\psi(y) \right)  
(D\psi(y),\ldots,D\psi(y)).
\end{equation*}
Now, this being taken at $y$, the $\|\cdot\|_{g_y^j}$-norm
is bounded by 
$e^{g_j^i(y)} \frac{a_{k,\alpha}}{\lambda^{k+\alpha}}
d^g_j(y,y')^\alpha \phi\circ \psi(y)$
%+ D^q\phi(\psi(y')).
%\left(
%(D\psi(y'),\ldots,D\psi(y'))-
%(D\psi(y),\ldots,D\psi(y))
%\right)  .
%\end{equation*}
and we obtain
as in the previous calculation:
\begin{equation*}
\frac{\|D^k\L_j\phi(y')-D^k \L_j \phi(y)\|_{g^j_y}}
{d_j^g(y,y')^\alpha} \leq 
\left( \frac{a_{k,\alpha}}{\lambda^{k+\alpha}} + \sum_{\ell=1}^{k}
a_\ell  \ p_{k,\ell}(k_T,k_g) 
\ + \ p_{k,0}(k_T,k_g)
\right) \L_j\phi(y), \ \ y\in V_j, 
\end{equation*}
This time, however, we have to distinguish if $k\geq 1$ in which case
we may choose $\sigma$ as above, and the case $k=0$ for which we
must choose $\sigma\in(\frac{1}{\lambda^\alpha},1)$ in order
to assure that
\begin{equation*}
\frac{\|D^q\L_j\phi(y')-D^q \L_j \phi(y)\|_{g_y^j}}{d_j^g(y,y')^\alpha} \leq 
\sigma {a_q}
\L_j\phi(y), \ \ y\neq y'\in V_j.
\end{equation*}
This proves the uniform cone-contraction. The last statement follows from Lemma \ref{lem:uniopen}.\Halmos

\subsection{Random physical measure and quenched linear response}\label{subsec:linresp}

We assume here that the metric tensor is uniformly
$C^\infty$ (cf. Hypothesis 
	\ref{hypAtlas} for the meaning of being uniform)
in all charts
and write $m$ for the associated volume form. 
In a local chart $V_j$, $j\in \J$, we have:
$dm_x = \rho_j(x) \; d^n x := \sqrt{\det(g_x^j}) \; d^n x$, $x\in V_j$
with $\rho_j$ being uniformly 
$C^\infty$ and uniformly bounded away from zero.\\

Going back to the setting of a random product \eqref{eq:randomprodexpmaps}
of $c^{r}$
uniformly expanding maps (section \ref{sec:expmaps}),
we consider $(T_{\omega},g_{\omega})_{\omega\in\Omega}\in\P_{k_T,k_g,k_\Sigma}$ where
the exponential of the weight is defined as the metric derivative of the map:
$e^{g_\omega} :=1/ \left| \det_g(\partial_x T_\omega) \right|$
(beware of the two different uses of the letter $g$).
In order to get physical measures we must further assume that the 
domain $D_{\omega,u}$ for each $T_{\omega,u}=T_\omega(u,\cdot)$ 
has full measure in $M$.

The regularity condition on $g_\omega$
implicitly impose a regularity condition on $T_w$
(or rather its inverse branches).
Omitting $\omega$ in the notation, let $j\in\J$
and $\psi_j^i : U\times V_j \rr M$, $i\in\cal A$
be a local inverse branch. Take a point $\xi_0=(u_0,y_0)\in V_j$ and
pick $m_i\in \J(\psi_j^i(\xi_0))$. Restricting to $B(\xi_0)\subset U\times V_j$
we get in local coordinates:
\begin{equation*}
g_\omega:=
\log  \frac{\rho_j(u,y)}{\rho_{m_i}(\psi_j^i(u,y))} 
- \log \left|\det \partial_y \psi_j^i(u,y)\right|, \ \ \
(u,y)\in B(\xi_0).
\end{equation*}
This should be uniformly bounded in $C^{r-1}$ norm.
Here,  as the $\rho$'s are uniformly $C^\infty$ and uniformly bounded from
above and below, the first term is bounded as soon as $\psi_j^i$ is.
For the second term, however, the bound (\ref{eq:condg})
imposes the auxiliary condition
that there should be $k_g'<+\infty$ so that
\begin{equation}
\|D \log  \left| \det \partial_y \psi_j^i \right|
\|_{(r-1,\alpha)} \leq k_g',
\ \ \ \ \forall i\in \calA, j\in \J.
\label{aux:cond}
\end{equation}
By change of variables (and the hypothesis that
$D_{\omega,u}$ has full measure in $M$):
\begin{equation*}
\int_M \bfone  \;dm_x = \int_M S_{\omega,u}(y)  \;dm_y
\end{equation*}
where
$S(y)=S_{\omega,u}(y) :=\sum_{x : T_{\omega,u}=y} e^{g(u,x)}$.
It follows that $S$ must take values both above and below one
(if it is not identically one).
Now, by condition (\ref{aux:cond}) ratios of $S(y)/S(y')$ are
uniformly bounded from above and below, whence we see that
condition (\ref{eq:condS}) is automatically verified.

Given the above family of maps and weights we may thus associate
as in the previous section
the transfer operator cocycle $\hBL_u$, 
acting on $X_{r-1}(M,\CX)$.
Then  for each $u\in\U$, $\hBL_u$ is a cone-contracting cocycle (in the sense of Section \ref{sec:cone-con}) of the cone $\C_{\vec{a}}$ for some well-chosen $\vec{a}$: this follows straightforwardly from Theorem \ref{unifconecontraction}.
\medskip

In this setting, Theorem \ref{thm:fixedptsection} apply, so we may construct $\field_u\in X_{r-1}(M,\CX)$, fixed point of $\bs\pi_u$. The coincidence\footnote{This is a somewhat particular case. In applications, often the linear form used for the outer regularity of $\C_{\vec{a}}$ and the left eigenvector of $\L_{\omega,u}$ do not coincide.} of the outer regularity form of $\C_{\vec{a}}$ and the left eigenvector of $\L_{\omega,u}$ ensures that 
\begin{itemize}
	\item For almost every $\omega\in\Omega$, every $u\in\U$, \[p_{\omega,u}=\langle\ell,\Lou f_{\omega,u}\rangle =\int_M\Lou f_{\omega,u}dm=\int_M f_{\omega,u}dm=1.\]
	In particular, the fixed point $\bs{\field}_u$ of $\bs{\pi}_u$ is also a fixed point for $\hBL_u$.
	\item The top characteristic exponent of this transfer operator cocycle is zero: by virtue of \eqref{eq:convenientformula}, one has \[\chi_u=\mathbb{E}[\log(\bs{p}_u)]=0.\]
\end{itemize}

By Theorem \ref{thm:regoftransferop}, we may apply Theorem \ref{thm:field}: for any $0\leq s<r-1$ the map $\field:\U\rightarrow X_{s}(M,\CX)$ is $C^t$ for any $0<t<r-1-s$.  For $\psi\in L^1(M)$ we then set
\begin{equation}\label{def:randacim}
\nu_{\omega,u}[\psi]=\int_{M}\psi \;f_{\omega,u}dm,
\end{equation}
which defines an absolutely continuous equivariant probability measure for $T_{\omega,u}$, i.e. for any $\omega\in\Omega$, any $u\in\U$, 
\begin{equation*}\label{randinvcond}
\nu_{\omega,u}[\psi\circ T_{\omega,u}]=\nu_{\tau\omega,u}[\psi].
\end{equation*} 
$\bs{\nu}_u$ is called the \emph{random a.c.i.m}~\footnote{A.c.i.m stands for absolutely continuous invariant measure. Depending on the context, it is also called the random S.R.B measure (for Sinai-Ruelle-Bowen), or the random physical measure.}.
\medskip

We may now formulate the main result of this section: a \emph{linear response formula} for the random a.c.i.m associated to a random product of expanding maps:

\begin{theorem}
\label{quenchedlinrespthm}
	Let $(T_{\omega,u})_{\omega\in\Omega,u\in\U}\in\P_{k_T,k_g,k_\Sigma}$ be a random family of uniformly expanding maps, and let $\hBL_u$ be the transfer operator cocycle it generates above $(\Omega,\tau)$, acting on $X_{r-1}(M,\CX)$.
	Let $\bs{\nu}_u$ be the random a.c.i.m introduced in \eqref{def:randacim}.
	\\For every observable $\psi\in L^1(M)$, the map $u\in\U\mapsto \bs\nu_{u}[\psi]\in L^\infty(\Omega,\CX)$ is differentiable at $u=u_0$, for every $u_0\in\U$ with
	\begin{equation}\label{linearresponseformula1}
	D_u\left[\int_M\psi \, d\nu_{\omega,u}\right]_{u=u_0}=\sum_{n=0}^\infty\int_M\psi\circ T^{(n)}_{\tau^{-n}\omega,u_0}P_{\tau^{-(n+1)}\omega,u_0}f_{\tau^{-(n+1)}\omega,u_0}dm
	\end{equation}
\end{theorem}

For the proof of this result, we will use an estimate on the speed of convergence of the random product $\Lnou$ towards its equivariant line, analogous to the spectral gap estimate for a deterministic cone contraction.

\begin{lemme}\label{speedofconvrandprodlemma}
	Let $(T_{\omega,u})_{\omega\in\Omega,u\in\U}\in\P_{k_T,k_g,k_\Sigma}$ be a random family of expanding maps, with dilation constants all bounded from below by some $\lambda>1$. 
	\\[1mm]
	Let $\hBL_u$ be the associated transfer operator cocycle above $(\Omega,\tau,\mathds{P})$, and let $\field_u\in X_{r-1}(M,\CX)$ be its fixed point.
	Then for every integer $0<s\leq r-1$ and every $\bs\phi\in X_s(M,\CX)$:
	\begin{equation}\label{speedofconvrandprod}
	\|\Lnou\phi_\omega-f_{\tau^n\omega,u}\int_M\phi_\omega dm\|_{C^s}\leq C\eta^{n-1}\|\phi_\omega\|_{C^s}
	\end{equation}
	where $C$ depends only on the cone $\C_{\vec{a}}$ and $\eta<1$.
\end{lemme}

\paragraph{Proof:} Let $n\geq 1$. We write \eqref{eq:majcontraction}, at $\bs\phi,f_{\tau^{-n}\omega}$, for some Bochner-measurable, essentially bounded family $(\phi_\omega)_{\omega\in\Omega}\in\C_{\vec{a}}\subset c^s(M)$, to get

\begin{equation*}
\left\|\L_{\tau^{-n}\omega,u}^{(n)}f_{\tau^{-n}\omega}-\frac{\L_{\tau^{-n}\omega,u}^{(n)}\phi_{\tau^{-n}\omega}}{\int_M\phi_\omega dm}\right\|_{C^s}\leq \frac{K\Delta}{2}\eta^{n-1}
\end{equation*}

with $K$ the sectional aperture of $\C_{\vec{a}}$, $\Delta=diam_{C_{\vec{a}}}\C_{\sigma\vec{a}}<\infty$, $0<\sigma<1$ and $\vec{a}$ given by Theorem \ref{unifconecontraction} and $\eta=\tanh\left(\frac{\Delta}{4}\right)<1$.
\\One obtains
\begin{equation*}
\left\|f_{\omega,u}-\frac{\L_{\tau^{-n}\omega,u}^{(n)}\phi_{\tau^{-n}\omega}}{\int_M\phi_\omega dm}\right\|_{C^s}\leq \frac{K\Delta}{2}\eta^{n-1}
\end{equation*}
which yields, with the "change of variables" $\omega\leftrightarrow\tau^{-n}\omega$ and once multiplied by $0<\int_M\phi_\omega dm\leq\|\phi_\omega\|_{C^s}$, \eqref{speedofconvrandprod}
\medskip

In the general case take $\phi\in c^s(M)$ of norm one.
As shown in Lemma \ref{lem:Ca-inner-outer} the cone
$\C_{\vec{a}}$ is
inner regular. In fact, as shown in the proof of that Lemma
there is
$\rho=\rho_{\vec{a}}>0$ so that
$B(\bfone ,\rho_{\vec{a}})\subset \C_{\vec{a}}$ (with $\bfone$ being
the constant function).
Then 
$\phi=
\frac{1}{2\rho} (\bfone+\rho \phi) -
\frac{1}{2\rho} (\bfone-\rho \phi)$ shows that we have a decomposition
$\phi=\phi_1-\phi_2$,
with $\left(\phi_1,-\phi_2\right)\in(\C_{\vec{a}})^2$ and 
$\|\phi_1\|_{C^s}+\|\phi_2\|_{C^s}\leq \frac{2}{\rho}\|\phi\|_{C^s}$. Thus,

\begin{align*}
\notag\left\|f_{\tau^n\omega,u}\int_M\phi_{\tau^n\omega}dm-\L_{\omega,u}^{(n)}\phi_\omega\right\|_{C^s}&= \left\|f_{\tau^n\omega,u}\int_M
(\phi_{1,\tau^n\omega}-\phi_{2,\tau^n\omega})dm-\L_{\omega,u}^{(n)}(\phi_{1,\omega}-\phi_{2,\omega})\right\|_{C^s}\\
&\leq \frac{K\Delta}{2}\eta^{n-1}\left[\|\phi_{1,\omega}\|_{C^s}+\|\phi_{2,\omega}\|_{C^s}\right]\\
&\leq \frac{K\Delta}{\rho}\eta^{n-1}\|\phi_\omega\|_{C^s}
\end{align*}

which conclude the proof in the general case $\phi\in c^s(M)$.\Halmos

\begin{rem}\label{limitwelldefined}
	\begin{itemize}
		\item From estimate \eqref{speedofconvrandprod}, we may draw the following conclusion: if $(\phi_\omega)_{\omega\in\Omega}\in c^s(M)$ is such that $\int_M\phi_\omega dm=0$, then for any $n\geq 1$, \[\|\L_{\tau^{-n}\omega,u}^{(n)}\phi_\omega\|_{C^s}\leq \frac{K\Delta}{\rho}\eta^{n-1}\|\phi_\omega\|_{C^s}.\]
		In particular, the limit $\sum_{n=0}^\infty \L_{\tau^{-n}\omega,u}^{(n)}\phi_\omega$ is well defined whenever $\int_M \phi_\omega dm=0$.
		\item Estimate \eqref{speedofconvrandprod} has far reaching consequences: in particular, it can be used to establish \emph{exponential decay of random correlations}, in the same way as one establishes exponential decay of correlations in the deterministic case. We refer to \cite[Thm B]{BKS96} for more details on this.
		\item It is possible to establish uniform exponential decay of random correlations for general weighted cocycles of transfer operators $(\L_{T_\omega,g_\omega})_{\omega\in\Omega}$: see \cite[Theorem 3.1]{Dubois08}. 
	\end{itemize}
\end{rem}

\paragraph{Proof of theorem \ref{quenchedlinrespthm}} 
For any $0\leq s< r-1$, it follows from Theorem \ref{thm:field} that the map $u\in\U\mapsto\field_u\in X_{s-1}(M)$ is differentiable with
\begin{equation}\label{invdensitydiff}
D_u\field_u=\left[\Id-Q_{u}(\field_u)\right]^{-1}P_u(\field_u),
\end{equation}
with $Q_u(\field_u)$, $P_u(\field_u)$ coming \eqref{eqforQ} and \eqref{eqforP}. In this particular case, one has 
\begin{align*}
(Q_{u}.)_\omega &=\L_{\tau^{-1}\omega,u}-\left(\int_M.dm\right)f_{\omega,u}\\
[P_u]_\omega(f_{\omega,u})&=(\partial_u\L_{\tau^{-1}\omega,u})f_{\tau^{-1}\omega,u}-\left(\int_M (\partial_u\L_{\tau^{-1}\omega,u}) f_{\tau^{-1}\omega,u}dm\right)f_{\omega,u}.
\end{align*}

The normalization $\int_M \field_udm =1$ implies that $\int_M D_u\field_u dm =0$. Similarly, $\int_M\hBL_u\field_udm=1$ yields that
\[0=D_u[\int_M\Lou f_{\omega,u}]=\int_M \partial_u\L_{\omega,u}f_{\omega,u}dm+\int_M\Lou D_uf_{\omega,u}dm= \int_M \partial_u\L_{\omega,u}f_{\omega,u}dm. \]
Hence, one obtains $Q_{\omega,u}D_uf_{\omega,u}=\L_{\tau^{-1}\omega,u}f_{\tau^{-1}\omega,u}$ and $(P_u)_\omega(f_{\omega,u})=\partial_u\L_{\tau^{-1}\omega,u}f_{\tau^{-1}\omega,u}.$
\medskip

A straightforward computation shows that \[(Q_u^n)_\omega=\L^{(n)}_{\tau^{-n}\omega,u}-\left(\int_M.dm\right)f_{\omega,u}\]
and thus by \eqref{speedofconvrandprod}, for any $0\leq s<r-1$, the sum $\sum_{n=0}^\infty Q_u^n$ converges in $X_s(M,\CX)$. Hence for any $\bs\phi\in X_{s}(M,\CX)$, one may define a bounded operator via
\begin{equation*}
\left(\left[\Id-Q_{u}\right]^{-1}\phi\right)_\omega:=\sum_{n=0}^\infty (Q_u^n)_\omega\phi_{\tau^{-n}\omega}.
\end{equation*}
\medskip

Let $u_0\in\U$. From remark \ref{limitwelldefined} and the previous discussion, one sees that \eqref{invdensitydiff} can be rewritten as

\begin{equation*}\label{linearresponseformula2}
D_uf_{\omega,u_0}=\sum_{n=0}^\infty\L^{(n)}_{\tau^{-n}\omega,u_0}(\partial_u\L)_{\tau^{-(n+1)}\omega,u_0}f_{\tau^{-(n+1)}\omega,u_0},
\end{equation*}
which is valid in $X_{s-1}(M,\CX)$. This last equation yields

\begin{align*}
D_u\nu_{\omega,u_0}[\psi]=\int_M\psi D_uf_{\omega,u_0}dm&=\int_M\sum_{n=0}^\infty\psi\L^{(n)}_{\tau^{-n}\omega,u_0}\partial_u\L_{\tau^{-(n+1)}\omega,u_0}f_{\tau^{-(n+1)}\omega,u_0}dm,
\end{align*}

which yields \eqref{linearresponseformula1} by using boundedness of the integral on $C^s(M)$ and the duality property of the transfer operator $\int_M\psi\Lou\phi dm=\int_M\psi\circ T_{\omega,u}\phi dm$.\Halmos
\medskip

Now that we have established the regularity of the map $u\in\U\mapsto\int_M \psi d\nu_{\omega,u}$ as an element of $L^\infty(\Omega,\mathbb C)$, we obtain for free the regularity of $u\in\U\mapsto R_{\psi}(u):=\mathds{E}(\int_M \psi d\nu_{\omega,u})$, i.e an annealed version of theorem \ref{quenchedlinrespthm}: 
\begin{theorem}
      \label{annealedlinrespthm}
Let $(T_{\omega,u})_{\omega\in\Omega,u\in\U}\in\P_{k_T,k_g,k_\Sigma}$
be a random family of uniformly expanding maps, let $(\hBL_u)_{u\in\U}$
be the transfer operator cocycle it generates above $(\Omega,\tau)$, acting on $X_{r-1}(M,\CX)$, and $\bs{\nu}_u$ be the associated random a.c.i.m.
	\\For every observable $\psi\in L^1(M)$, the map $u\in\U\mapsto \mathds{E}\left[\int_M\psi d\bs{\nu}_u\right]$ is differentiable at $u=u_0$, for every $u_0\in\U$ with
	\begin{equation}\label{linearresponseformula3}
	D_u\left[\mathds{E}\left(\int_M\psi d\bs{\nu}_{u}\right)\right]_{u=u_0}=\sum_{n=0}^\infty\int_{\Omega}\int_M\psi\circ T^{(n)}_{\omega,u_0}\partial_u\L_{\tau^{-1}\omega,u_0}f_{\tau^{-1}\omega,u_0}dmd\mathds{P}.
	\end{equation}
\end{theorem}

\subsection{Proof of Theorem \ref{exquenchedlinresponse} (for
	Examples \ref{ex:intro} and \ref{ex:introtwo})}
\label{sub:examples}

Consider the type of map used in 
Example \ref{ex:intro} in the introduction
  of the form:
$T_u(x)=M x + \kappa(u,x) \, \mod \Z^d$,
where for simplicity we have omittend
reference to $\omega\in \Omega$ in the notation.
We want to show that for  a measurable collection of such
maps we may apply
Theorem \ref{quenchedlinrespthm} and
Theorem \ref{annealedlinrespthm}. For this we need to show
that $T_u$ is uniformly $\lambda$-expanding and 
that the map together
with the natural weight
verify the uniform bounds in
section \ref{sec:dynapp}.\\

For the atlas in Hypothesis \ref{hypAtlas}
we may take as charts e.g.\ 
a collection of $4^d$ cubes $\V_1,...,\V_{4^d}$ of side-length $\frac12$
placed in a uniform way on $\T^d$. Then any $\delta_0\in (0,1/4)$ 
is a Lebesgue number and we may identify lengths on charts by
the lengths in corresponding cubes $V_j$. 
The metric in each chart is the standard Euclidean metric.
From the conditions in the introduction on the matrix $M$ and the function
$\kappa$ we note that
for any $v\in \R^d$ we have $|\partial_x T_u(x)v|\geq \lambda |v|$
with $\lambda = \frac{1}{\theta_0}-\theta_1>1$, so $T_u$ is indeed uniformly
$\lambda$-expanding with $\lambda>1$.

Let $n=|\det(M)|\geq 2$.
The map $x\in \T^d\mapsto M x\in \T^d$
is an $n$-fold covering map of the torus. 
By a deformation and continuity 
argument the same holds for our map $T_u$.
For every $1\leq j\leq 4^d$ we thus have  $n$ local inverses
of $T_u$: $\psi_j^i: \V_j \rr \T^d$, $i\in \calA=\{1,...,n\}$.
By the above all local inverses of $T_u$ are
$\frac{1}{\lambda}$-Lipschitz
contracting. 
Applying the implicit function theorem to $T_u(\psi_j^i(u,y))=y$
and using the specific form for $T_u$
we get for each inverse branch:
\begin{equation*}
\partial_y \psi_j^i(u,y) = 
\left( M+ \partial_x \kappa(u,\psi_j^i(u,y) \right)^{-1}
\end{equation*}
\begin{equation*}
\partial_u \psi_j^i(u,y) = -   
\left( M+ \partial_x \kappa(u,\psi_j^i(u,y) \right)^{-1}
\partial_u \kappa(u,\psi_j^i(u,y) .
\end{equation*}

From these we obtain 
Fa\`a di Bruno like algebraic
formulae. Every $q$-th order derivative of $\psi$ will be
given in terms of derivatives up to order $q$ of $\kappa$ and 
up to order $q-1$ of $\psi$
itself. Since $\|(M+\partial_x \kappa)^{-1}\|\leq \frac{1}{\lambda}$
 it follows
by a recursive argument that if $\kappa$ is uniformly
bounded in $c^{r_0}$-norm then so is every $\psi$.
This shows that our inverse branches verify
condition (\ref{eq:kT}). 

For condition
(\ref{eq:condg}),
the weight for the natural measure is given by 
$g(x) =-\log \left| \det DT(x) \right|$ 
which when composed with an inverse branch
on $V_j$ takes the form:
\begin{eqnarray}
g_j^i(u,y) &=& - \log \left| \det \partial_y \psi_j^i(u,y)\right| 
             \nonumber \\
   &=& - \log n -
   \log \det\left(1 + M^{-1} \partial_x \kappa\right).
   \label{formulagij}
\end{eqnarray}
% Now, we have
% \[ \left|\det\left(M+\partial_x \kappa\right)\right| = 
% n \left(1 + \det\left(1 + M^{-1} \partial_x \kappa\right)\right).  \]
An eigenvalue of $M^{-1}\partial_x \kappa$ never exceeds 
$\theta_0\theta_1 <1$ in absolute value. Therefore,
  \begin{equation}
     0<(1-\theta_0\theta_1)^d \leq 
 \det\left(1 + M^{-1} \partial_x \kappa)\right) \leq
   (1+\theta_0\theta_1)^d .
     \label{bound:detM}
     \end{equation}
    As the constant $n$ disappears when
   taking  derivatives,
   we conclude that 
$|Dg_j^i|$ is uniformly bounded. 
Repeating the arguments from above we see that 
$\|D(g_j^i)\|_{(k-1,\alpha)}$ is uniformly bounded as well, proving (\ref{eq:condg}).

The upper and lower bounds sum of weights, condition
(\ref{eq:condS}), follows from formula (\ref{formulagij})
for  $g_j^i$ and the bound (\ref{bound:detM}),
 yielding the uniform bound:
 \[
 \frac{1}{(1+\theta_0\theta_1)^d}
 \leq 
     \sum_{x: T_u(x)=y} e^{g(u,x)} \leq 
 \frac{1}{(1-\theta_0\theta_1)^d}
  \].

We are thus in  position to
apply Theorems
\ref{quenchedlinrespthm} and
\ref{annealedlinrespthm}, which achieves our goal for this example.\\

For our Example \ref{ex:introtwo}, 
the manifold consists of one chart $V=(0,1)$.
We consider auxiliary functions of the type
$\kappa\in C^{r_0}(\R\times [0,1])$ with  $\|\kappa\|_{r_0} \leq K$,
$\partial_x \kappa(u,x)\geq \theta> \frac{1}{2}$ and $\kappa(u,0)=0$,
$\kappa(u,1)=1$.
The parametrized expanding maps $T_u : D_u \rr V=(0,1)$ given by
$T_u=1/\kappa - \lfloor 1/\kappa \rfloor$,
is defined whenever $1/\kappa$ is not an integer, so that $D_u$ is $V$ minus a countable set.
For fixed $u$ there are countably many inverse branches, which we may index by $i\in \N$ the corresponding
integer value of $\frac{1}{\kappa}$ for that branch.
\\We thus have:
$T_u(\psi^i(u,y))= \frac{1}{\kappa(u,\psi^i(u,y))}-i=y$,
$i\in \calA=\N$. Since 
$\partial_x T_u=-\partial_x \kappa/\kappa^2$ and $\partial_u T_u=-\partial_u \kappa/\kappa^2$ we get for each inverse
branch ($i\in \N$) the two identities:
\begin{equation*}
\partial_y \psi^i(u,y) = \frac{-1}{(i+y)^2\partial_x \kappa(u,\psi^i(u,y))}
\ \ \mbox{and} \ \
\partial_u \psi^i(u,y) = -\frac{\partial_u \kappa(u,\psi^i(u,y))}{\partial_x \kappa(u,\psi^i(u,y))} 
\end{equation*}
As in the previous example the first condition \eqref{eq:kT} is clearly verified.
\\For the weight condition note that $\partial_x \kappa$ takes values in $[\frac{1}{K},2]$ implying that 
$\log(\partial_x\kappa(u,\psi_j^i(y)))$ is $C^{r_0-1}$ bounded, and similarly the derivative of $\log (i+y)$ is $C^{r_0-2}$ bounded.
\\By summability of $\left(\frac{1}{i^2}\right)_{i\geq 1}$ \eqref{eq:condS} also holds. 
\\The only lacking element is the uniform contraction
(since $\partial_y\psi(u,y)$ need not be smaller than one).
For this we introduce the metric in the chart $(0,1)$:
$g(x)=\frac{1}{(1+x)^2}$, or equivalently given by
the line element $ds=\frac{1}{1+x}|dx|$. Our atlas verifies
Hypothesis (\ref{hypAtlas})
and for 
the metric derivative we have:
\begin{equation*}
\|\partial_y
\psi^i(u,y)\|_{g_y} = \frac{1+y}{1+\psi^i(u,y)} 
\left| \partial_y \psi^i(u,y)\right|.
\end{equation*}
From $\partial_x\kappa\geq \theta$ we get
$1-\kappa = \int_x^1 \partial_x \kappa \, dx \geq \theta(1-x)$
or $\kappa = \frac{1}{y+i} \leq \theta x + (1-\theta)$. Then
as $\theta>1/2$ and $x=\psi^i(u,y)\in [0,1]$:
\begin{equation*}
\|\partial_y\psi^i(u,y)\|_{g_y} 
\leq  \frac{1+y}{1+x}\times \frac{1}{\theta (i+y)^2}
\leq \frac{1}{(1+x)\theta} \frac{1}{i+y} \leq
\frac{\theta x + 1-\theta}{\theta x+\theta} \leq \frac{1}{2\theta} <1.
\end{equation*}
We may then apply our theorems as wanted.
Note that for $\theta\in(0,\frac12)$ we may find $\kappa$
satisfying all other conditions and
for which the associated $T_u$ has an attractive
fixed point.
\Halmos\\

\subsection{Application: Hausdorff dimension of repellers for 1D expanding maps}\label{subsec:Hausdorffdimandstuff}

In this section we are interested in the random product of one-dimensional maps, with uniform dilation but not necessarily defined everywhere. More precisely, we are interested in the following class of systems:

\begin{definit}\label{def:cookie-cutters}
	Let $I_1,\dots I_N\subset[0,1]$ be disjoint intervals, and $r\geq 2$. A $C^r$ map $T:I_1\cup\dots\cup I_N\rightarrow[0,1]$ is called a \textbf{cookie-cutter} if it satisfies the following conditions:
	\begin{itemize}
		\item There exists some $\lambda>1$ such that $\inf|T'|\geq\lambda$
		\item For each $i\in\{1,\dots,n\}$, $T(I_i)=[0,1]$
	\end{itemize}
	If $T$ is a cookie-cutter, we introduce its \textbf{repeller},
	\begin{equation*}
	\Lambda:=\{x\in I_1\cup\dots\cup I_N, T^n(x)~\text{is well-defined for all n}  \}=\bigcap_{i=1}^\infty T^{-i}([0,1])
	\end{equation*}
	We will denote by $CC^r([0,1])$ the set of all $C^r$ cookie-cutters.
\end{definit}
In other words, a cookie-cutter is a one-dimensional expanding map with full branches. It is a well-known fact that the repeller associated to such a map is a Cantor set. We now define perturbed cookie-cutters, in the following way:

\begin{definit}\label{def:ccperturb}
	
	Let $\U$ be an open subset of some Banach space $\B$, and let $\psi_i:\U\times[0,1]\rightarrow (0,1)$, $i\in\{1,\dots,N\}$ be $C^r$ maps such that
	
	\begin{itemize}
		\item For every $i\in\{1,\dots,N\}$, $\|\partial_x\psi_i\|_{\infty}\leq 1/\lambda<1$
		\item For every $i\in\{1,\dots,N\}$, every $u\in\U$, the intervals $I_{i,u}=\psi_{i,u}([0,1])$ are pairwise disjoint.
	\end{itemize}
	
	This data defines a cookie-cutter map $T_u$ on $I_{1,u},\dots,I_{N,u}$, by $T_u=\psi_{i,u}^{-1}$ on $I_{i,u}$. We call it a perturbed cookie-cutter.
\end{definit}

The question we want to study is the following: if one were to choose at each step a random cookie-cutter, and then perturb it in the sense of definition \ref{def:ccperturb}, does the Hausdorff dimension of the repeller change in a smooth way?
\medskip

To answer that question, we will use a random version of Bowen formula, which connects the transfer operator cocycle's top characteristic exponent and the Hausdorff dimension of the associated (random) repeller (cf. theorem \ref{Bowenformula}). 
\\More precisely, one considers a random product $T_{\omega,u}^{(n)}:=T_{\tau^{n-1}\omega,u}\circ\dots\circ T_{\omega,u}$ of perturbed cookie-cutters, i.e we assume that each $T_{\omega,u}\in\P_{k_T,k_g,k_\Sigma}$ is a perturbed cookie-cutter in the sense of Definition \ref{def:ccperturb}.
\\Associated to this random product is a random repeller, defined by
\begin{equation}\label{randrepeller}
\Lambda_{\omega,u}:= \bigcap_{i=1}^\infty\left(T_{\omega,u}^{(i)}\right)^{-1}\left([0,1]\right). 
\end{equation}
Given $t\geq 0$, we also define the transfer operator $\L_{\omega,t,u}$ by
\begin{equation*}\label{perturbedtransferop2}
\L_{\omega,t,u}\phi(x):=\sum_{T_{\omega,u}y=x}\frac{1}{|T'_{\omega,u}(y)|^t}\phi(y)
\end{equation*}

It follows from Theorem \ref{unifconecontraction} that $\hBL_{t,u}$ is a cone-contracting cocycle in the sense of Section \ref{sec:cone-con}, for $\C_{\vec{a}}$. 
\\It is also clear from the definition that $\L_{\omega,t,u}$ depends analytically on $t\geq 0$ (up to considering a small complex extension of $t$), so that it follows from \cite[Theorem 10.2]{Rugh10} that for every fixed $u_0\in\U$, the map $t\mapsto\chi_{t,u_0}$ is analytic.
\medskip

We also introduce the following quantities:
\begin{align*}
M_n(\omega,t,u)&:=\sup_{y\in\Lambda_{n,\omega,u}}\L_{\omega,t,u}^{(n)}{1}(y)\\
m_n(\omega,t,u)&:=\inf_{y\in\Lambda_{n,\omega,u}}\L_{\omega,t,u}^{(n)}{1}(y)
\end{align*}

where $\Lambda_{\omega,n,u}:=\bigcap_{i=n}^\infty\left(T_{\omega,u}^{(i)}\right)^{-1}([0,1])$, and finally we let \[-\infty\leq\underline{P}(\omega,t,u):=\liminf\frac{1}{n}\log(m_n(\omega,t,u))\leq\overline{P}(\omega,t,u):=\limsup\frac{1}{n}\log(M_n(\omega,t,u))\leq +\infty\]

Those last quantities exists by super-multiplicativity (resp. sub-multiplicativity) and Kingman's ergodic theorem, and are $\mathds{P}$-almost surely constant by ergodicity of $\tau$. 
\\One can show that those quantities almost surely agree, their common value being $\chi_{t,u}$ the top characteristic exponent of the random product, and that $t\geq 0\mapsto\chi_{t,u}+t\log(\lambda)\in\mathbb{R}$ is decreasing (see \cite[Lemma 3.5 and Theorem 4.4]{Rugh08}).
\\Furthermore, this decreasing map admits a unique zero that coincide with the (a.s) Hausdorff dimension of the random repeller $\Lambda_{\omega,u}$ (see \cite[Theorem 4.4 and 5.3]{Rugh08}):

\begin{theorem}[\cite{Rugh08} Theorem 5.3]\label{Bowenformula}
	Let $\tau$ be an invertible and ergodic map of $(\Omega,\mathds{P})$. Let $(T_{\omega})_{\omega\in\Omega}$ be a random product of cookie-cutters, such that $\mathds{E}\left[\log\|T'_\omega\|_\infty\right]<\infty$. 
	\\Then $\mathds{P}$-almost surely the Hausdorff dimension of the random repeller $\Lambda_\omega$ is given by the unique zero $z(T)$ of the top characteristic exponent $\chi_t$ of the transfer operator cocycle $\hBL_t$.
\end{theorem}

For a proof, we refer to \cite[§4-5]{Rugh08}. The question is now the dependence of that zero on the parameter $u$ :

\begin{theorem}\label{HausdorffdimisC1}
	Let $(T_{\omega,u})_{\omega\in\Omega,u\in\U}$ be a random product of perturbed cookie-cutters.
	\\Then the Hausdorff dimension of the random repeller defined by \eqref{randrepeller} is $C^s$ with respect to $u\in\U$.
\end{theorem}

\paragraph{Proof of theorem \ref{HausdorffdimisC1}}
Theorem \ref{Bowenformula} entails that the almost-sure Hausdorff dimension of $\Lambda_{\omega,u}$ is given by some $z(u)$ such that $\chi_{z(u),u}=0$.
\\From theorem \ref{Bowenformula} and \ref{thm:field}, one has that
\begin{itemize}
	\item For every $u\in\U$, $\chi_{z(u),u}=0$
	\item The map $(t,u)\mapsto\chi_{t,u}$ is $C^{r-2}$.
	\item $\partial_t\chi(t,u)\leq -\log(\lambda)<0$
\end{itemize}
Hence the result follows from the implicit function theorem.
\Halmos

\appendix

\section{Real cone contraction theory}
\label{appendix:conecontraction}

As the literature is somewhat disparate on the subject,
we provide here a catalog of more or less standard
results on cone contractions. 
Most of these results already occur in similar forms in e.g.
\cite{Bir57,Liv95,Dubois08}. In a few cases, we have added
(short) proofs where appropriate.
In the following, $(E,\|.\|)$ will denote a real Banach space.
%\subsection{Properties of Birkhoff cones}

\begin{definit}\label{def:realcone}
	Let $\C\subset E$. We say that $\C$ is a proper closed convex cone 
	(for short, a Birkhoff cone) if
	\begin{itemize}
		\item $\mathbb{R}_+\C=\C$, i.e $\C$ is stable by multiplication with a positive scalar.
		\item $\C$ is a closed and convex subset of $E$. 
		\item
		%We will say that the cone is \textbf{proper} if
		$\C\cap(-\C)=\{0\}$  (the cone is proper).
	\end{itemize}
	
	%We also define the \emph{dual} of the real closed convex cone $\C$ to be the set of non-zero functionals on $\C$, i.e
	
	%\begin{equation}\label{def:dualcone}
	%\C':=\{m\in E', \langle m,x\rangle\geq 0~\forall x\in\C^*\}
	%\end{equation}
	%where $\C^*=\C\backslash\{0\}$. We also call
\end{definit}

\begin{definit}\label{Bconeprop}
	
	Let $\C\subset E$ be a Birkhoff cone.  We say that $\C$ is
	
	\begin{enumerate} 
		\item \textbf{inner regular} iff $\C$ has non-empty interior.
		Equivalently, there is $\rho>0$ so that:
		\begin{equation*}
		\C(\rho) = \{ x\in \C: B_E(x,\rho \|x\|)\subset \C\}
		\label{app inner reg}
		\end{equation*}
		is non-trivial (contains other points than the origin).
		\item 
		\textbf{outer regular} if there
		is $\ell\in E'$, $\|\ell\|=1$ and $K<+\infty$ such that for every $u\in \C$: 
		\begin{equation}
		\frac{1}{K} \|u\| \leq
		\langle \ell,u\rangle 
		\leq \|u\|.
		\label{condouter}
		\end{equation}
	\end{enumerate}
	We will say that $\C$ is \textbf{regular} if it is both inner and outer regular.
	We define:
	\begin{equation*}
	\C_{\ell=1}(\rho)=\{u \in \C(\rho): \langle \ell,u\rangle=1\}.
	\label{def:Clrho}
	\end{equation*}
\end{definit}

\begin{definit}\label{Hmetric}
	Let $\C\subset E$ be a Birkhoff  cone, and let $x,y\in\C^*$. We define
	the projective \textbf{Hilbert metric}:
	\begin{align*}
	\delta(x,y)&=\inf\{t>0,~tx-y\in\C\}\\
	d_\C(x,y)&=\log(\delta(x,y)\delta(y,x))
	\in [0,+\infty].
	\end{align*}
	
	Hilbert's original (but equivalent) definition was through cross-ratios
	which may be formulated as follows:
	When $x,y\in \C^*$ are non-colinear,
	it is always possible to normalize $x$ and $y$ so that 
	$I(x,y) = \{t\in \R :  (1+t)x+(1-t)y \in \C\} = [t_1,t_2]$ 
	is a bounded (closed) interval.
	with $t_1\leq -1<1\leq t_2$. One then has for their projective distance,
	relative to the cone:
	\begin{equation}
	d_\C(x,y) = \log \frac{t_2+1}{t_2-1} \frac{t_1-1}{t_1+1}
	\in [0,+\infty]
	\label{Hilbert 2nd}
	\end{equation}
\end{definit}

\begin{lemme}\label{logHmetric}
	Let $\C\subset E$ be a regular Birkhoff cone, and let $d_\C$
	be the associated Hilbert metric. 
	Suppose that
	$B(x_1,r_1)\subset\C$ and $B(x_2,r_2)\subset \C$
	for some $r_1,r_2>0$. Then 
	\begin{equation*}
	d_{\C}(x_1,x_2)
	\leq 
	\log \left(1 + \frac{\|x_1-x_2\|}{r_1}\right)
	+ \log \left(1 + \frac{\|x_1-x_2\|}{r_2}\right) 
	\leq \left( \frac{1}{r_1} + \frac{1}{r_2} \right)
	\|x_1 - x_2\|.
	\end{equation*}
\end{lemme}
\Proof
Inclusion of the two balls imply that 
$I(x_1,x_2) \supset 
\left [-1 - \frac{2r_1}{\|x_1-x_2\|}, +1+ \frac{2r_2}{\|x_1-x_2\|} \right]$
which together with (\ref{Hilbert 2nd}) yields the result.
\Halmos

\begin{corol}
	For $x,y\in
	\C_{\ell=1}(\rho)$:
	$\ds d_\C(x,y) 
	\leq 2\log\left( 1 + \frac{\|x-y\|}{\rho}\right)
	\leq \frac{2}{\rho}\|x-y\|$.\\
	When $x\in\C_{\ell=1}(\rho)$ and $\|u\|<\rho$ we have:
	$\ds d_\C(x,x+u) \leq \frac{2 \|u\|}{\rho-\|u\|}
	= \frac{2}{\rho}\|u\| + o(\|u\|)$.
	\label{corol dC bound}
	\label{CORA5}
\end{corol}

%In the applications, it is of foremost importance to be able to compare the Hilbert metric with the distance given by the norm of the ambient Banach space. It is the object of this next lemma:
\begin{lemme}[\cite{Dubois08}, Appendix A]\label{Hmetricetnorme}
	Let $\C\subset E$ be an outer regular Birkhoff cone
	with $\ell\in E'$ as above.
	Then for all $x,y\in\C^*$:
	\begin{equation*}
	\left\|\frac{x}{\langle \ell,x\rangle}-
	\frac{y}{\langle \ell,y\rangle}\right\|\leq\frac{K}{2}d_\C(x,y)
	\end{equation*}
\end{lemme}

%We can now state the main result of this subsection: a linear operator sending a cone inside another shrinks the Hilbert metric.

% \begin{theorem}[Birkhoff's theorem]\label{Birkhoffthm}
% Let $E_1,~E_2$ be Banach spaces, and $\C_1\subset E_1$, $\C_2\subset E_2$ be proper, closed convex cones. Let $\L:E_1\rightarrow E_2$ be a linear map, such that $\L(\C_1^*)\subset\C_2^*$. Let $\Delta=diam_{\C_2}(\L(\C_1^*))\in[0,+\infty]$. Then one has 
% \begin{equation*}
% d_{\C_2}(\L x,\L y)\leq\tanh\left(\frac{\Delta}{4}\right)d_{\C_1}(x,y)
% \end{equation*}
% \end{theorem}

\begin{theorem}[Birkhoff's theorem, \cite{Bir57}]\label{myBirkhoffthm}
	Let $\C\subset E$ be a Birkhoff cone and let $\L\in L(E)$ be a
	contraction of $\C^*$, i.e. 
	such that $\L(\C^*)\subset\C^*$.
	Setting $\Delta=diam_{\C}(\L(\C^*))\in[0,+\infty]$ we have
	for $x,y\in \C^*$:
	\begin{equation*}
	d_{\C}(\L x,\L y)\leq
	\left( \tanh\frac{\Delta}{4}\right)d_{\C}(x,y)
	\end{equation*}
\end{theorem}

\begin{corol}
	\label{CORA8}
	Let $\C$ be a regular cone and $n\in\mathbb N$. Let $(L_i)_{1\le i\le n}\in L(E)$ be cone contractions, i.e $L_i(\C^*)\subset\C^*$ for any $1\le i\le n$ and set $\Delta_i=diam_{\C}(L_i(\C^*))\in[0,+\infty]$.
	Then for all $x,y\in\C^*$:
	\begin{equation}\label{eq:majcontraction}
	\left\|\frac{L_n\dots L_1 x}{\langle \ell,L_n\dots L_1 x\rangle}-
	\frac{L_n\dots L_1 y}{\langle \ell,L_n\dots L_1y\rangle}\right\|
	\leq\frac{K}{2} \left(\prod_{i=1}^n\tanh \frac{\Delta_i}{4}\right) d_\C(x,y)
	\end{equation}
\end{corol}

\begin{lemme}
	\label{lemma:upperlower}
	%\label{lemma:upper_lower}
	Let $\C$ be a regular cone with associated $\ell$, $K<+\infty$
	and $\rho>0$ as in definition \ref{Bconeprop}.
	Let $\L\in L(E)$ with $L(\C)\subset \C$.
	Then for every $x \in \C(\rho)^*$:
	\begin{equation*}
	\frac{\rho}{K} \|\L\| \; \|x\|   \leq
	\langle \ell, \L x \rangle \leq
	\|\L\| \; \|x\|  \ \ \mbox{and} \ \
	\frac{\rho}{K} \|\L\| \leq
	\frac{\langle \ell, \L x \rangle}
	{\langle \ell, x \rangle}
	\leq
	K \|\L\| .
	\end{equation*}
\end{lemme}
\Proof Suppose $x\in \C(\rho)$, $\|x\|=1$
and consider  $\|u\|\leq 1$. 
Then $x\pm \rho u \in \C$ so also 
$\L(x\pm \rho u) \in \C$. By outer regularity:
$\langle \ell,\L (x\pm \rho u) \rangle \geq 
\frac{1}{K}\|\L (x\pm \rho u)\|$.
Therefore, 
\begin{eqnarray*}
	2\rho \|\L u\| & \leq&  \|\L(x+\rho u)\|+\|\L(x-\rho u)\| \\
	&\leq& K \langle \ell , \L (x+\rho u) + \L (x-\rho u) \rangle \\
	&=& 2K \langle \ell , \L x \rangle 
\end{eqnarray*}
Thus, $\frac{\rho}{K} \|\L \| \leq \langle \ell,\L x\rangle$
from which we deduce the left-most inequality. The rest follows from \eqref{condouter}. \Halmos\\

%This result has important consequences for the study of spectral properties of cone-preserving operators: we will state two of them. 

%First, a quasi-compact operator preserving a Birkhoff cone has a spectral gap: this is the content of the \emph{Krein-Rutman theorem}:
%\begin{theorem}[\cite{KR50}, Theorem 6.3]
%Let $\C$ be an inner regular Birkhoff cone of a real Banach space $E$. Let $\L:E\to E$ be a quasi-compact operator, such that $\L(\C^*)\subset Int(\C)$. Then $\L$ admits a spectral gap. 
%\end{theorem}

%Second, it is noteworthy that the conditions on the spectrum of the cone-contracting operator can be replaced by geometrical assumptions on the cone: a linear operator preserving a reproducing Birkhoff cone with bounded subsectional aperture (e.g, if the cone is regular) has a \emph{spectral gap}, without condition on its spectrum: this a theorem of Birkhoff \cite{Bir57}.
%
%\begin{theorem}
%Let $\C\subset E$ be a reproducing Birkhoff cone with bounded sectional aperture. We assume that $\L: E\to E$ is a strict cone contraction, i.e $\L(\C^*)\subset \C^*$ and $diam_{\C}(\L(\C^*))<+\infty$. Then $\L$ has a spectral gap.
%\end{theorem}

%We end this subsection with an intuitive and useful lemma; as it is not clearly stated anywhere in the literature to the best of our knowledge, we provide a proof here.

\begin{lemme}\label{lem:uniopen}
	Let $\C$ be a regular Birkhoff cone with 
	associated linear functional $\ell\in E'$ and $K<+\infty$
	as in Definition \ref{Bconeprop}.
	Let $\C_1\subset \C$ be a subcone of finite diameter,
	$\Delta=diam_{\C}(\C_1^*)<+\infty$.
	Suppose that there is $x\in \C_1$, $\langle \ell,x\rangle=1$ with
	$B(x,r)\subset \C_1$.
	%Let $\C_1\subset \C_2$ be two regular Birkhoff cones, such that $\Delta:=diam_{\C_2}\C_1^*<+\infty$.
	Then for every $y\in\C_1$, 
	\begin{equation*}
	B_E\left(y,
	\frac{1}{K} r e^{-\Delta}
	\|y\|\right)\subset\C.
	\end{equation*}
\end{lemme}
\Proof
Pick  $x,y\in \C_{\ell=1}$ and set
$u_t=\frac12 ( (1+t)y+(1-t)x)$. Then
$\{t\in \R :  u_t \in \C\} = [t_1,t_2]$ 
is a bounded (closed) interval (definition \ref{Hmetric})
with $t_1\leq -1<1\leq t_2$. 
Then $\Delta \geq d_{\C}(x,y) \geq \log \frac{t_2+1}{t_2-1}$
or $\frac{t_2-1}{t_2+1}\geq e^{-\Delta}$.
For $|u|<1$ we have $x+r u\in \C$. Also $u_{t_2}\in \C$ so
the following convex combination is also in $\C$:
\[ \frac{(t_2-1) (x+ru)+2u_{t_2}}{t_2+1} 
= y+r\frac{t_2-1}{t_2+1} u \in \C.\]
Here $\frac{1}{K}\|y\| \leq \langle \ell,y\rangle=1$ so 
$B(y,\rho\|y\|)\subset \C$
with $\rho = \frac{1}{K}re^{-\Delta}$.\Halmos\\

\section{Bochner and strong measurability}
\label{appendix:measurability}
In our setup we need measurability of quantities related
to sections and operators. It is close to standard Bochner
measurability and strong measurability in the sense of e.g.
\cite[Appendix A]{GTQ14}
but not quite the same so we bring here a brief account of
the notions we use.

In this appendix $(X,|\cdot|_X)$ denotes a Banach space and
$(\Omega,\calF)$ a non-empty space equipped with a $\sigma$-algebra.
\begin{definit}
	\label{def:Bmeas}
	\begin{enumerate}
		\item
		A map $\phi:\Omega\rr X$ is said to be $\sigma$-simple if it has
		a countable image and is measurable (with respect to $\calF$).
		We write $S_\sigma(\Omega,X)$ for the set of such functions.
		\item
		A map $\psi:\Omega\rr X$ is said to be Bochner-measurable if  it
		may be written as the uniform limit of a sequence of $\sigma$-simple
		functions.
		We write $L^\infty(\Omega;X)$ for the set
		of Bochner measurable maps that are uniformly bounded.
		It is a Banach space under the uniform norm.
		\item
		Let $A\subset X$ be a non-empty set.
		We say that $\BL : \Omega \rr L(X)$ is strongly 
		Bochner measurable on $A$ provided 
		that $\omega \in \Omega \mapsto \BL(\omega) a \in X$ is Bochner-measurable
		for every $a\in A$.
		\item
		We say that $\BL$ in the previous definition is measurably bounded
		provided there is a measurable function $\rho:\Omega \mapsto [0,+\infty)$
		so that $\|\BL(\omega)\|_{L(X)} \leq \rho(\omega)$
		for every $\omega\in \Omega$.
		
		$\BL$ is of course  bounded if  
		$\sup_{\omega\in \Omega} \|\BL(\omega)\|_{L(X)} <+\infty$.
		
	\end{enumerate}
\end{definit}

\begin{rem}
	A pointwise limit of measurable functions is measurable
	\cite[VI,\textsection 1,M7]{Lang93}, so the above definition a Bochner-measurable
	function is equivalent to saying that $\psi$ is measurable with
	image having a countable dense subset.
\end{rem}

\begin{prop}
	\label{Prop-Bochner} 
	Let $\phi:\Omega \rr A\subset X$ be Bochner measurable
	with values in $A$ and let $\BL:\Omega \rr L(X)$ be measurably
	bounded and strongly Bochner measurable on $A$. Then
	\begin{equation*}
	\omega\in \Omega \mapsto \BL(\omega) \phi(\omega)\in X
	\end{equation*}
	is also Bochner measurable.
\end{prop}
\Proof Let $\rho$ be as in the last part of Definition \ref{def:Bmeas}.
For $m=0,1,...$ we set $\Omega_m = \rho^{-1} ([m,m+1))$ which provides
a measurable partition of $\Omega$. Let $\epsilon>0$.
By Bochner measurability, and for each $m\geq 0$
we may find
sequences $x_{m,k}\in A\subset X$ and
$E_{m,k}\in \calF$, $k\geq 1$ so that $(E_{m,k})_{k\geq 1}$ form
a measurable partition and
$|\phi(\omega)-x_{m,k}|\leq \frac{\eps}{2(m+1)}$. By strong measurability
on $A$ we have for each $m\geq 0, k\geq 1$ sequences
$y_{m,k,\ell}\in X$ and $F_{m,k,\ell}\in \calF$, $\ell\geq 1$
so that $(F_{m,k,\ell})_{\ell\geq 1}$ forms a measurable
partition and $|\BL(\omega)x_{m,k}-y_{m,k,\ell}|_X \leq \eps/2$
for every $\omega\in F_{m,k,\ell}$. Then for
$\omega\in G_{m,k,\ell}=\Omega_m\cap E_{m,k}\cap F_{m,k,\ell}$
we have
\[ 
|\BL(\omega)\phi(\omega) -y_{m,k,\ell}|_X \leq 
|\BL(\omega)(\phi(\omega)-x_{m,k})|
+ |\BL(\omega)x_{m,k} -y_{m,k,\ell}|_X \leq
\eps.\]
The $G_{m,k,\ell}$-collection gives a measurable 
partition of $\Omega$ and the conclusion follows.
\Halmos

\section{Differential Calculus}
\label{app:diffcalc}
We provide a listing  of
some of the more or less standard results in differential calculus
which we are using.
Let $U$ and $V$ denote open convex sets in Banach spaces
$B_U$ and $B_V$, respectively,
and let $Z$ be a fixed Banach space. 
We write $|\phi|_0$ 
(or $\|\phi\|_0$) to denote
a uniform norm on the relevant domain of definition.\\

Let $r>0$ and $\delta_1 \in (0,1]$ be fixed number in this
section.  We write $\undr=(k,\alpha)$, with
$k\in \N_0$ and $\alpha\in (0,1]$ such that
$r=k+\alpha$. $C^\undr$-norms between open convex subsets
of Banach spaces are defined in the standard way.
By $c^\undr(U;Z)$ we understand the closure
of $C^{s}(U;Z)$ in $C^r(U;Z)$ with $s>r$ (the closure is independent
of the choice of $s$).

Given $A\in L(E_1\times \cdots \times E_n; Z)$, a multilinear form
with $E_1,...,E_n$ being Banach spaces, 
there is a natural isomorphism obtained by singling out the $j$'th
Banach space, yielding:
$A_j\in L(E_1\times \cdots \widehat{E}_j
\cdots \times E_n; L(E_j;Z))$.
In many places in this article we make use
of the telescopic principle which asserts that 
\begin{equation}\label{eq:telescopic}
\|A(M_1,\ldots,M_n) - A(N_1,\ldots,N_n) \|\leq
\sum_{j=1}^n \|A_j(M_1,...,M_{j-1},N_{j+1},\ldots,N_n)\| \; 
\|M_j - N_j\|.\end{equation}

%This principle is used very frequently when obtaining the estimates below
%for H\"older constants.

\begin{prop}
	\label{prop:banalg}
	Let $\phi_1,\phi_2\in c^\undr(V;\CX)$, $\phi\in c^\undr(V;Z)$ and
	$\psi\in c^\undt (U;V)$, $t=r\vee 1$. Then there are constants
	$C_{r,1},C_{r,2},C_{r,3},C_{4,r}$ depending only upon $r$ and
	the chosen norms so that:\\
	
	\begin{tabular}{l r c l c l}
		1. & 
		$\phi_1\phi_2$ & $\in$ & $ c^\undr(V;\CX)$  & and &
		$\|\phi_1\phi_2\|_r \leq C_{1,r} \|\phi_1\|_r \|\phi_2\|_r$.\\[2mm]
		2. &
		$e^{\phi_1}$ & $\in$ & $ c^\undr(V;\CX)$ & and &
		$\|e^{\phi_1}\|_r \leq 
		C_{2,r} |e^{\phi_1}|_0 \; 
		(1 + \|\phi_1\|_{r})^r$\\[2mm]
		3. &
		$\phi\circ \psi $ & $\in $ & $c^\undr (U; Z)$ & and &
		$\|\phi \circ \psi\|_r \leq 
		C_{3,r}  \|\phi_1\|_r (1 + \|D\psi\|_{t-1})^r$,\\[3mm]
	\end{tabular}    
	
	and also:\\
	
	\begin{tabular}{lll}
		4.  &  &
		$\|e^{\phi_1}-e^{\phi_2}\|_r \leq 
		C_{4,r} \; (|e^{\phi_1}|_0 \vee |e^{\phi_2}|_0)\;
		\;\|\phi_1-\phi_2\|_r
		\;(1 + \|\phi_1\|_{r}\vee \|\phi_2\|_r)^r$\\[2mm]
	\end{tabular}
\end{prop}
\Proof

1. For $C^k$-functions with $k$ an integer, this is standard.
When $\undr=(k,\alpha)$ with $\alpha \in(0,1]$  we have
$ D^k(\phi_1\phi_2) = 
(D^k \phi_1) \phi_2 +  
\phi_1 (D^k \phi_2)+R_{k-1}$, where $R_{k-1}$
is a bilinear form in $\phi_1$ and $\phi_2$
involving derivatives of order at most $k-1$.
The $\delta_1$-local H\"older constant may then be estimated
using
the MVT 
on the last term
and the above-mentioned telescopic principle
to obtain:
$h_{\delta_1}^\alpha(D^k (\phi_1\phi_2)) \leq
h_{\delta_1}^\alpha(D^k \phi_1)|\phi_2|_0  +
|\phi_1|_0 h_{\delta_1}^\alpha(D^k \phi_2) + |D R_{k-1}|_0$,
The last term is bounded by the $C^k$ norm of $\phi_1\phi_2$.

2. Write $D^q (e^{\phi_1}) = e^{\phi_1} (D\phi_1 + D)^q \bfone$,
$0\leq q \leq k$ and
develop. For the H\"older constant for
$D^k(e^{\phi_1})$ we use the same argument as above.

3. The so-called Fa\`a di Bruno formula gives a combinatorial
expression for $D^k (\phi\circ \psi)$. Exhibiting the $k$'th order
derivative one has:
\begin{equation*}
D^k (\phi\circ \psi) =
D^k \phi \circ \psi. \underbrace{(D\psi,\ldots,D\psi)}_{k\ {\rm times}} +
D \phi \circ \psi. D^k \psi + R_{k-1}.\label{eq:Faa}
\end{equation*}

Again $R_{k-1}$ involves only derivatives up to order $k-1$
and may be treated using the MVT.
For the H\"older constant of the $k$'th derivative,
let $ K = \|D\psi\|_0 $. If $y,y'\in U$ are at a distance at most $\delta_1$,
then the $\psi$-image of this path has length (in $V$) at most 
$\whdel_1=K\delta_1$. Lemma
\ref{holdereq}  then shows that
$h_{\delta_1}^\alpha (D^k\phi\circ \psi) \leq \lceil K \rceil^{1-\alpha}
h_{\delta_1}^\alpha (D^k\phi)$.

4. For the last inequality, the MVT yields
$|e^{\phi_1}-e^{\phi_2}|_0 \leq 
L (|e^{\phi_1}|_0 \vee |e^{\phi_2}|_0)\; \;|\phi_1-\phi_2|_0$.
Developing $D^q(e^{\phi_1}-e^{\phi_2})$ and using the telescopic 
principle yields the wanted bound.
The above calculations were done
for $C^r$ functions. But the uniform bounds
implies that the result easily
carries over to $c^\undr$ functions as well.\Halmos\\

\subsection{Regularity when extracting a parameter}\label{appendix:regularity}
In this section we will show how the regularity
of a function of two variables behaves when extracting one
variable as a parameter. The notation is as above.
We have the following smoothness result when extracting
a variable as a parameter:

\begin{theorem} 
	\label{lemma:injection}
	We equip the product $B_U\times B_V$
	with the max norm.
	Let $r,s,t \geq   0$ with $t=r-s>0$. 
	We have the following canonical continuous injections:
	\begin{eqnarray}
	\phi \in C^{r}(U\times V;Z) &\lhrr&
	\hphi \in C^{{\undt}}(U; C^{s}(V;Z))\\
	\phi \in c^{r}(U\times V;Z) &\lhrr&
	\hphi \in C^{{\undt}}(U; c^{s}(V;Z))
	\label{eqn:smoothinj}
	\end{eqnarray}
	under the natural identification: $\hphi_u(x) := \phi(u,x)$, $u\in U, x\in V$.
	\label{thm:parextract}
\end{theorem}
\paragraph{Proof:}
Consider the first inclusion. Denote $Y_s= C^{s}(V;Z))$
and let $\phi\in C^{r}(U\times V;Z)$. 
Our claim is that $\hphi_u\in Y_s$
and that the map $u\in U \mapsto \hphi_u\in Y_s$ is $C^{\undt}$
with $t=r-s>0$. \\

Case 1: We first show this for $0\leq  s < r\leq 1$.
Here $\|\phi\|_{r}= |\phi|_0 \vee h^r_{\delta_1}(\phi)$ is simply
the local $r$-H\"older norm.
For fixed $u\in \U$,
clearly $|\hphi_u|_0\leq |\phi|_0$ and 
$h^s_{\delta_1}(\hphi_u) \leq h^r_{\delta_1}(\hphi_u)\leq h^r_{\delta_1}(\phi)$ (since $0<\delta_1\leq 1$).
so $\hphi_u\in Y_s$
with $\|\hphi_u\|_{s}\leq \|\phi\|_r$.
To check  the regularity w.r.t.\ $u$ pick $u_1,u_0\in U$
with $0<|u_1-u_0|\leq \delta_1$ and set 
$\Delta=\hphi_{u_1}-\hphi_{u_0}\in Y_s$.
We have
$| \Delta |_0
\leq  h^r_{\delta_1}(\phi) 
|u_1-u_0|^{\,r}$.
\\To check H\"older regularity with respect to $x$ consider
$x_0,x_1\in V$ with $0<|x_1-x_0|\leq \delta_1$. Then 
$\Delta(x_1)-\Delta(x_0)=
\phi(u_1,x_1)-
\phi(u_1,x_0)-
\phi(u_0,x_1)+
\phi(u_0,x_0)
$
may be estimated using H\"older-regularity either with respect to $u$
or to $x$.
This yields (the middle term is maximal
when 
$|x_1-x_0|= |u_1-u_0|$):
\begin{equation*}
\frac{|\Delta(x_1)-\Delta(x_0)|}
{ |x_1-x_0|^{\,s}}
\leq 2h^r_{\delta_1}(\phi) 
\frac{ |x_1-x_0|^{\,r} \wedge |u_1-u_0|^{\,r}}
{ |x_1-x_0|^{\,s} }
\leq 2h^r_{\delta_1}(\phi) 
|u_1-u_0|^{r- s}.
\end{equation*}
Thus,
$h_{s} (\Delta) \leq
2 h^r_{\delta_1} (\phi)
|u_1-u_0|^{r- s}$
from which:
$ \|\hphi_{u_1} - \hphi_{u_0} \|_{Y_s} \leq
2\|\phi\|_r |u_1-u_0|^{r- s}$. 
So
$\hphi\in C^{(0,r-s)}(U;Y_s)$ with $\|\hphi\|_{r-s}\leq 2 \|\phi\|_r$.\\

Case 2: Consider now when $0 \leq s\leq 1< r \leq 1+s$.
Again for fixed $u$: $\hphi_u\in Y_s$ with $\|\hphi_u\|_s\leq \|\phi\|_r$.
With $\Delta$ as above we have
$|\Delta|_0\leq |D\phi|_0 \; |u_1-u_0| \leq \|\phi\|_r\;  |u_1-u_0|^s$.
\\Let $u(\tau)=\tau u_1+(1-\tau)u_0$, $0\leq \tau\leq 1$ be the segment 
joining $u(0)=u_0$ and $u(1)=u_1$ (included in $U$ by convexity). Applying the MVT, one gets:
\begin{align*}
|\Delta (x_1) -\Delta(x_0)| & \leq
\int_0^1 
\left|\frac{d}{d\tau} ( \hphi_{u_\tau}(x_1)-\hphi_{u_\tau}(x_0)) \right|d\tau \\
&\leq \int_0^1 |D\hphi_{u_\tau}(x_1)-D\hphi_{u_\tau}(x_0)) | \; |u'(\tau)| d\tau \\
&\leq h^{r-1}_{\delta_1}(D \phi) |x_1-x_0|^{\; r-1} |u_1-u_0|.
\end{align*}
Interchanging the roles of $x_1,x_0$ and $u_1,u_0$ we also have:
\begin{equation*}
|\Delta(x_1)-\Delta(x_0)| 
\leq h^{r-1}_{\delta_1}(D \phi) |u_1-u_0|^{r-1}|x_1-x_0|,
\end{equation*}
and then (again the middle term is maximal for $|u_1-u_0| = |x_1-x_0|$):
\begin{align*}
\frac{ |\Delta(x_1) \!-\! \Delta(x_0)|} {|x_1-x_0|^{s}}
&\leq  
h^{r-1}_{\delta_1}(D\phi) 
\frac{
	|x_1 \!-\! x_0|^{\; r-1} |u_1 \!-\! u_0| \wedge |u_1 \!-\! u_0|^{\; r-1} |x_1 \!-\! x_0| }
{|x_1-x_0|^{\; s}}\\
&\leq h^{r-1}_{\delta_1}(D\phi) |u_1-u_0|^{\; r-s}
\leq \|\phi\|_r \; |u_1-u_0|^{\; r-s}.
\end{align*}
It follows that $\|\hphi_{u_1}-\hphi_{u_0}\|_{Y_s}
\leq \|\phi\|_r |u_1-u_0|^{\;r- s}  $ and
$\hphi\in C^{(0,r-s)}(U;{Y_s})$ with $\|\hphi\|_{r-s}\leq  \|\phi\|_r$. 
Note that when $r-s=1$ the conclusion is that
\begin{equation*}
\phi \in C^{(1,s)}(U\times V;Z) \lhrr
\hphi \in C^{(0,1)}(U; C^{s}(V;Z)),
\end{equation*}
i.e. in general, $u\mapsto \hphi_u$ need not be differentiable in this case, but it is Lipschitz continuous.
\\Higher order regularity can be reduced to the above two cases. 
\\To see this let $0\leq \beta < \alpha \leq 1$ and $r=k+m+\alpha$,
$s=m+\beta$ with $k,m\in \N_0$. For $\phi\in C^{r}(U\times V;Z)$ and fixed $u\in\U$ it is clear that $\hphi_u\in C^{(m,\beta)}(V;Z)$.
For the regularity of the map $u\mapsto \hphi_u$ 
we leave intermediate derivatives to the reader and
consider only the highest order.
Note that the first case treated above yields the following injection and identification (using natural isomorphisms for the linear maps involved):
\begin{eqnarray*}
	\partial_u^k \partial_x^m \phi 
	& \in &
	C^{(0,\alpha)}(U\times V;L(B_U^k \times B_V^m; Z)) \\
	&\lhrr &  
	C^{ (0,\alpha-\beta)}(U; C^{(0,\beta)}(V;L(B_U^k \times B_V^m; Z))) \\
	&\simeq &  
	C^{ (0,\alpha-\beta)}(U; L(B_U^k ;
	C^{(0,\beta)}(V;L(B_V^m; Z)))). \\
	&= &  
	C^{ (0,\alpha-\beta)}(U; L(B_U^k ; W))
\end{eqnarray*}
with $W=C^{(0,\beta)}( V;L(B_V^m; Z))$.
We observe that the latter
precisely gives the
identification with
$\partial_x^m ( ( \partial_u^k \hphi_u))$
and implies that 
$ \partial_u^k \hphi_u \in 
C^{ (0,\alpha-\beta)}(U;L(B_U^k; C^{s} (V;Z)))$
whence that 
$ \hphi_u \in C^{ \undt}  
(U;C^{s}(V;Z)) $ 
with $\undt=(k,\alpha-\beta)$
as we wished to show.
\\In the case $0<\alpha\leq \beta \leq 1$ and $r=k+m+1+\alpha$,
$s=m+\beta$ we consider again $\partial_x^k \partial_u^m \phi$
which reduces the necessary injection
to the second case treated above. In either case the norm
increases at most by a factor of $2$.
\medskip

For the second injection,  if $\phi$ is the $C^r$
limit of smooth functions then 
the induced 
function $\hphi_u$ is the $C^s$ limit of smooth functions. 
The norm-estimates carry over from before.
%\\Note that in general, one may not replace
%$C^{\undt}$ by $c^{\undt}$ in the injection \eqref{eqn:smoothinj}.
\Halmos\\

\begin{rem}
	For non-compact $V$ and 
	for integer values $n\geq k\geq 0$,
	there is in general no natural injection of
	$C^n (U \times V;\R)$ into $C^{n-k}(U ; C^{k}(V;\R))$.
	An easy (counter) example is for $n=k=0$, 
	$U=V=(0,1)$, where you may
	consider e.g.\ $\phi(u,x)=\sin(u/x)$.
\end{rem}

\begin{corol}
	\label{cor:compodiff}
	Assume that $W\subset \B_W$ is an open convex subset of a Banach space.
	Let $r>1$ and $\beta\in(0,1]$.
	Let $\psi_0,\psi_1\in C^r(V;W)$ and $\phi\in C^{r+\beta}(W;Z)$.
	We write 
	$M=\|\psi_0\|_r \vee \|\psi_1\|_r \vee 1$.
	Suppose that 
	$\|\psi_0 -\psi_1 \|_r \leq \delta_1$. Then we have
	\begin{equation*}
	\|\phi\circ \psi_1 - \phi\circ \psi_0\|_{r} \leq
	C_{r,\beta} \|\phi\|_{r+\beta} \; M^r  \,
	\|\psi_0 -\psi_1 \|_r^\beta
	\end{equation*}
	with $C_{r,\beta}$ a constant that
	depends only upon $r$, $\beta$ and the choice of norms.
\end{corol}
\Proof We first show this under the additional assumption that
$\psi_0,\psi_1\in C^{r+\beta}$. 
Let $a=\|\psi_0 -\psi_1 \|_r$. 
We may assume that $a>0$ or else there is nothing to show.
We define for $(t,y) \in [0, a] \times V$ the following
linear interpolation between $\psi_0$ and $\psi_1$ (allowed since $W$ is convex):
\begin{equation*}
\tpsi(t,y) = 
\frac{t}{a} \psi_1(y) + \left(1-\frac{t}{a}\right) \psi_0(y)
= \psi_0(y) + 
\frac{t}{a} \left(\psi_1(y) - \psi_0(y)\right).
\end{equation*}
One has $\tpsi\in C^{r+\beta}([0,a]\times V;W)$\footnote{
	$[0,a]\times V$ is not open, but the construction of $C^r$ functions
	works equally well on a space obtained by intersecting an open convex set
	with a closed half-space.}.
Note that $\partial_t \tpsi= \frac{1}{a}(\psi_1-\psi_0)$
has $(r+\beta-1)$-norm not greater than 1. It follows 
that $\|\tpsi\|_{r+\beta}\leq
\|\psi_0\|_{r+\beta}\vee\|\psi_1\|_{r+\beta}\vee 1$. 
It follows from Proposition \ref{prop:banalg} that
$F=\phi\circ \tpsi$ has $(r+\beta)$-norm bounded by
$K=C'_{r,\beta}\|\phi\|_{r+\beta} M_{r+\beta}^r$.
\\By our
parameter-extraction theorem \ref{thm:parextract}
we deduce that
$t\in [0,a] \mapsto (F(t,\cdot) \in C^{r}(U;Z))$  is $C^{(0,\beta)}$
with at most twice the indicated bound for the norm.
But then the H\"older bound implies:
\begin{equation*}
\|\phi\circ\psi_1-\phi\circ \psi_0\|_r =
\|F(a,\cdot) - F(0,\cdot)\|_r \leq 2K |a-0|^\beta
= 2K 
\|\psi_0 -\psi_1 \|_r^\beta
\end{equation*}
as we wanted to show. Returning to the general case,
consider the telescopic form for the derivative:
\[D(\phi\circ \psi_1)- D(\phi\circ \psi_0)
= (D\phi\circ \psi_1 -D\phi\circ\psi_0).D\psi_1
+ (D\phi\circ \psi_0).(D\psi_1-D\psi_0)\]
Here,
$D\phi\in C^{r+\beta-1}$ and $\psi_1\in C^{r}\subset C^{r+\beta-1}$
(since $\beta\leq 1$). The first part then applies
(with $r+\beta-1$ instead of $r+\beta$) and shows that
$\|D\phi\circ \psi_1 -D\phi\circ\psi_0\|_{r-1}\leq
C_{r-1,\beta} \|D\phi\|_{r+\beta-1} \|\psi_1-\psi_0\|^\beta$.
The last term trivially verifies the same type of bound. 
From this we deduce the result for the difference without derivatives.
\Halmos\\

\begin{prop}\label{prop:extractregforop}
	Let us consider $r=k+\alpha$ with $k\in\mathbb N$ and $0<\alpha\leq 1$, and $\undr=(k,\alpha)$, $\B,X,Y$ three Banach spaces and $\U\subset\B$ an open subset. Under the natural identifications, we have the following injections:
	\begin{align*}
	L(X,C^\undr(\U,Y))&\hookrightarrow C^\undr(\U,L(X,Y))
	%L(X,c^\undr(\U,Y))&\hookrightarrow c^\undr(\U,L(X,Y))
	\end{align*}
\end{prop}
\Proof We present a proof by induction on $k\in \mathbb N$. 
\\For $\undr=(0,\alpha)$, we assume that we have an operator $\L_u:X\to Y$, satisfying: there is a $C>0$, such that for any $\phi\in X$, any $u\not=v,~u,v\in\U$,
\begin{align*}
\L_u\phi&\in C^{0,\alpha}(\U,Y)\\
\|\L_u\phi\|_Y&\leq C\|\phi\|_X\\
\|\L_u\phi-\L_v\phi\|_Y&\leq C\|\phi\|_X\|u-v\|^\alpha_\B
\end{align*}
Define $\hat\L$ as the map $u\in\U\mapsto \L_u\in L(X,Y)$. Then it is easy to see that under the previous assumptions, $\hat\L$ is a $C^\alpha$ map 
%(the $c^\alpha$ case is similary).
\\Let us assume that the wanted injection is established at rank $k-1$, and consider $r=k+\alpha$, an operator $\L_u:X\to Y$, with $\L_u\phi\in C^\undr(\U,Y)$ and $\|\L_u\phi\|_{C^r(U,Y)}\leq C\|\phi\|_X$.
\\For any $\phi\in X$, we may consider the partial derivative (w.r.t $u$) of $\L_u\phi$, $\partial_u(\L_u\phi)\in L(\B,Y)$, which is, by assumption, a $C^{(k-1,\alpha)}$ map w.r.t $u\in\U$, with $\|\partial_u\L_u\phi\|_{C^{k-1,\alpha}}\leq C\|\phi\|_X$. But by induction hypothesis, this means that the map $\hat\L:\U\to L(X,Y)$ admits a derivative which is $C^{k-1,\alpha}$, i.e that $\hat\L$ is $C^{k,\alpha}$.\Halmos
\subsection{Bochner measurable smooth sections}
\label{app:Crsections}
Let $(\Omega,\F)$  be a measurable space. In this section there
is no measure involved.
We write
\begin{equation*}
X_{r}(U;Z) :=
X_{(k,\alpha)}(U;Z) 
= L^\infty(\Omega; c^{(k,\alpha)}(U;Z))
\end{equation*}
for 
a Bochner measurable map from $\Omega$ to the Banach space
$Y=c^{(k,\alpha)}(U;Z)$, with $r=k+\alpha>0$, $\alpha\in(0,1]$, $k\in\N_0$.
The norm of $\Bphi\in X_r(U;Z)$ is the uniform norm:
$\|\Bphi\|_{X_r} = \sup_{\omega\in\Omega} \|\Bphi_\omega\|_r$.
Operations in the following proposition is understood to take 
place fiber-wise, e.g.\ for fixed $\omega\in\Omega$,
$(\Bphi\circ\Bpsi)_\omega :=
\Bphi_\omega \circ\Bpsi_\omega$.
We have:

\begin{prop}
	\label{prop:Xbanalg}
	Let $\Bphi_1,\Bphi_2\in X_r(V;\CX)$, $\Bphi\in X_r(V;Z)$ and
	$\Bpsi\in X_s (U;V)$, $s=r\vee 1$. Then there are constants
	$C_{r,1},C_{r,2},C_{r,3}$ depending only upon $r$ and
	the chosen norms so that:\\
	%\begin{enumerate}
	%\item $\phi_1\phi_2\in C^\undr(V;\CX)$ and
	%$\|\phi_1\phi_2\|_r \leq C_{1,r} \|\phi_1\|_r \|\phi_2\|_r$.
	%\item
	%$e^{\phi_1}\phi_2\in C^\undr(V;\CX)$ and $\|e^{\phi_1}\phi_2\|_r \leq 
	%C_{2,r} |e^{\phi_1}|_0 \; \|\phi_2\|_r (1 + \|D\phi_1\|_{r-1})^r$
	%\item $\phi\circ \psi \in C^\undr (U; Z)$ and 
	%$\|\phi \circ \psi\|_r \leq 
	%C_{3,r}  \|\phi_1\|_r (1 + \|D\psi\|_{r-1})^r$.
	%\end{enumerate}
	
	\begin{tabular}{l r c l c l c l}
		1. & 
		$\Bphi_1 \Bphi_2$ & $\in$ & $ X_r(V;\CX)$  & and &
		$\|\Bphi_1\Bphi_2\|_{X_r}$ & $\leq $&$
		C_{1,r} \|\Bphi_1\|_{X_r} \|\Bphi_2\|_{X_r}$.\\[2mm]
		2. &
		$e^{\Bphi_1}$ & $\in$ & $ X_r(V;\CX)$ & and &
		$\|e^{\Bphi_1}\|_{X_r} $&$ \leq $&$
		C_{2,r} |e^{\Bphi_1}|_0 \;  
		(1 + \|D\Bphi_1\|_{X_{r-1}})^r$\\[2mm]
		3. &
		$\Bphi\circ \Bpsi $ & $\in $ & $X_r (U; Z)$ & and &
		$\|\Bphi \circ \Bpsi\|_{X_r}$ &$ \leq $&$
		C_{3,r}  \|\Bphi\|_{X_r} (1 + \|D\Bpsi\|_{X_{s-1}})^r$,
		\ \ $(r>1)$.
	\end{tabular}    
\end{prop}
\Proof As operations are fiber-wise we clearly have the stated bounds
on the norms. The only issue is Bochner-measurability. We show this
for the first case: Let
$M  >  \|\Bphi_1\|_{X_r} \vee  \|\Bphi_2\|_{X_r}$.
Given $\eps>0$ we may find a countable measurable partition
$(\Omega_m)_{m\in \N}$  so that
for every $m\in\N$, $\omega,\omega'\in \Omega_m$, $i=1,2$ we have:
$\|(\Bphi_i)_\omega -  (\Bphi_i)_{\omega'}\|_r
\leq \frac{\eps}{2 C_{1,r} M}$. Then by the telescopic principle
and the above bounds:
\begin{equation*}
\|(\Bphi_1 \Bphi_2)_\omega -  (\Bphi_1 \Bphi_2)_{\omega'}\|_r
\leq 2 C_{1,r} \frac{\eps}{2 C_{1,r} M} M \leq \eps,
\end{equation*}
implying that $\Bphi_1\Bphi_2$ is Bochner measurable in
the sense of definition
\ref{def:Bmeas}. The other two statements follow in the
same way. \Halmos\\

\begin{lemme}
	\label{lem:inftyinject}
	One has the following injections:
	\begin{equation*}
	\Bphi\in L^\infty(\Omega; C^\undr(V; Z)) \lhrr 
	\Bhphi \in C^{\undr}(V; L^\infty(\Omega;Z)),
	\end{equation*}
	\begin{equation*}
	\Bphi\in L^\infty(\Omega; c^\undr(V; Z)) \lhrr 
	\Bhphi \in c^{\undr}(V; L^\infty(\Omega;Z)),
	\end{equation*}
	under the natural fiber-wise identification:
	\begin{equation*}
	\Bhphi(v)(\omega) := \Bphi(\omega)(v),   \ \ \
	\omega\in \Omega, v\in V.
	\end{equation*}
\end{lemme}
\Proof
Given $\eps>0$ we find a countable measurable partition
$(\Omega_i)_{i\in \N}$  so that
for every $i\in\N$, $\omega,\omega'\in \Omega_i$, we have:
$\|(\Bphi)_\omega -  (\Bphi)_{\omega'}\|_r
\leq {\eps}$. Pick also for each $i\geq1$: $\omega_i\in\Omega_i$
and set $f_i=\Bphi(\omega_i)$. Define:
\begin{equation*}
\Bhf(v)(\omega) =\Bf(\omega)(v)
:=\sum_i \bfone_{\Omega_i}(\omega) f_i(v).
\end{equation*}
Then $\Bf$ is a $\sigma$-simple $\eps$-uniform approximation to
$\Bphi$. Clearly, $\Bhf$ takes values in $Y=L^\infty(\Omega,Z)$.
Also $\partial_v^q \Bhf(v) \in L(B_V^q; Y)$ and 
$\|\partial_v^q \Bhf -
\partial_v^q \Bhphi\|\leq \eps$, $0\leq q\leq k$
and similarly for the $\alpha$-H\"older estimate for the $k$'th
derivative. Thus $\Bhphi \in c^{\undr}(V; L^\infty(\Omega;Z))$
and it has the same norm as $\Bphi$. \Halmos\\

%\begin{rem} Note that the reverse inclusion need not hold above:
%With $\Omega=V=(0,1)$, 
%$Z=\R$, $\alpha\in(0,1]$,
%$\Bhphi : (v,\omega)\in V\times \Omega
%\mapsto |v-\omega|^\alpha$ defines a function in 
%$L^\infty(\Omega; c^\undr(V; Z))$ but
%$\Bphi \not\in c^{\undr}(V; L^\infty(\Omega;Z))$,
%since $\|\Bphi(\omega)-\Bphi(\omega')\|\geq 2$
%whenever $\omega\neq \omega'$ so 
%$\omega\mapsto \Bphi(\omega)$ is not measurable.
%\end{rem}

We conclude this section with a key ingredient for our applications section:
\begin{prop}
	\label{prop:Bparamextr}
	With the notation as in Theorem \ref{lemma:injection} and this section, 
	we have the
	following injection of norm at most 2:
	\begin{equation*}
	\Bphi\in X_r(U\times V; Z) \lhrr 
	\Bhphi \in C^{{r-s}}(U; X_s(V;Z)),
	\end{equation*}
	under the natural fiber-wise identification:
	\begin{equation*}
	\Bhphi(u)(\omega)(v) := \Bphi(\omega)(u,v),   \ \ \
	\omega\in \Omega, u\in U, v\in V.
	\end{equation*}
\end{prop}
\Proof
Combining 
\ref{lemma:injection}
and
\ref{lem:inftyinject} 
we have the following injections:
\begin{align*}
X_r(U\times V; Z)  \ \
& = \ \ \  L^\infty (\Omega; c^\undr(U\times V; Z)) \\
&  \lhrr \ \  L^\infty (\Omega; C^{r-s}(U; c^\unds(V; Z)) \\
&  \lhrr  \ \  C^{r-s}(U; L^\infty(\Omega; c^\unds(V;Z)))\\
&  = \  \ \  C^{r-s}(U; X_s(V;Z)). 
\end{align*}

\section{Graded differential calculus}
\label{app:graded-diff-calc}
An essential ingredient in differential calculus is
the Leibnitz principle: When e.g.\ 
$f,g$ are $C^r(\R)$-functions for $r\geq 1$ ($r$ not necessarily an integer)
then so is their product and 
one has a formula for the derivative of the product
$(f\cdot g)'=f' \cdot  g + f \cdot  g'$. The derivative
is then $C^{r-1}$ and one may iterate the derivation
formula when $r\geq 2$.
The aim here is to develop a similar theory for graded
differential calculus, in particular the Leibniz principle,
when $f,g$ are replaced
by linear operators depending on a parameter $u$ but where 
regularity with respect to the
parameter only appears when
downgrading the codomain (the image space) or upgrading the
domain within a certain scale of Banach spaces. 
The upshot of this appendix is to show that the resulting regularity when
performing algebraic operations on graded differential 
operators is as good as one could possibly hope for. In particular, we prove Lemma \ref{lem:induc}.
\\

We will stick to the notation of section \ref{sec:loss-of-reg}.
More precisely, 
let ${\cal X}=(X_t)_{t\in (0,r_0]}$ denote the 
scale of Banach spaces. By this we mean a parametrized family
of Banach spaces 
coming with a family of bounded linear (downgrading) operators
$j_{s,r}\in L(X_r;X_s)$, $0< s\leq r\leq r_0$. We assume that each
operator is injective and has dense image and that the
collection satisfies the  transitivity condition: $j_{s,s}={\rm Id}$
and 
$j_{s,c}j_{c,r}=j_{s,r}$ whenever $0< s\leq c\leq r\leq r_0$.
%A family $(\phi_r)_{r\in[a,b]}$ with each
%$\phi_r\in X_r$ is $j$-equivariant if $\phi_s=j_{s,r} \phi_r$ for
%all $a\leq s<r \leq b$.

\begin{exemple}
	An instructive  example 
	to have in mind is
	$X_t = C^t(S^1)$, $t \in (0,r_0]$ with
	$j_{s,r}: C^r(S^1) \rr C^s(S^1)$ 
	being the natural embedding for $0< s\leq r\leq r_0$.
	\label{ex:downgrade}
\end{exemple}

We let $\B$ denote a  Banach space
and let $U\subset \B$ be a non-empty open convex subset.

\begin{definit}
	\label{def:j-equiv}
	Let $\myD\in \N_0=\{0,1,\ldots\}$,
	$\gamma>0$ with $\gamma+\myD\leq r_0$. 
	We associate to the integer $\myD$ the following
	set of ordered pairs:
	$I_\myD=
	\{(s,r)\in (0,r_0]^2 : s + \myD \leq r\}$.
	Consider a family  $\M$ of
	bounded linear operators 
	$M_{s,r}(u) \in L(\B^\myD; L(X_r,X_s))$, $(s,r)\in I_\myD$ 
	and parametrized by $u\in U$.
	We say that the family $\M$
	is ($j$-)equivariant and $(\gamma,\myD)$-regular provided
	that:
	\begin{enumerate}
		\item 
		For every $(s,r),(s',r')\in I_\myD$ with $s<s'$, $r<r'$
		and $u\in U$:
		\ $\ds  j_{s,s'} M_{s',r'}(u)
		= M_{s,r}(u) j_{r,r'}$.
		\item The map \ 
		$u\in U \mapsto M_{s,r}(u)\in L(\B^\myD; L(X_r;X_s))$ \ 
		is $C^t$ for all $(s,r)\in I_\myD$ and 
		$0 \leq t<\gamma \wedge (r-s-\myD)$.
	\end{enumerate}
\end{definit}

%In order to keep the notational complexity down,
%we assume in the following that the parameter space
%$\B=\R$ and that $U$ is an
%open interval of $\R$. This allows us below to take derivatives
%without introducing extra linear operators with respect to $\B$.
%All the results below also hold in the general case,
%cf. Remark \ref{rem:simplify}.

Keeping the same notation as in the previous definition we define:
\begin{definit}
	\label{def:j-leftequi}
	Consider a family
	$\calN$ of functions 
	$N_{s}(u) \in L(\B^\myD; X_s)$, $0 \leq s< r_0-\myD$, parametrized
	by $u\in U$. We say that $\calN$ 
	is left-equivariant and $(\gamma,\myD)$-regular  provided
	that
	for all $0\leq s<r< r_0-\myD$, $u\in U$:
	$N_{s}(u) = j_{s,r} N_{r}(u)$,
	and the map 
	$u\in U \mapsto N_{s}(u)\in L(\B^\myD;X_s)$ \ 
	is $C^t$ for all 
	$0 \leq t<\gamma\wedge (r_0-\myD-s)$.
\end{definit}

\begin{lemme}
	\label{lem:derivedfam}
	Let $\M$ be an equivariant $(\gamma,\myD)$-regular family
	with $\gamma>1$. We define the derived family 
	$\partial_u\M$ given by:
	$\partial_u M_{s,r}(u)\in
	L(\B; L(\B^\myD;L(X_s,X_r)))\equiv 
	L(\B^{\myD+1};L(X_s,X_r))$ 
	for
	all $(s,r)\in I_{\myD+1}$. This derived
	family is equivariant and
	$(\gamma-1,\myD+1)$-regular.
	Conversely, 
	suppose that $\M$ is $(1+\alpha,\myD)$-regular with 
	$\alpha>0$ and that derived family $\partial_u \M$ is 
	$(\gamma',\myD+1)$-regular, then
	setting $\gamma=\alpha\vee \gamma'+1$, we have that
	$\M$ is $(\gamma,\myD)$-regular. A similar statement holds
	for a left-equivariant family $\calN$.
\end{lemme}
\Proof The first statement is obvious from definitions.
For the second, we may assume $\gamma'>\alpha$ or else it is trivial.
Suppose that $\M$ is $(1+\alpha,\myD)$-regular with $\alpha>0$
and let $\partial_u \M$ be the derived family.
If $0<s<r\leq r_0$ with $r-s\leq r_0-\myD$ and
$1 < t< t_*=(\gamma'+1)\wedge(r-s-\myD)$.
Then $u\mapsto \partial_u M_{s,r}(u)$ is $C^{t-1}$
and consequently $u\mapsto M_{s,r}(u)$ is $C^t$ as we wanted to show.\Halmos
\\

%In Example \ref{ex:downgrade} 
%with $b-a>1$, the derivative operator
%$D = d/dx : C^b(S^1) \rr C^a(S^1)$ gives rise to 
%a $j$-equivariant
%and $1$-regular family over $[a,b]$.
The main reason for introducing
equivariant, $(\gamma,\myD)$-regular
families comes from  the stability
under products:

\begin{prop}
	\label{prop:prod}
	Let $\M^1$  and $\M^2$ 
	%$(M^1_r(u))_{r\in[a,b]}$ and $(M^2_r(u))_{r\in[a,b]}$, $u\in U$ 
	be two families of 
	$j$-equivariant, 
	$(\gamma_1,\myD_1)$-regular, respectively $(\gamma_2,\myD_2)$-regular,
	operators. Suppose that $\myD=\myD_1+\myD_2<r_0$ and
	set $\gamma=\gamma_1\wedge \gamma_2 \wedge (r_0-n)>0$.
	Then there is a well-defined product family
	$\M=\M^1\star \M^2$ obtained
	by declaring for $(s,r)\in I_\myD$:
	\begin{equation*}
	M_{s,r}(u):= M^1_{s,c}(u) M^2_{c,r}(u)\in 
	L(\B^{\myD_1}\times \B^{\myD_2}; L(X_r;X_s))
	\simeq L(\B^{\myD}; L(X_r;X_s))
	\end{equation*}
	with $c$ being any number in the non-empty interval
	$(s+\myD_1,r-\myD_2)$.
	This product family is 
	equivariant and $(\gamma,\myD)$-regular.
	Similarly, with $\M^1$ as above and $\calN^2$ a left-equivariant and
	$(\gamma_2,\myD_2)$-regular family,
	the product
	$N_s(u):=M^1_{s,c} N^2_c(u)$ defines a family
	$\calN=\M^1 \star \calN^2$ which is left-equivariant
	and $(\gamma,\myD)$-regular.
\end{prop}

\paragraph{Proof:} If $c<c'$ are two numbers in the above interval,
then we have by equivariance (all operators being well-defined):
\begin{equation*}
M^1_{s,c}(u) M^2_{c,r}(u)  = 
M^1_{s,c}(u) \left[ j_{c,c'} M^2_{c',r}(u) \right]  = 
\left[ M^1_{s,c}(u) j_{c,c'} \right] M^2_{c',r}(u)  = 
M^1_{s,c'}(u) M^2_{c',r}(u) ,
\end{equation*}
showing that the product does not depend upon the choice
of $c$. 
Let
$(s,r)\in I_\myD$ and $t_*=\gamma\wedge (r-s-\myD)>0$.
Regularity will be shown by induction in $\gamma$.  
First, assume that $\gamma \in (0,1]$.
Then, in particular, $t_*\leq 1$.
We will show that
$u\mapsto M_{s,r}(u)\in L(X_r;X_s)$ is $t$-H\"older for every
$0<t<t_*$. Let $0<\eps< t_*$ and set
$s':=s+\myD_1+\eps/3 < r':=r-\myD_2-\eps/3$.
We then have by equivariance:
\begin{equation*}
M_{s,r}(u) =
M^1_{s,s'}(u) j_{s',r'} M^2_{r',r}(u).
\end{equation*}
Note that $r'-s-\myD_1=r-s'-\myD_2=r-s-\myD-\eps/3$.
Also let $c=(s'+r')/2$.
When $u,u+h\in U$, using a telescopic sum, equivariance
and H\"older continuity we have the following identity: 
\begin{align}
M_{s,r}(u+h)&-M_{s,r}(u)
= 
M^1_{s,s'}(u+h) j_{s',r'} M^2_{r',r}(u+h)-
M^1_{s,s'}(u) j_{s',r'} M^2_{r',r}(u) \label{eq:M}\\
&= 
\left(M^1_{s,r'}(u+h)-M^1_{s,r'}(u)\right) 
M^2_{r',r}(u)
\label{eq:M12} \\
& +
\left(M^1_{s,c}(u+h)-M^1_{s,c}(u)\right) 
\left(M^2_{c,r}(u+h)-M^2_{c,r}\right)
\nonumber \\
& +
M^1_{s,s'}(u)  \
\left( M^2_{s',r}(u+h)-M^2_{s',r}(u) \right)
\nonumber \\
&= 
\O(h^{t_*-\eps/3}) \O(1) +
\O(h^{t_*/2-\eps/3})  \O(h^{t_*/2-\eps/3})  +
\O(1) \O(h^{t_*-\eps/3}) = 
\O(h^{t_* - 2\eps/3}).
\nonumber 
\end{align}
We may here let $\eps\rr 0$ and obtain the claim for this case.\\

For higher order regularity,
let $k\geq 1$ be an integer and suppose that the proposition has
been proven whenever
$0<\gamma \leq k$.  For $k=1$ this was done above. 
So assume now that $\gamma\in (k,k+1]$.
For $(s,r)\in I_\myD$ we set $t_*=\gamma\wedge(r-s-\myD)$.
We may assume that $t_*\in (k,k+1]$ as well
(or else there is nothing to show).
We write $t_*=k+\alpha$ with $\alpha\in (0,1]$.
We want to show that $u\mapsto M_{s,r}(u)$ is $C^t$ for all
$0<t<t_*$.
Set $s'=s+\myD_1+\alpha/3$ and $r'=r-\myD_2-\alpha/3$. 
Then as before $r'-s' =r-s-\myD-2\alpha/3 > k+\alpha/3$.
Since $k\geq 1$ we obtain derived families
$\partial_u \M^1$ and $\partial_u \M^2$ as described in
Lemma \ref{lem:derivedfam}.
We have e.g. $M^1_{s,s'}(u)j_{s',r'}= M^1_{s,r'}(u)$ so by the MVT we get:
\begin{eqnarray*}
	\lefteqn{
		\left| M^{1}_{s,r'}(u+h)-M^{1}_{s,r'}(u)
		- h \cdot \partial_u M^{1}_{s,r'}(u) \right|_{L(X_s;X_{r'})}  } 
	\hspace*{10mm} \\
	& \leq& 
	|h| \sup_{\tau\in[0,1]} 
	\left| 
	\partial_u M^{1}_{s,r'}(u+\tau h)
	-\partial_u M^{1}_{s,r'}(u) \right|_{L(X_s;X_{r'})}  \\
	& = & O(h^1) O(h^{\alpha/3}) =
	O(h^{1+\alpha/3}) \ .
\end{eqnarray*}
With a similar expansion for $M^2$ 
and using 
H\"older estimates for the middle term we expand
\eqref{eq:M12} to get:
\begin{align*}
M_{s,r}(u+h)-M_{s,r}(u)
&= 
(h.\partial_uM^1_{s,r'})(u) M^2_{r',r}(u)+
M^1_{s,s'}(u) 
\; (h.\partial_u M^2_{s',r})(u) + \O(h^{1+\alpha/3}),
\end{align*}
showing that $M_{s,r}$ is differentiable
with derivative
\begin{align*} \partial_u M_{s,r}(u) = 
\partial_u M^1_{s,r'}(u)  M^2_{r',r}(u)+
M^1_{s,s'}(u) \partial_u M^2_{s',r}(u).
\end{align*}
Now, in this expression we may again use equivariance to write
\begin{align*} \partial_u M_{s,r}(u) = 
\partial_u M^1_{s,c_1}(u)  M^2_{c_1,r}(u)+
M^1_{s,c_2}(u) \partial_u M^2_{c_2,r}(u).
\end{align*}
with $c_1\in (s+1+\myD_1, r-\myD_2)$ and
$c_2\in (s+\myD_1, r-1-\myD_2)$. The first term is the product of
two j-equivariant families that are
$(\gamma-1,n_1+1)$ and $(\gamma,n_2)$ regular, respectively.
%$(\myD_1+1)$-regular and $\myD_2$-regular, respectively.
Since $(\gamma-1)\wedge(r-s-n-1)=t_*-1\leq k$
%$\gamma'=r-s-(\myD_1+\myD_2) = \gamma-1\leq k$
we may apply the induction hypothesis on this term to
conclude that this first product is $(\gamma-1,\myD+1)$-regular.
Similarly for the second term. 
Thus $u\mapsto \partial_u M_{s,r}(u)$
is $C^{t}$ for every $t<\gamma-1$, whence 
$u\mapsto  M_{s,r}(u)$
is $C^{t}$ for every $t<\gamma$ as we wanted to show 
(see Lemma \ref{lem:derivedfam}).
The proof in the left-equivariant case follows the same path.
\Halmos\\

\begin{lemme}
	Let $(Q_{s,r}(u))$
	be an equivariant, $(\gamma,0)$-regular family with the
	additional property that $\bfone-Q_s(u)$ is invertible 
	for all $0<s\leq r_0$ and having a uniformly bounded
	inverse 
	$R_s(u) = (\bfone-Q_{s}(u))^{-1}$
	when $\epsilon<s\leq r_0$ for any $\eps>0$. 
	Then the family of operators $R_{s,r}(u)=\Bj_{s,r} R_r(u)$,
	$0<s<r\leq r_0$
	is again equivariant and
	$(\gamma,0)$-regular.
\end{lemme}
\Proof This boils down to the resolvent identity combined
with equivariance. We have e.g. for $u,u+h\in U$:
\begin{equation*}
R_{s,r}(u+h)-R_{s,r}(u) =
R_{s}(u+h) \left( Q_{s,r}(u+h)-Q_{s,r}(u)\right)R_{r}(u).
\end{equation*}
H\"older-continuity then follows using regularity of the middle term.
When $t_*=\gamma\wedge(r-s-1)>1$ we may again develop
the middle term and conclude that $R_{s,r}(u)$ is differentiable
with derivative:
\begin{equation*}
\partial_u R_{s,r}(u) =
R_{s}(u) \left( \partial_u Q_{s,r}(u)\right)R_{r}(u)
\in L(\B;L(X_r;X_s)).
\end{equation*}
Here we have a product of 3 operators being 
$(t_*,0)$, $(t_*-1,1)$ and $(t_*,0)$-regular, respectively.
The product is then itself
$(t_*-1,1)$-regular and therefore $R_{s,r}(u)$ is 
$(t_*,0)$-regular as we wanted to show. \Halmos\\

\paragraph{Proof of Lemma \ref{lem:induc}:}
First note that in Theorem \ref{thm:field}
the collection of operators 
$L_{s,r}(u):=\BL_{s,u}\Bj_{s,r}$
with $0<s\leq r\leq r_0$,
forms
an equivariant family $\L$ of $(r_0,0)$-regular operators over $(0,r_0]$.
The derived family 
$\left(\partial_u(\BL_{s,u}\Bj_{s,r_0})\right)_{s\in(0,r_0]}$
is then $(r_0-1,1)$-regular.

Under Hypothesis $\calH(\gamma)$ the family of fixed fields
$(\field_s(u))_{s\in(0,r_0]}$ is
left-equivariant
and $(\gamma,0)$-regular. From Proposition \ref{prop:prod}
it follows that the family of products 
$\left(\partial_u(\BL_{s,u}\Bj_{s,r_0}) \field_{r_0}(u)\right)_{s\in(0,r_0]}$
is $(\gamma\wedge(r_0-1),1)$-regular. This is in fact the principal term
in the definition \ref{eqforP} of $P_{s,r_0,u}(\field_{r_0}(u))$ which
is therefore also $(\gamma\wedge(r_0-1),1)$-regular:
this shows the first claim in Lemma \ref{lem:induc}.
In a similar way, using Proposition \ref{prop:piCt}
and $\calH(\gamma)$ we see that 
$M_{s,r}(u) = Q_{s,u}(\field_{s}(u))\Bj_{s,r}$ is 
equivariant and $(\gamma,0)$-regular.
This implies the second claim in the Lemma. \Halmos\\

\section*{Acknowledgment} JS was supported by the European Research Council (ERC) under the European Union's Horizon 2020 research and innovation program (grant agreement No 787304).

\begin{small}
\bibliography{biblio}

\begin{thebibliography}{10}

\bibitem{BRS17}
Wael Bahsoun, Marks Ruziboev, and Beno\^{\i}t Saussol.
\newblock Linear response for random dynamical systems.
\newblock {\em Adv. Math.}, 364:107011, 44, 2020.

\bibitem{BahsounSaussol16}
Wael Bahsoun and Beno\^it Saussol.
\newblock Linear response in the intermittent family: differentiation in a
  weighted {$C^0$}-norm.
\newblock {\em Discrete Contin. Dyn. Syst.}, 36(12):6657--6668, 2016.

\bibitem{Ba00}
Viviane Baladi.
\newblock {\em Positive transfer operators and decay of correlations},
  volume~16 of {\em Advanced Series in Nonlinear Dynamics}.
\newblock World Scientific Publishing Co., Inc., River Edge, NJ, 2000.

\bibitem{Ba14}
Viviane Baladi.
\newblock Linear response, or else.
\newblock In {\em Proceedings of the {I}nternational {C}ongress of
  {M}athematicians---{S}eoul 2014. {V}ol. {III}}, pages 525--545. Kyung Moon
  Sa, Seoul, 2014.

\bibitem{Ba16}
Viviane Baladi.
\newblock {\em Dynamical zeta functions and dynamical determinants for
  hyperbolic maps}, volume~68 of {\em Ergebnisse der Mathematik und ihrer
  Grenzgebiete. 3. Folge. A Series of Modern Surveys in Mathematics}.
\newblock Springer, Cham, 2018.
\newblock A functional approach.

\bibitem{BKS96}
Viviane Baladi, Abdelaziz Kondah, and Bernard Schmitt.
\newblock Random correlations for small perturbations of expanding maps.
\newblock {\em Random Comput. Dynam.}, 4(2-3):179--204, 1996.

\bibitem{BaladiSmania08}
Viviane Baladi and Daniel Smania.
\newblock Linear response formula for piecewise expanding unimodal maps.
\newblock {\em Nonlinearity}, 21(4):677--711, 2008.

\bibitem{BaladiSmaniaalt10}
Viviane Baladi and Daniel Smania.
\newblock Alternative proofs of linear response for piecewise expanding
  unimodal maps.
\newblock {\em Ergodic Theory Dynam. Systems}, 30(1):1--20, 2010.

\bibitem{BaladiSmania12}
Viviane Baladi and Daniel Smania.
\newblock Linear response for smooth deformations of generic nonuniformly
  hyperbolic unimodal maps.
\newblock {\em Ann. Sci. \'Ec. Norm. Sup\'er. (4)}, 45(6):861--926 (2013),
  2012.

\bibitem{BaladiTodd16}
Viviane Baladi and Mike Todd.
\newblock Linear response for intermittent maps.
\newblock {\em Comm. Math. Phys.}, 347(3):857--874, 2016.

\bibitem{Bir57}
Garrett Birkhoff.
\newblock Extensions of {J}entzsch's theorem.
\newblock {\em Trans. Amer. Math. Soc.}, 85:219--227, 1957.

\bibitem{BCV16}
T.~Bomfim, A.~Castro, and P.~Varandas.
\newblock Differentiability of thermodynamical quantities in non-uniformly
  expanding dynamics.
\newblock {\em Adv. Math.}, 292:478--528, 2016.

\bibitem{CSG11}
Micka\"el~D. Chekroun, Eric Simonnet, and Michael Ghil.
\newblock Stochastic climate dynamics: random attractors and time-dependent
  invariant measures.
\newblock {\em Phys. D}, 240(21):1685--1700, 2011.

\bibitem{Crim19}
Harry {Crimmins}.
\newblock {Stability of hyperbolic Oseledets splittings for quasi-compact
  operator cocycles}.
\newblock {\em arXiv e-prints}, December 2019.
\newblock \url{https://arxiv.org/abs/1912.03008}.

\bibitem{RDLL99}
R.~de~la Llave and R.~Obaya.
\newblock Regularity of the composition operator in spaces of {H}\"older
  functions.
\newblock {\em Discrete Contin. Dynam. Systems}, 5(1):157--184, 1999.

\bibitem{D17}
Thai~Son Doan.
\newblock On analyticity for {L}yapunov exponents of generic bounded linear
  random dynamical systems.
\newblock {\em Discrete Contin. Dyn. Syst. Ser. B}, 22(8):3113--3126, 2017.

\bibitem{Dol04}
Dmitry Dolgopyat.
\newblock On differentiability of {SRB} states for partially hyperbolic
  systems.
\newblock {\em Invent. Math.}, 155(2):389--449, 2004.

\bibitem{Drag2018}
D.~Dragi{\v{c}}evi{\'{c}}, G.~Froyland, C.~Gonz{\'a}lez-Tokman, and S.~Vaienti.
\newblock A spectral approach for quenched limit theorems for random expanding
  dynamical systems.
\newblock {\em Communications in Mathematical Physics}, Jan 2018.

\bibitem{DFGV19}
D.~Dragi\v{c}evi\'{c}, G.~Froyland, C.~Gonz\'{a}lez-Tokman, and S.~Vaienti.
\newblock A spectral approach for quenched limit theorems for random hyperbolic
  dynamical systems.
\newblock {\em Trans. Amer. Math. Soc.}, 373(1):629--664, 2020.

\bibitem{DF18}
Davor Dragi\v{c}evi\'c and Gary Froyland.
\newblock H\"older continuity of {O}seledets splittings for semi-invertible
  operator cocycles.
\newblock {\em Ergodic Theory Dynam. Systems}, 38(3):961--981, 2018.

\bibitem{Dubois08}
Lo\"{i}c Dubois.
\newblock Real cone contractions and analyticity properties of the
  characteristic exponents.
\newblock {\em Nonlinearity}, 21(11):2519--2536, 2008.

\bibitem{El13}
Japp Eldering.
\newblock {\em Normally Hyperbolic Invariant Manifolds}.
\newblock Atlantis Series in Dynamical Systems. Atlansits Press, first edition,
  2013.

\bibitem{FGTQ14}
Gary Froyland, Cecilia Gonz\'alez-Tokman, and Anthony Quas.
\newblock Stability and approximation of random invariant densities for
  {L}asota-{Y}orke map cocycles.
\newblock {\em Nonlinearity}, 27(4):647--660, 2014.

\bibitem{FGTQ15}
Gary Froyland, Cecilia Gonz\'alez-Tokman, and Anthony Quas.
\newblock Stochastic stability of {L}yapunov exponents and {O}seledets
  splittings for semi-invertible matrix cocycles.
\newblock {\em Comm. Pure Appl. Math.}, 68(11):2052--2081, 2015.

\bibitem{FLQ10}
Gary Froyland, Simon Lloyd, and Anthony Quas.
\newblock Coherent structures and isolated spectrum for {P}erron-{F}robenius
  cocycles.
\newblock {\em Ergodic Theory Dynam. Systems}, 30(3):729--756, 2010.

\bibitem{FLQ13}
Gary Froyland, Simon Lloyd, and Anthony Quas.
\newblock A semi-invertible {O}seledets theorem with applications to transfer
  operator cocycles.
\newblock {\em Discrete Contin. Dyn. Syst.}, 33(9):3835--3860, 2013.

\bibitem{GalatoloGiulietti17}
S.~Galatolo and P.~Giulietti.
\newblock A linear response for dynamical systems with additive noise.
\newblock {\em Nonlinearity}, 32(6):2269--2301, 2019.

\bibitem{GalatoloPollicott}
Stefano Galatolo and Mark Pollicott.
\newblock Controlling the statistical properties of expanding maps.
\newblock {\em Nonlinearity}, 30(7):2737--2751, 2017.

\bibitem{GalatoloSedro20}
Stefano Galatolo and Julien Sedro.
\newblock Quadratic response of random and deterministic dynamical systems.
\newblock {\em Chaos}, 30(2):023113, 15, 2020.

\bibitem{CSG08}
Michael Ghil, Micka\"el~D. Chekroun, and Eric Simonnet.
\newblock Climate dynamics and fluid mechanics: natural variability and related
  uncertainties.
\newblock {\em Phys. D}, 237(14-17):2111--2126, 2008.

\bibitem{GTQ14}
Cecilia Gonz\'alez-Tokman and Anthony Quas.
\newblock A semi-invertible operator {O}seledets theorem.
\newblock {\em Ergodic Theory Dynam. Systems}, 34(4):1230--1272, 2014.

\bibitem{GTQ18}
Cecilia {Gonz{\'a}lez-Tokman} and Anthony {Quas}.
\newblock {Stability and Collapse of the {L}yapunov spectrum for
  {P}erron-{F}robenius Operator cocycles}.
\newblock {\em arXiv e-prints}, June 2018.
\newblock \url{https://arxiv.org/abs/1806.08873}.

\bibitem{GL06}
S\'ebastien Gou\"ezel and Carlangelo Liverani.
\newblock Banach spaces adapted to {A}nosov systems.
\newblock {\em Ergodic Theory Dynam. Systems}, 26(1):189--217, 2006.

\bibitem{GL08}
S\'ebastien Gou\"ezel and Carlangelo Liverani.
\newblock Compact locally maximal hyperbolic sets for smooth maps: fine
  statistical properties.
\newblock {\em J. Differential Geom.}, 79(3):433--477, 2008.

\bibitem{GR81}
Mikhael Gromov.
\newblock Groups of polynomial growth and expanding maps.
\newblock {\em IHES Publ. Math.}, 1981(53):53--73, 1981.

\bibitem{HairerMajda10}
Martin Hairer and Andrew~J. Majda.
\newblock A simple framework to justify linear response theory.
\newblock {\em Nonlinearity}, 23(4):909--922, 2010.

\bibitem{Hennion91}
Hubert Hennion.
\newblock D\'erivabilit\'e du plus grand exposant caract\'eristique des
  produits de matrices al\'eatoires ind\'ependantes \`a coefficients positifs.
\newblock {\em Ann. Inst. H. Poincar\'e Probab. Statist.}, 27(1):27--59, 1991.

\bibitem{KL98}
Gerhard Keller and Carlangelo Liverani.
\newblock Stability of the spectrum for transfer operators.
\newblock {\em Ann. Scuola Norm. Sup. Pisa Cl. Sci. (4)}, 28(1):141--152, 1999.

\bibitem{Lang93}
Serge Lang.
\newblock {\em Real and Functional Analysis, 3rd edition}.
\newblock Springer, 1993.

\bibitem{LePage89}
\'Emile Le~Page.
\newblock R\'egularit\'e du plus grand exposant caract\'eristique des produits
  de matrices al\'eatoires ind\'ependantes et applications.
\newblock {\em Ann. Inst. H. Poincar\'e Probab. Statist.}, 25(2):109--142,
  1989.

\bibitem{Liv95}
Carlangelo Liverani.
\newblock Decay of correlations.
\newblock {\em Ann. of Math. (2)}, 142(2):239--301, 1995.

\bibitem{LucariniL12}
V.~Lucarini.
\newblock Stochastic perturbations to dynamical systems: a response theory
  approach.
\newblock {\em J. Stat. Phys.}, 146(4):774--786, 2012.

\bibitem{Ruelle82}
David Ruelle.
\newblock Characteristic exponents and invariant manifolds in {H}ilbert space.
\newblock {\em Ann. of Math. (2)}, 115(2):243--290, 1982.

\bibitem{Ru90}
David Ruelle.
\newblock An extension of the theory of {F}redholm determinants.
\newblock {\em Inst. Hautes \'Etudes Sci. Publ. Math.}, 72:175--193 (1991),
  1990.

\bibitem{Ru97}
David Ruelle.
\newblock Differentiation of {SRB} states.
\newblock {\em Comm. Math. Phys.}, 187(1):227--241, 1997.

\bibitem{Ru98}
David Ruelle.
\newblock Nonequilibrium statistical mechanics near equilibrium: computing
  higher-order terms.
\newblock {\em Nonlinearity}, 11(1):5--18, 1998.

\bibitem{Ru97erratum}
David Ruelle.
\newblock Correction and complements: ``{D}ifferentiation of {SRB} states''
  [{C}omm. {M}ath. {P}hys. {\bf 187} (1997), no. 1, 227--241].
\newblock {\em Comm. Math. Phys.}, 234(1):185--190, 2003.

\bibitem{Rugh08}
Hans~Henrik Rugh.
\newblock On the dimensions of conformal repellers. {R}andomness and parameter
  dependency.
\newblock {\em Ann. of Math. (2)}, 168(3):695--748, 2008.

\bibitem{Rugh10}
Hans~Henrik Rugh.
\newblock Cones and gauges in complex spaces: spectral gaps and complex
  {P}erron-{F}robenius theory.
\newblock {\em Ann. of Math. (2)}, 171(3):1707--1752, 2010.

\bibitem{Sedro}
Julien Sedro.
\newblock A regularity result for fixed points, with applications to linear
  response.
\newblock {\em Nonlinearity}, 31(4):1417, 2018.

\bibitem{SelleyTanzi20}
Fanni~M. {Selley} and Matteo {Tanzi}.
\newblock {Linear Response for a Family of Self-Consistent Transfer Operators}.
\newblock {\em arXiv e-prints}, January 2020.
\newblock \url{https://arxiv.org/abs/2001.01317}.

\end{thebibliography}
\bibliographystyle{plain}
\end{small}
\end{document}